\documentclass[twoside,leqno]{article}
\usepackage{latexsym}
\newtheorem{thm}{\sc Theorem}[section]
\newtheorem{dft}{\sc Definition}[section]

\newtheorem{cor}[thm]{\sc Corollary}
\newtheorem{pro}[thm]{\sc Proposition}

\newtheorem{lem}[thm]{\sc Lemma}
\newtheorem{cnj}[thm]{\sc Conjecture}

\def\rem{{\sc Remark}\ }
\def\rems{{\sc Remarks}\ }

\def\Rom #1{\uppercase\expandafter{\romannumeral #1}}
\def\dem{{\sc Proof}}

\def\thebibliography#1{\section*{\centerline{References}}\list
 {[\arabic{enumi}]}{\settowidth\labelwidth{[#1]}\leftmargin\labelwidth
 \advance\leftmargin\labelsep
 \usecounter{enumi}}
 \def\newblock{\hskip .11em plus .33em minus -.07em}
 \sloppy
 \sfcode`\.=1000\relax}

\title{\Large\bf Conditional Bounds for Prime Gaps with Applications}
\author{\sc Jacques Grah}
\date{}

\begin{document}
\maketitle

\begin{abstract}
We posit that $d_n^2 < 2p_{n+1}$ holds for all $n\geq 1$, where $p_n$ represents
the $n$th prime and $d_n$ stands for the $n$th prime gap i.e. $d_n := p_{n+1} - p_n$. Then,
the presence of a prime between successive perfect squares, as well as the validity of
$\Delta_n := \sqrt{p_{n+1}} - \sqrt{p_n} < 1$ are derived. Next, $\pi(x)$ being the number
of primes  $p$ up to $x$, we deduce $\pi(n^2-n) < \pi(n^2) < \pi(n^2+n)$ $(n\geq 2)$.
In addition, a proof of $\pi((n+1)^k) - \pi(n^k) \geq \pi(2^k)$ \ $(k\geq 2, n\geq 1)$ is worked out.
The vanishing nature of $\Delta_n$ as $n$ goes to infinity is set, and used afterwards to achieve both
$\displaystyle{\lim_{n\rightarrow\infty}d_n/\sqrt{p_n} = 0}$
and the twin prime conjecture. Also, question about the estimate
$p_n < 2j_n^2 \ (n\geq 6)$, where $j_n$ counts the twin prime pairs up to $p_n$, is raised.
Finally, we put forward the conjecture that any rational number $r$ $(0\leq r \leq 1)$ represents an
accumulation point of the sequence $\left(\{\sqrt{p_n}\}\right)_{n\geq 1}$, where $\{x\}$
acts for the fractional part of $x$.
\end{abstract}

\section{Introduction}
\setcounter{equation}{0}
\renewcommand{\theequation}{\thesection.\arabic{equation}}
Let $p_n$ be the $n$th prime number and $d_n$ stand for the $n$th prime gap defined by
$d_n = p_{n+1} - p_n$, where $d_1 =1$, $d_2 = 2$, etc.
The study of the rate of growth of $d_n$ is intimately related to the distribution of primes in
small intervals. Many results comparing $d_n$ to $p_n^{\theta}$ ($\theta$ a real number)
have been published. A result by Tschebycheff well known as Bertrand's postulate set $\theta$ to $1$:
$d_n < p_n$. Later on, Hoheisel was the first to prove the existence of a value of $\theta < 1$
such that $d_n < p_n^{\theta}$ for all $n$ sufficiently large, see [5].
 Subsequent values of $\theta$ getting closer and closer to $1/2$ have been established ever since.
See [2] for a comprehensive description of the race. Regarding the record to date, Baker, Harman and Pintz [3]
came up with $\theta = 0.525$.
Not surprisingly, the conjecture along these lines requires a constant $c$ such that

\vskip 7 pt \noindent
\begin{eqnarray}
d_n < c(p_n)^{1/2} \quad (n\geq 1).
\end{eqnarray}

There are a bunch of defiant questions regarding primes. One of them related to $d_n$
has to do with the difference of the square roots of consecutive primes. Actually, in [1] Andrica
conjectured that $\displaystyle{\Delta_n := \sqrt{p_{n+1}} - \sqrt{p_{n}} < 1}$ holds for all $n$.
The infinitude of primes of shape $n^2+1$ is another example, as is the decreasing nature of the sequence
$\displaystyle{(p_n^{1/n})_{n\geq 1}}$. Other tenacious
concerns revolve about primes and perfect squares.
For instance, the existence of a prime between the squares of successive integers $n$ and $n+1$,
so-called Legendre's conjecture, is  still waiting for an answer. This is surely implied by
Andrica's guess that is equivalent to $d_n < 2\sqrt{p_n} + 1$.  At the same time, when $n$ and $n+1$ are
replaced in order by
successive primes $p_n$ and $p_{n+1}$ in Legendre's conjecture, then emerges an assertion by Brocard that
$\displaystyle{\pi(p_{n+1}^2) - \pi(p_n^2) \geq 4}$  ($n>1$).
The double inequality $\pi(n^2 - n) < \pi(n^2) < \pi(n^2 + n) \ (n > 1)$, stated in 1882 by
Oppermann is yet unproved.  These and their interconnections are exposed in [6].

\vskip 10pt
Consider $n\geq 2$ such that  $d_n^2 < 2p_{n+1}$. Then, $(p_n - d_n)p_{n+1} = p_n^2 - d_n^2$
represents the largest odd multiple of $p_{n+1}$ up to $p_n^2$.
In fact, $(p_n - (d_n-2))p_{n+1} = p_n^2 - (d_n^2 - 2p_{n+1})$
exceeds $p_n^2$. Therefore, the distance $d_n^2$ from $p_n^2$ down to this odd multiple of $p_{n+1}$
satisfies either $p_{n+1} < d_n^2 < 2p_{n+1}$ or is in $[1,p_{n+1}[$.
Otherwise, there exist at least three values of $n$: $4$, $9$ and $30$ for which $d_n^2 > p_{n+1}$.
Presuming that for all $n\geq 1$, $d^2_n$ is less than $2p_{n+1}$,
arises $d_n < 2(p_n)^{1/2} \ (n\geq 1)$. So, we submit

\vskip 5pt\noindent{\cnj\ }{\sl $d_n \ <  \ \sqrt{2}\left(p_{n+1}\right)^{1/2} \quad (n\geq 1)$. }

\noindent
\underline{\hspace{5cm}}

\begin{small}

\noindent
2010 Mathematics Subject Classification: 11N05, 11B05, 11N32 \\
Key Words: Prime gaps, Andrica's, Brocard's, Legendre's and Oppermann's conjectures, Primes between conse-cutive
powers, Twin prime conjecture, Prime-generating polynomials, accumulation points.
\end{small}

\vskip 10pt\noindent The conjecture has equivalent forms and consequences to which the rest of the section is devoted.

\vskip 10pt \noindent {\cor\ }{\sl $\displaystyle{d_n^2 \leq 2(p_{n+1}-1) < 4(p_n - 1/2)\quad (n\geq 1).}$
 }

\vskip 10pt \noindent {\cor\ }
{\sl The open interval $]p_n^2 - d_n^2 , p_n^2[$ contains no multiple of $p_{n+1}$,
if and only if,
$\displaystyle{d_n < \sqrt{p_{n+1}}}$. Equivalently, $(p_n-d_n+1)p_{n+1} = p_n^2 - d_n^2 + p_{n+1}$
is the largest multiple of $p_{n+1}$ less than $p_n^2$, if and only if,
$\displaystyle{(p_{n+1})^{1/2} < d_n < \sqrt 2(p_{n+1})^{1/2}}$.}

\vskip 10pt \noindent {\cor\ }{\sl $\displaystyle{d_n < 2p_n^{1/2} \quad (n\geq 1)}$. There is a prime in the interval $(x , x +2\sqrt x) \ (x\geq 1)$.
 }

\vskip 10pt Here are other four equivalent expressions to $\displaystyle{d_n^2 < 2p_{n+1}}$.

\vskip 6pt\noindent {\thm\ }{ \sl For $n\geq 1$, either of the next statements is equivalent to the conjecture.
\begin{enumerate}
\item[1.] \quad $p_n > d_n(d_n/2 - 1)$, that is, $\displaystyle{\frac{p_n}2 > \sum_{i=1}^{d_n/2-1}i}$.
\item[2.] \quad $\displaystyle{\sum_{i=1}^{d_n}i \ - \ \sum_{i=1}^nd_i \ < \ \frac{d_n}2 \ + \ 2}.$
\item[3.] \quad $\displaystyle{\left(p_{n+1}\right)^{1/2}\ < \ \left(p_n + 1/2\right)^{1/2}
+ \sqrt{2}/2.}$
\item[4.] \quad $\displaystyle{\left(1 + 1/(2p_n)\right)p_{n+1} \ < \ \left(\sqrt{p_n} +
(\sqrt{2}/2)\sqrt{p_{n+1}/p_n}\right)^2.}$
\end{enumerate}
}

\vskip 10pt\noindent\dem\quad $d_n^2 < 2p_{n+1} = 2(p_n + d_n)$, i.e. $p_n > d_n(d_n/2 -1 )$, say, the first estimate.
To establish the second inequality, we rephrase the conjecture
using the identities:
$\displaystyle{p_{n+1} = 2 + \sum_{i=1}^nd_i}$ and $d_n^2 = 2(1 + 2 + \cdots + d_n) - d_n$.
As to the third one, replacing $d_n$ by $p_{n+1}-p_n$ in
$d_n < \sqrt{2}\left(p_{n+1}\right)^{1/2}$, gives
$\displaystyle{p_{n+1} - \sqrt{2}(p_{n+1})^{1/2} < p_n }$ or $\displaystyle{\left(\sqrt{p_{n+1}} -
\sqrt{2}/2\right)^2  < p_n + 1/2}$; hence the expected estimate. Clearly,
$\displaystyle{p_{n+1} < p_n + }$ $\displaystyle{\sqrt{2}(p_{n+1})^{1/2} = \left(\sqrt{p_n} +
(\sqrt{2}/2)\sqrt{p_{n+1}/p_n}\right)^2 - \frac{p_{n+1}}{2p_n}}$,
 which is the last estimate.
\hspace{\fill}$\Box$

\vskip 10pt Consider the identity $\displaystyle{
 p_{n+1} \ = \ k_nd_n \ + \ r_n, \ \mbox{where} \  1 \leq r_n  < d_n }$. Thus,
 the integers $k_n$ and $r_n$ are defined by the Division Algorithm of $p_{n+1}$ by the $n$th
prime gap $d_n$. Surely $r_n$ is odd and coprime to $k_n$ and to $d_n$.
We give anew other equivalent formulations to the estimate $d_n^2 < 2p_{n+1}$.

\vskip 10pt\noindent {\thm\ }
{\sl For $n\geq 2$, either of the following propositions is equivalent to $d_n^2 < 2p_{n+1}$.
\begin{enumerate}
\item[1.] There are $d_n/2$ odd integers not exceeding $p_n$ whose product with $p_{n+1}$ is greater than $p_n^2$.
That is, the smallest odd multiple of $p_{n+1}$ greater than $p_n^2$ is $(p_n - (d_n-2))p_{n+1}$.
\item[2.] $\displaystyle{(d_n - 2)p_{n+1}}$ is the largest even multiple of $p_{n+1}$ up to $d_np_n$.
\item[3.] If $p_{n+1}  =  k_nd_n  +  r_n$ with $1 \leq r_n  < d_n$, then $k_n$ is greater than or equal to $d_n/2$.
\end{enumerate}
}

\vskip 5pt\noindent\dem\quad Let $m\geq 0$ such that the estimate $(p_n-2m)p_{n+1} > p_n^2$ holds. Then,
 $\displaystyle{m < \frac{d_n}2 - \frac{d_n^2}{2p_{n+1}} < \frac{d_n}2}.$
If $2m = d_n - 2$, then $(p_n - 2m)p_{n+1} = p_n^2 - d_n^2 + 2p_{n+1}$ which is greater than $p_n^2$, if and
only if, $2p_{n+1} - d_n^2 > 0$. While for $2m = d_n$ we have $(p_n - d_n)p_{n+1} = p_n^2 - d_n^2 < p_n^2$.
Thus, $m = 0, 1, \cdots , \frac{d_n}2 -1$. For the second statement,
$(d_n - 2)p_{n+1} < d_np_n$ is equivalent to $d_n^2 < 2p_{n+1}$. As to the last proposition,
let $p_{n+1} =  k_nd_n + r_n$. Then the estimate $p_{n+1} =  k_nd_n + r_n > d_n^2/2$ holds,
if and only if, $k_n + r_n/d_n > d_n/2$, where $r_n/d_n < 1$. That is, if and only if, $k_n\geq d_n/2$.
\hspace{\fill}$\Box$

\vskip 10pt Lastly, there is a well known result in the Theory of quadratic residue to the effect
that the set $E_n = \{1^2, 2^2, 3^2, \cdots , ((p_n-1)/2)^2\}$ where $p_n$ is an odd prime,
is a reduced residue system modulo $p_n$. It is instructive to note the equivalence between
Conjecture 1.1 and the fact that $d_1^2\in E_2$ and $(d_n^2/4)^2$ must belong to $E_{n+1}$
for every $n\geq 2$.

\section{A Few Consequences of the Conjecture}
\setcounter{equation}{0}
\renewcommand{\theequation}{\thesection.\arabic{equation}}
$d_n < \sqrt 2(p_{n+1})^{1/2}$ is used
to bound both the difference of inverses and ratios of successive primes.

\vskip 10pt\noindent {\pro\ }
{\sl $\displaystyle{\frac 1{p_n} - \frac 1{p_{n+1}} \leq \frac 16} \ (n\geq 1)$,
where equality is reached if and only if $n = 1$.}

\vskip 10pt\noindent\dem\quad
$\displaystyle{\frac 1{p_n} - \frac 1{p_{n+1}} = \frac{d_n}{p_np_{n+1}}}$
and $\displaystyle{\frac{d_1}{p_1p_2} = \frac 16}.$ If $p_n$ and $p_{n+1}$ are twin primes,
then $\displaystyle{\frac{d_n}{p_np_{n+1}}\leq \frac 2{15} < \frac 16}.$
Suppose $d_n=2m \geq 4$. Using $2m^2 < p_{n+1}$  (and $2m(m-1) < p_n $),
one has $4m^3(m-1) < p_np_{n+1}$, or equivalently
$\displaystyle{\frac{d_n}{p_np_{n+1}} = \frac{2m}{p_np_{n+1}} < \frac 1{2m^2(m-1)} \leq \frac 18 < \frac 16}.$
That completes the proof. \hspace{\fill}$\Box$

\vskip 10pt \noindent {\pro\ }
{\sl $\displaystyle{\frac{p_{n+1}}{p_n}} \leq 5/3$ \ $(n\geq 1)$, with equality when $n=2$. }

\vskip 10pt\noindent\dem\quad $\displaystyle{\frac{p_{n+1}}{p_n} = \frac{p_n + d_n}{p_n} = 1 +
\frac{d_n}{p_n}.}$ If $n = 1$, $\displaystyle{\frac{p_2}{p_1} = \frac 32 < 5/3. }$ For $d_n = 2$,
$\displaystyle{\frac{p_{n+1}}{p_n} = }$  $\displaystyle{ 1 + \frac 2{p_n} \leq 5/3.}$
Also, $\displaystyle{\frac{p_5}{p_4} = \frac{11}7 < 5/3. }$ Suppose $d_n = 2m\geq 4$ and $n\geq 6$.
Then, using $d_n < \sqrt{2p_{n+1}}$ leads to $\displaystyle{\frac{p_{n+1}}{p_n} - 1 < }$
\\
$\displaystyle{\frac{\sqrt 2}{\sqrt{p_n}}\sqrt{\frac{p_{n+1}}{p_n}} =
\frac{\sqrt 2}{\sqrt{p_n}}\sqrt{1 + \frac{2m}{p_n}} \leq \sqrt{\frac 2{13}}\sqrt{1 + \frac{2m}{p_n}}}$,
since $n\geq 6$. Next, $p_{n+1} < 2p_n$ along with $2m^2 < p_{n+1}$, yield
 $\displaystyle{\frac 1{2p_n} <  \frac 1{p_{n+1}} < }$ $\displaystyle{\frac 1{2m^2}}$, that is
$\displaystyle{\frac{2m}{p_n} <  \frac{4m}{p_{n+1}} < \frac 2{m}}.$ Substituting this into
the formula above, since $m\geq 2$ we get
$\displaystyle{\frac{p_{n+1}}{p_n} < 1 + \sqrt{\frac 2{13}}\sqrt{1 + \frac 2m} < 1 + \frac 2{\sqrt{13}} < 5/3,}$
and the proof is complete. \hspace{\fill}$\Box$

\vskip 10pt\noindent{\rems\ }Let $n\geq 2$. It comes from the proof that
$\displaystyle{\frac{p_{n+1}}{p_n}}$ is equal to $\displaystyle{1 + \frac 2{p_n}}$ if $d_n = 2$, and less than
$\displaystyle{1 + \frac 2{\sqrt{p_n}}}$ otherwise, provided $n\geq 5$. Combining both cases, we have
$\displaystyle{\frac{p_{n+1}}{p_n} <  1 + \frac 2{\sqrt{p_n}} }$ $(n\geq 2)$. In particular,
$\displaystyle{\frac{p_{n+1}}{p_n} < 3/2 }$ $(n\geq 7)$ and we verify that $p_{n+1}/p_n < 3/2$ holds for $n=5$ and $ n = 6$.

\vskip 10pt A consequence of the proposition is
\begin{eqnarray}
p_n \ \geq \ d_{n+1} \quad (n\geq 1), \ \mbox{ where equality holds for } n = 1.
\end{eqnarray}

\vskip 7pt\noindent
This is really true for $n = 1, 2$ and $3$. If $p_n$ is less than $d_{n+1}$ then
$p_n^2 < d_{n+1}^2 < 2p_{n+2} < 4p_{n+1}$, i.e. $p_n < 4\frac{p_{n+1}}{p_n} < 20/3 < 7 = p_4$,
according to the proposition. But, $p_1$, $p_2$ and $p_3$ satisfy the opposite inequality:
$p_n\geq d_{n+1}$.

\vskip 10pt Another application of $d_n^2 < 2p_{n+1}$ is Legendre's conjecture.

\vskip 10pt \noindent {\thm\ }{\sl There is a prime between two successive perfect squares. }

\vskip 10pt\noindent\dem\quad Assume $N\geq 4$ such that $[N^2 , (N+1)^2]$
is devoid of any prime. Then, there exists an integer $k \geq 5$ such that
$\displaystyle{d_k = p_{k+1} - p_k > (N+1)^2 - }$ $\displaystyle{ N^2 + 1 = 2(N+1).}$
Next, let $\mu_{k+1}$ represents the fractional part of $\sqrt{p_{k+1}}$. Then there exists some $r\geq 1$
such that $\sqrt{p_{k+1}} = (N+r)+\mu_{k+1} \geq N+1+\mu_{k+1}$. Thus,
$\displaystyle{d_k > 2(N+1) = \sqrt{2}(p_{k+1})^{1/2} + (2 - \sqrt{2})(p_{k+1})^{1/2} - 2\mu_{k+1}}$,
where, since $k \geq 5$, $(2 - \sqrt{2})(p_{k+1})^{1/2} - 2\mu_{k+1}$ is positive. Then, follows the
contradiction $d_k > \sqrt{2}(p_{k+1})^{1/2}.$
 \hspace{\fill}$\Box$

\vskip 10pt A further consequence of $\displaystyle{d_n<\sqrt 2(p_{n+1})^{1/2}}$ is
$\displaystyle{\lim_{n\rightarrow\infty}d_n/p_n = 0}$. Of course, a
 proof of this limit based on the Prime Number Theorem already exists.

\vskip 10pt\noindent{\cor\ }{\sl $\displaystyle{\lim_{n\rightarrow\infty}d_n/p_n = 0}$.}

\vskip 10pt\noindent\dem\quad $\displaystyle{d_n/p_n < \sqrt{2p_{n+1}}/p_n < 2/\sqrt{p_n}}$
yielding the claimed limit.
\hspace{\fill}$\Box$

\vskip 10pt Moving on to Andrica's question, the postulate implies $d_n < \sqrt{2p_{n+1}} < 2\sqrt{p_n}$,
which leads to $\displaystyle{\sqrt{p_{n+1}} - \sqrt{p_n} = \frac{d_n}{\sqrt{p_{n+1}}+\sqrt{p_n}}
< \frac{d_n}{2\sqrt{p_n}} < 1}$, establishing the inequality. Similarly,
the left-hand side in the fourth item of Theorem 1.5 is greater than $p_{n+1}$, while the right-hand side,
is less than $(\sqrt{p_n} + 1)^2$, that is $\displaystyle{p_{n+1} < (\sqrt{p_n} + 1)^2}$;
confirming again Andrica's guest. This proves also the equivalence between $d_n < 2\sqrt{p_n} + 1$
and $\sqrt{p_{n+1}} - \sqrt{p_n} < 1$ stated in the introduction. The third statement of the same theorem
can be used to justify Andrica's inequality. Indeed, its right-hand side is surely less than
$\sqrt{p_n} + 1/4 + \sqrt{2}/2$, i.e. $\sqrt{p_{n+1}} < \sqrt{p_n} + 0.957\cdots$. Lastly,
$\displaystyle{\frac{p_{n+1}}{p_n} < 1 + \frac 2{\sqrt{p_n}}}$ from remarks on the preceding page induces
$\displaystyle{\sqrt{\frac{p_{n+1}}{p_n}} < 1 + \frac 1{\sqrt{p_n}}}$, hence $\sqrt{p_{n+1}} - \sqrt{p_n} < 1$.
A tentatively refined bound in Andrica's statement will be set in the next subsection.

\subsection{$d_n < 2\sqrt{p_n}\quad $ Revisited}
Replacing $p_{n+1}$ by $p_n$ into $d_n/2 \leq (p_{n+1}-1)/4$
 yields a relation that holds no more for every $n\geq 1$.
At least three values of $n$, \ $n=1, 2, 4$ for which $2d_n \leq p_n - 1$ is false.
To better substitute $p_{n+1}$ by $p_n$ into $d_n < \sqrt{2p_{n+1}}$, let $n\geq 1$
satisfy $2d_n - r_n = p_n$ for some $r_n$ \ $(1\leq r_n < d_n)$. Then, $2d_n + (d_n-r_n) = p_{n+1}$,
so that in accordance with Theorem 1.6(3),
$k_n = 2\geq d_n/2$ i.e. $d_n\leq 4$ which infers $n=2$ or $n=4$, values for which
$d_n = \frac{p_n+1}2$. Likewise, $d_n + 1 = p_n$ gives again $d_n \leq 4$; so $n=1$
or $n = 2$. The situation of $2d_n + 1 = p_n$ leads to $d_n \leq 6$; thereby $d_n = (p_n - 1)/2$ for $n=3$.
Therefore,
\begin{eqnarray}
d_n\quad \leq\quad \frac{p_n - 1}2 \quad (n\neq 1, 2 \mbox{ and } 4).
\end{eqnarray}

\vskip 3pt\noindent
We guess that with this, a proof of preceding conjecture could be adapted to

\vskip 10pt\noindent {\cnj\ }{\sl $\displaystyle{d_n < \sqrt 2\left(p_n\right)^{1/2}}$ \ holds for all
values but $n = 4$. }

\vskip 10pt\noindent This conjecture, just like the previous, has equivalent expressions and consequences.

\vskip 10pt \noindent {\cor\ }{\sl $d_n^2 \leq 2(p_n - 1) \quad (n\neq 4)$. }

\vskip 10pt \noindent {\cor\ }{\sl The interval $(x , x+\sqrt 2x^{1/2})$ contains a prime number,
provided $x$ is not in the range $7 \leq x \leq 7.2 041 684 766$. }

\vskip 10pt The next propositions are equivalent to Conjecture 2.5.

\vskip 7pt\noindent {\thm\ }
{\sl Let $n\geq 2$ not equal to $4$. Then the next assertions are equivalent to
$\displaystyle{d_n < \sqrt{2}\left(p_n\right)^{1/2}}.$
\begin{enumerate}
\item[1.] $(p_{n+1} + d_n)p_n$ is the largest odd multiple of $p_n$ up to $p_{n+1}^2$.
\item[2.] There are $\displaystyle{d_n}$ odd multiples of $p_n$ in the open interval
$\displaystyle{]p_n^2 , p_{n+1}^2[}$.
\item[3.] The sum of positive integers up to $d_n$ is less than the arithmetic mean of the bounding primes $p_n$
and $p_{n+1}$. That is to say, $\displaystyle{\frac{d_n}2 + \sum_{i=1}^{d_n-1}i \ < \ p_n}.$
\end{enumerate}
}

\vskip 7pt\noindent\dem\quad For an integer $m \geq 0$ satisfying  $(2m + 1)p_n < p_{n+1}^2$, one has \\
$\displaystyle{m <
\frac{p_{n+1}^2 - p_n}{2p_n} = \frac{(p_n + d_n)^2 - p_n}{2p_n} = \frac{p_{n+1} + d_n - 1}2 + \frac{d_n^2}{2p_n}.}$
Thus, $m = 0, 1, \cdots , $ $(p_{n+1} + d_n - 1)/2$, if and only if,  $d_n^2 < 2p_n$.  For the second statement,
the number of odd multiples of $p_n$ in $]p_n^2,p_{n+1}^2[$ is the integral part of
$\displaystyle{\frac{p_{n+1}^2 - p_n^2}{2p_n} = \frac{d_n(2p_n + d_n)}{2p_n} =}$
$\displaystyle{\frac{2p_nd_n + d_n^2}{2p_n}}$, that is
$\displaystyle{\frac{p_{n+1}^2 - p_n^2}{2p_n}  \ = \ d_n + \frac{d_n^2}{2p_n}}$, whose integral part
equals $d_n$, if and only if, $d_n^2 < 2p_n$.
With respect to the last premise, the inequality $\displaystyle{d_n^2 < 2p_n}$ is
equivalent to $\displaystyle{d_n^2 + d_n < 2p_n + d_n = p_n + p_{n+1}},$ that is,
$\displaystyle{d_n(d_n + 1) = 2\sum_{i=1}^{d_n}i < p_n + p_{n+1}}$.
Rewriting this last inequality, yields
 $\displaystyle{2d_n + 2\sum_{i=1}^{d_n-1}i < p_n + p_{n+1} = 2p_n + d_n}$. That is,
 $\displaystyle{\frac{d_n}2 + \sum_{i=1}^{d_n-1}i \ < \ p_n}.$
\hspace{\fill}$\Box$

\vskip 10pt \noindent  For $n\neq 4$,
$\displaystyle{\sqrt{p_{n+1}} - \sqrt{p_{n}} < \frac{d_n}{2\sqrt{p_n}} < \frac{\sqrt2(p_n)^{1/2}}{2\sqrt{p_n}} =
\frac {\sqrt{2}}2}$. Also, $\sqrt{11} - \sqrt{7} = 0.67087... $ $< 7/10 < \frac{\sqrt 2}2$.
Thus, one improves Andrica's conjecture.

\vskip 10pt \noindent {\thm\ }
{\sl $\displaystyle{\sqrt{p_{n+1}} - \sqrt{p_{n}} \leq \sqrt{11} - \sqrt{7} < 7/10 < \sqrt{2}/2 \quad (n\geq 1).}$
}

\vskip 5pt \noindent {\cor\ }
{\sl Any sequence of $k$ consecutive composite integers contains at most one perfect square.
This holds particularly for $(k+1)!+2$, $(k+1)!+3$, $\cdots$, $(k+1)!+k+1$.
}

\vskip 10pt Andrica's conjecture implies the one by Legendre. Suppose that there is no
prime between $N^2$ and $(N+1)^2$ for a certain $N\geq 2$. Then, there exists
$k$ such that $p_k < N^2 < (N+1)^2 < p_{k+1}$, implying that
$\displaystyle{\sqrt{p_{k+1}} - \sqrt{p_k} > (N+1) - N =1}$, contradicting
Andrica's statement. Evidently, $p_k < N^2 < (N+1)^2 < p_{k+1}$ denies the corollary too, as we have
a sequence of at least $2(N+1)$ consecutive composite integers with two perfect squares.

\vskip 10pt \noindent Invoke preceding theorem to narrow down the interval $(x , x + 2\sqrt{x})$
without any restriction on $x$.

\vskip 10pt \noindent {\pro\ }{\sl The interval $(x , x + 8/5\sqrt{x})$ \ $(x\geq 1)$ contains a prime, i.e.
$d_n < 8/5\sqrt{p_n} \ \ (n\geq 1)$.
}

\vskip 10pt \noindent\dem\quad The proposition is true for $n = 1, 2$ and $3$. Consider the function $f$ of
a real variable $x\geq 1$, by
$\displaystyle{f(x) = (x + 8/5\sqrt{x})^{1/2} - \sqrt{x} = \frac{8/5}{(1 + \frac 8{5\sqrt{x}})^{1/2} + 1} < 4/5}.$
Then, $\displaystyle{f'(x) = \frac{\sqrt x + 4/5 - (x + 8/5\sqrt x)^{1/2}}{2\sqrt{x(x + 8/5\sqrt x)}} = }$ \\
$\displaystyle{\frac 8{25(\sqrt x + 4/5 + (x + 8/5\sqrt x)^{1/2})\sqrt{x(x + 8/5\sqrt x)}}}$
represents the derivative of $f$ at $x$.
It is positive on $]0 , \infty[$; hence
$f$ is increasing and bounded by $4/5$. In addition, $f(4) = 0.68328\cdots > \sqrt{11} - \sqrt 7 = 0.6708\cdots \geq
\sqrt{p_{n+1}} - \sqrt{p_{n}} \ (n\geq 1)$. If for $x\geq 4$ there is no prime between $x$ and $x + 8/5\sqrt x$, then there exits
an integer $k$ such that $\sqrt{p_{k+1}} - \sqrt{p_k} > \sqrt{x + 8/5\sqrt{x}} - \sqrt{x} = f(x) \geq f(4) > \sqrt{11} - \sqrt 7$.
The inequality $\sqrt{p_{k+1}} - \sqrt{p_k} > \sqrt{11} - \sqrt 7$ gives the desired contradiction. The proof is then complete.
\hspace{\fill}$\Box$

\subsection{Bound for $d_n$ Versus Bound for $\sqrt{p_{n+1}} - \sqrt{p_n}$}
There is a relation between bounds for $d_n$ and bounds
for the difference $\sqrt{p_{n+1}} - \sqrt{p_{n}}$.

\vskip 10pt \noindent {\pro\ }
{\sl Let $a > 0$ stand for a real number and $n_0\geq 1$ such that
$\displaystyle{d_n < a\left(p_n\right)^{1/2} \ (n\ge n_0)}$.  Then, for every integer $n\geq n_0$,
$\displaystyle{\sqrt{p_{n+1}} - \sqrt{p_{n}} < a/2}$
and $\displaystyle{\sqrt{p_{n+1}} + \sqrt{p_{n}} < 2\left(p_n\right)^{1/2} + \frac a2}$. }

\vskip 10pt \noindent\dem\quad The first estimate comes from $\displaystyle{\sqrt{p_{n+1}} -
\sqrt{p_{n}} = \frac{d_n}{\sqrt{p_{n+1}} + \sqrt{p_{n}}} < \frac{d_n}{2\sqrt{p_n}}}$. As to the second,
 it proceeds via $\displaystyle{\sqrt{p_{n+1}} - \sqrt{p_{n}} < a/2}$ to which we add
$2\sqrt{p_n}$.  \hspace{\fill}$\Box$

\vskip 15pt  The converse of the proposition is false. For $a = \sqrt 2$,
$\sqrt{11} - \sqrt 7 = 0.6708\cdots < \frac{\sqrt 2}2$ but $d_4 > \sqrt 2\times \sqrt 7 = 3.7416\cdots$.
However, we do have the following.

\vskip 10pt \noindent {\pro\ }
{\sl Let $a > 0$ be a real number and $n\geq 1$ an integer such that
$\displaystyle{\sqrt{p_{n+1}} - \sqrt{p_{n}} < a}$. Then,
 $\displaystyle{d_n < 2a\left(p_n\right)^{1/2} + a^2}$ and
$\displaystyle{\sqrt{p_{n+1}} + \sqrt{p_{n}} < 2\left(p_n\right)^{1/2} + a}$.}

\vskip 10pt\noindent\dem\quad
$\displaystyle{d_n < a\left(\sqrt{p_{n+1}} + \sqrt{p_{n}}\right) < a\sqrt{p_{n}} + a\left(a +
\sqrt{p_{n}}\right) = a^2 + 2a\sqrt{p_{n}}}$ results from calling on the hypothesis twice;
and thereby the desired estimates.
 \hspace{\fill}$\Box$

\vskip 10pt\noindent For example, $a = \sqrt 2/2$ leads to
$\displaystyle{d_n < \sqrt 2(p_n)^{1/2} + 1/2}$ for every $n\geq 1$.

\vskip 10pt Now, we firstly establish and next improve upon a result: $p_{n+2} < p_n + p_{n+1}, (n\geq 2)$
by Ishikawa (see [6], p. 185). For the proof use the identity
$p_n + p_{n+1} = p_n + p_{n+2} - d_{n+1} = p_{n+2} + (p_n - d_{n+1})$ and (2.1).

\vskip 10pt \noindent {\thm\ }{\sl $\displaystyle{p_{n+2} < (p_n + 1) + 2\sqrt 2(p_n)^{1/2}\quad (n\geq 1)}$.
That is, there are at least two primes between $p_n$ and
$\displaystyle{p_n + 1 + 2\sqrt{2p_n} = \left(\sqrt{p_n} + \sqrt 2\right)^2 - 1}$.}

\vskip 10pt\noindent\dem\quad Obviously, $n = 1, 2, 3$ and $4$ verify the theorem.
For $n\geq 5$, one has $\displaystyle{p_{n+1} - p_n < \sqrt 2(p_n)^{1/2}}$ to which we
add $\displaystyle{p_{n+2} - p_{n+1} < \sqrt 2(p_{n+1})^{1/2}}$. Calling in proposition 2.12
with $a = \sqrt 2$, gives
\\
$\displaystyle{
p_{n+2} <  p_n  + \sqrt 2\left(\sqrt{p_n} + \sqrt{p_{n+1}}\right) < p_n + }$
 $\displaystyle{2\sqrt 2(p_n)^{1/2} + 1}$, hence
the inequality.   \hspace{\fill}$\Box$

\vskip 10pt \noindent {\rem\ }
Empirically, $p_n + 2\sqrt 2(p_n)^{1/2} = \left(\sqrt{p_n} + \sqrt 2\right)^2 - 2 > p_{n+2}$.
To amend the theorem, we proceed with
 $\displaystyle{p_{n+2} - p_{n+1} < \sqrt 2(p_{n+1})^{1/2} = \sqrt{2p_n}\left(1 + \frac{d_n}{p_n}\right)^{1/2}}$
that is bounded above by $\displaystyle{\sqrt{2p_n}\left(1 + \frac{d_n}{2p_n}\right)}$. This yields
$\displaystyle{p_{n+2} - p_n < \sqrt{2p_n}\left(1 + \frac{d_n}{2p_n}\right) +
d_n < 2\sqrt{2p_n} + \frac{d_n}{\sqrt{2p_n}}}$. We come through, for $n\geq 5$, \\
$\displaystyle{p_{n+2} < p_n + 2\sqrt{2p_n} + \frac{d_n}{\sqrt{2p_n}} = \left(\sqrt{p_n} +
\sqrt 2\right)^2 - 2 + \frac{d_n}{\sqrt{2p_n}}}$, where $\displaystyle{\frac{d_n}{\sqrt{2p_n}} < 1}$.

\section{The Effectiveness of the Bound}
\setcounter{equation}{0}
\renewcommand{\theequation}{\thesection.\arabic{equation}}
Define $\Delta_n := \sqrt{p_{n+1}} - \sqrt{p_n}$. Then,
$d_n = \Delta_n(\sqrt{p_{n+1}} + \sqrt{p_n})$ and $\Delta_n < 7/10 < \sqrt 2/2$. Thus,
\begin{eqnarray}
2\sqrt{p_{n}}\Delta_n \ < \ d_n \ < \ 2\sqrt{p_{n+1}}\Delta_n  \quad (n\geq 1).
\end{eqnarray}
On the one hand, the difference between the two bounds equals $2\Delta_n^2$ that
is less than $1$, so:
\begin{eqnarray}
2\sqrt{p_{n}}\Delta_n < d_n < 2\sqrt{p_{n}}\Delta_n + 2\Delta_n^2.
\end{eqnarray}
This leads to the fact that $d_n$ and the integral part of the upper bound in (3.1) must not only share the
same parities but be equal too.

\vskip 10pt Given a real number $x$, let $\lfloor{x}\rfloor$ and $\{x\}$ stand for the integral and the fractional parts
of $x$, respectively.

\vskip 10pt \noindent {\pro\ }
{\sl $\displaystyle{ d_n\ = \ \lfloor{2\sqrt{p_{n+1}}\Delta_n}\rfloor}$. As a result, for $n > 1$,
$\displaystyle{\lfloor{2\sqrt{p_n}\Delta_n}\rfloor}$  is of odd parity that is,
$\displaystyle{d_n\ = \ \lfloor{2\sqrt{p_n}\Delta_n}\rfloor + 1}$.}

\vskip 10pt \noindent Rewriting (3.2) as $\displaystyle{(\sqrt{p_{n}}\Delta_n - 1/2) - 1/2 < \frac{d_n}2 - 1
< (\sqrt{p_{n}}\Delta_n - 1/2) - (1/2 - \Delta_n^2),}$ where $0 < 1/2 - \Delta_n^2  < 1/2$,
induces immediately

\vskip 10pt\noindent {\pro\ }
{\sl $\displaystyle{\frac{d_n}2 - 1 \ = \ \lfloor{\sqrt{p_n}\Delta_n - 1/2}\rfloor \quad (n\geq 2).}$
 }

\vskip 10pt\noindent {\cor\ }
{\sl $\{\sqrt{p_n}\Delta_n - 1/2\} < 1/2$ \ $(n\geq 2)$.
}

\vskip 10pt\noindent\dem\quad The proposition gives $d_n/2 - 1 = \sqrt{p_n}\Delta_n - 1/2 - \{\sqrt{p_n}\Delta_n - 1/2\}$,
whereas Proposition 3.1 leads to $d_n/2 - 1 = \sqrt{p_n}\Delta_n - 1/2 - \{2\sqrt{p_n}\Delta_n\}/2$. Yet, follows
$\{\sqrt{p_n}\Delta_n - 1/2\} = \{2\sqrt{p_n}\Delta_n\}/2 < 1/2$.
\hspace{\fill}$\Box$

\vskip 10pt Accordingly, for $n>1$, $\lfloor{\sqrt{p_n}\Delta_n}\rfloor$ equals
$d_n/2 - 1$, the number of odd composite integers within $[p_n , p_{n+1}]$. That is,
$d_n/2 - 1 = \sqrt{p_n}\Delta_n - \left(1/2 + \{\sqrt{p_n}\Delta_n - 1/2\}\right)$. The corollary and
the identity: $d_n = (\sqrt{p_{n+1}} + \sqrt{p_{n}})\Delta_n = \sqrt{p_{n+1}}\Delta_n + \sqrt{p_{n}}\Delta_n$,
whose leftmost term is an integer, imply $\{\sqrt{p_{n+1}}\Delta_n\} + \{\sqrt{p_n}\Delta_n\} = 1$.
Consequently,
$\lfloor{\sqrt{p_{n+1}}\Delta_n}\rfloor$ represents the number of even integers within $]p_n , p_{n+1}]$,
and all this leads to the following.

\vskip 10pt\noindent{\thm\ }{\sl $\displaystyle{
d_n \ = \ \lfloor{\sqrt{p_{n+1}}\Delta_n}\rfloor \ + \ \lfloor{\sqrt{p_n}\Delta_n}\rfloor \ + \ 1 \ (n\geq 1)}$.
That is, for every integer $n\geq 2$,
\\
$\displaystyle{  d_n \ \ = \  2\lfloor{\sqrt{p_{n+1}}\Delta_n}\rfloor }$
$\displaystyle{ \ = \ 2(\lfloor{\sqrt{p_n}\Delta_n}\rfloor \ + \ 1)}$. Furthermore,
$\displaystyle{\{\sqrt{p_{n+1}}\Delta_n\}}$ is between $0$ and $1/2$; while the quantities
$\lfloor{\sqrt{p_{n+1}}\Delta_n}\rfloor$ and $\lfloor{\sqrt{p_n}\Delta_n}\rfloor$ are of opposite parities.
  }

\vskip 10pt\noindent We draw, from the theorem, new bounds on
$\{\sqrt{p_n}\Delta_n\}$ and $\{\sqrt{p_{n+1}}\Delta_n\}$.

\vskip 10pt\noindent{\thm\ }{\sl $\displaystyle{\Delta_n^2 + 2\{\sqrt{p_n}\Delta_n\} = 2 \ (n\geq 2).}$
Equivalently, $\displaystyle{\Delta_n^2 = 2\{\sqrt{p_{n+1}}\Delta_n\}}$ $\displaystyle{(n\geq 2).}$
For $n\geq 1$,  $\displaystyle{\{\sqrt{p_{n+1}}\Delta_n\}}$ is less than $1/4$,
meanwhile $\displaystyle{\{\sqrt{p_n}\Delta_n\}}$ exceeds $3/4$ for $n\geq 2$. }

\vskip 10pt\noindent\dem\quad Since $\lfloor{\sqrt{p_n}\Delta_n}\rfloor = d_n/2 - 1$ and
$\lfloor{\sqrt{p_{n+1}}\Delta_n}\rfloor = d_n/2$,
\begin{eqnarray*}
1 & = & \lfloor{\sqrt{p_{n+1}}\Delta_n}\rfloor - \lfloor{\sqrt{p_n}\Delta_n}\rfloor \ = \
\sqrt{p_{n+1}}\Delta_n - \sqrt{p_n}\Delta_n + \{\sqrt{p_n}\Delta_n\} - \{\sqrt{p_{n+1}}\Delta_n\} \\
& = & \Delta_n^2 + \{\sqrt{p_n}\Delta_n\} - \{\sqrt{p_{n+1}}\Delta_n\} \
 = \  \Delta_n^2 + 2\{\sqrt{p_n}\Delta_n\} - 1  \\
 & = & \Delta_n^2 + 1 - 2\{\sqrt{p_{n+1}}\Delta_n\},
\end{eqnarray*}
for $\displaystyle{\{\sqrt{p_{n+1}}\Delta_n\} + \{\sqrt{p_n}\Delta_n\} = 1}$.
This establishes the two identities. The statements on fractional parts,
 derive from $\displaystyle{\Delta_n^2 < 1/2}$ along with respectively
$\Delta_n^2 = 2\{\sqrt{p_{n+1}}\Delta_n\}$ and $\Delta_n^2 = 2 - 2\{\sqrt{p_n}\Delta_n\}$.
\hspace{\fill}$\Box$

\vskip 10pt Note that $\displaystyle{\frac{d_n}2 = \sqrt{p_n}\Delta_n + \frac{\Delta_n^2}2}$ $(n\geq 2)$
implies $\displaystyle{\frac{d_n}{p_n} = \frac{2\Delta_n}{\sqrt{p_n}} + \frac{\Delta_n^2}{p_n} <
\frac 2{\sqrt{p_n}}}$, using $\Delta_n < \sqrt{2}/2$.
It yields again $\displaystyle{\lim_{n\rightarrow\infty}d_n/p_n = 0}$.

\vskip 10pt\noindent{\cor\ }{\sl $\displaystyle{\{1 + 2\sqrt{p_n}\Delta_n\}  =  2\{\sqrt{p_n}\Delta_n\} - 1 =
1 - \{2\sqrt{p_{n+1}}\Delta_n\} \ (n\geq 2)}$
}

\vskip 10pt\noindent\dem\quad Firstly, $\displaystyle{d_n = 2 + 2\sqrt{p_n}\Delta_n - 2\{\sqrt{p_n}\Delta_n\} =
1 + 2\sqrt{p_n}\Delta_n - \left( 2\{\sqrt{p_n}\Delta_n\} - 1\right)}$. Next, the previous theorem
implies $1/2 < 2\{\sqrt{p_n}\Delta_n\} - 1 < 1$, thereby showing that
$2\{\sqrt{p_n}\Delta_n\} - 1 = \{1 + 2\sqrt{p_n}\Delta_n\} = \{2\sqrt{p_n}\Delta_n\}$, and proves the first identity.
Second identity comes from $\displaystyle{1 = \{\sqrt{p_{n+1}}\Delta_n\} + \{\sqrt{p_n}\Delta_n\}}$
and $\displaystyle{\{\sqrt{p_{n+1}}\Delta_n\} < 1/4}$.
\hspace{\fill}$\Box$

\vskip 12pt\noindent We condense the different expressions of $d_n$.

\begin{eqnarray}
\frac{d_n}2 \ = \ \left\{
\begin{array}{lcl}
\sqrt{p_{n+1}}\Delta_n - \{\sqrt{p_{n+1}}\Delta_n\} & = & \sqrt{p_{n+1}}\Delta_n - 1 + \{\sqrt{p_n}\Delta_n\}    \\
\sqrt{p_n}\Delta_n + 1 -  \{\sqrt{p_n}\Delta_n\} & = &  \sqrt{p_n}\Delta_n + \{\sqrt{p_{n+1}}\Delta_n\}    \\
\sqrt{p_n}\Delta_n + \Delta_n^2  - \{\sqrt{p_{n+1}}\Delta_n\} & = &  \sqrt{p_n}\Delta_n + 1/2 - \{\sqrt{p_n}\Delta_n - 1/2\}.
\end{array}
\right.
\end{eqnarray}

\vskip 10pt Theorems 3.4 and 3.5 generate forthwith the following inequalities.

\vskip 10pt\noindent{\cor\ }{\sl $\displaystyle{\sqrt{p_n}\Delta_n < \frac{d_n}2 < \sqrt{p_{n+1}}\Delta_n
< \frac{d_n}2 + \frac{1}4 < \sqrt{p_n}\Delta_n + \frac{1}2 < \frac{\sqrt{2p_n}}2 + \frac{1}2 \ \ (n\geq 1).}$
}

\vskip 10pt\noindent\dem\quad $d_n/2 = \sqrt{p_n}\Delta_n + \Delta_n^2/2$ implies the first estimate.
The second estimate comes from $d_n/2 = \lfloor{\sqrt{p_{n+1}}\Delta_n\rfloor}$, The third stems from
$\{\sqrt{p_{n+1}}\Delta_n\} < 1/4$, whereas $\Delta_n^2 < 1/2$ along with $d_n/2 = \sqrt{p_n}\Delta_n + \Delta_n^2/2$
yield the fourth inequality. In the end, $\Delta_n < \sqrt{2}/2$ induces the rightmost estimate.
\hspace{\fill}$\Box$

\vskip 10pt\noindent $\displaystyle{\frac{d_n}2 = \sqrt{p_n}\Delta_n + \frac{\Delta_n^2}2 \ }$
gives $\displaystyle{ \ \sqrt{p_n}\Delta_n < \sqrt{p_n}\Delta_n + \frac{\Delta_n^2}2  <
\sqrt{p_n}\Delta_n + \Delta_n^2 < \sqrt{p_n}\Delta_n + \frac{3\Delta_n^2}2 < \sqrt{p_n}\Delta_n \ + }$  \\
$\displaystyle{2\Delta_n^2 < \sqrt{p_n}\Delta_n + \frac{3\Delta_n^2}2 + 1/4}$,
namely, another form of the corollary.

\vskip 12pt
In a different connection, since $\displaystyle{\{\sqrt{p_n}\Delta_n\}> 3/4}$, adding $1/4$ to $\displaystyle{\sqrt{p_n}\Delta_n}$
increases $\displaystyle{\lfloor{\sqrt{p_n}\Delta_n}\rfloor}$ by one to equal
$\displaystyle{\lfloor{\sqrt{p_{n+1}}\Delta_n}\rfloor}$.
So, can we insert $\displaystyle{\sqrt{p_n}\Delta_n + 1/4}$ between $d_n/2$ and $\displaystyle{\sqrt{p_{n+1}}\Delta_n}$ into
the estimates of the corollary? Clearly, $\displaystyle{\sqrt{p_n}\Delta_n + 1/4}$ exceeds $d_n/2$ thus the challenge is on
$\displaystyle{\sqrt{p_{n+1}}\Delta_n}$ side. Is $\displaystyle{\sqrt{p_n}\Delta_n + 1/4}$ less than
$\displaystyle{\sqrt{p_{n+1}}\Delta_n}$? This seems to be the case for finitely many values of $n$,
for instance $n = 2, 4, 6, 9, 11$ and $30$.
Bringing together this matter and the estimates (3.2), we rightly wonder whether, except for these values of $n$,
$\displaystyle{\Delta_n^2}$ is less than $1/4$. If as expected the limit of $\Delta_n$
vanishes when $n$ tends to infinity, then $\Delta_n^2 > 1/4$ holds for a finite number of integers $n$.

\vskip 10pt Another concern is next:
{\sl Are there finitely many integers $n$, like $n = 30$, satisfying
$\displaystyle{\{\sqrt{p_{n+1}}\Delta_n\} > 2/d_n}$? That is to say, does the product of the integral and fractional parts
 of $\displaystyle{\sqrt{p_{n+1}}\Delta_n}$ exceed $1$?}
An affirmative answer would lead to $\displaystyle{\lim_{n\rightarrow\infty}\left(\sqrt{p_{n+1}} - \sqrt{p_{n}}\right) = 0}$.
In fact, there would exist an integer $N$ such that $\displaystyle{\Delta_n^2 < 4/d_n \ (n > N)}$. Such an answer would
also furnish $\displaystyle{\{\sqrt{p_n}\Delta_n\} > 1 - 2/d_n}$ for every integer $n > N$, and thus
$\displaystyle{\lim_{n\rightarrow\infty}\{\sqrt{p_n}\Delta_n\} = 1}$.

\vskip 10pt\noindent Working on fractional part of $\sqrt{p_n}$ brings up issues
about some limits. Out of (3.3), the identity
$\displaystyle{\frac{d_n}2 = \sqrt{p_n}\Delta_n + \{\sqrt{p_{n+1}}\Delta_n\}}$
suggests the link:
$\displaystyle{\lim_{n\rightarrow\infty}\Delta_n = \frac 12\lim_{n\rightarrow\infty}\frac{d_n}{\sqrt{p_n}}}$.
It is conjectured that $\displaystyle{\lim_{n\rightarrow\infty}\Delta_n = 0}$. In turn,
$\displaystyle{\Delta_n = \sqrt{2}\left(\{\sqrt{p_{n+1}}\Delta_n\}\right)^{1/2}}$, thus
\begin{eqnarray}
\lim_{n\rightarrow\infty}\left(\sqrt{p_{n+1}} - \sqrt{p_{n}}\right) = \frac 12\lim_{n\rightarrow\infty}
\frac{d_n}{\sqrt{p_n}} = \sqrt 2\lim_{n\rightarrow\infty}\sqrt{\{\sqrt{p_{n+1}}\Delta_n\}}.
\end{eqnarray}
Moreover, this hints again at $\displaystyle{\lim_{n\rightarrow\infty}\{\sqrt{p_n}\Delta_n\} = 1}$,
since $\displaystyle{\{\sqrt{p_{n+1}}\Delta_n\} + \{\sqrt{p_n}\Delta_n\} = 1}$.

\vskip 10pt  Let $n<m$ such that $d_n = d_m = 2s$ for a fixed integer $s\geq 1$. Then,
$\displaystyle{\frac{d_n}2 = \sqrt{p_{n+1}}\Delta_n - \{\sqrt{p_{n+1}}\Delta_n\}}$ $\displaystyle{ = s = \frac{d_m}2 =
\sqrt{p_{m+1}}\Delta_m - \{\sqrt{p_{m+1}}\Delta_m\}}$. Since $d_n = d_m$ and $\sqrt{p_{n+1}} < \sqrt{p_{m+1}}$,
necessarily, $\Delta_n > \Delta_m$. Indeed, $d_n = \Delta_nD_n = \Delta_mD_m$
with $D_m > D_n = \frac{\Delta_m}{\Delta_n}D_m.$ As a result, $\Delta_m < \Delta_n.$
So, if $2s$ is assumed by infinitely many $n$, then there exists a sequence
$\left(s_n \right)$ such that $d_{s_n} = 2s$, for every $s_n$; and the sequence
$\displaystyle{\left(\Delta_{s_n}\right)_{s_n>1}}$
is decreasing. It surely vanishes, for $\displaystyle{\Delta_{s_n} = \frac{2s}{D_{s_n}}}$. Then, so does the sequence
$\displaystyle{\left(\{\sqrt{p_{s_n+1}}\Delta_{s_n}\}\right)_{s_n>1}}$
of fractional parts. The implication is that $\displaystyle{\sqrt{p_{s_n+1}}\Delta_{s_n}}$ tends to $s$ as $s_n$
tends to infinity. Moreover, on two successive values of $s_n$, say $s_n+1$ and $s_n+2$,
we have $\displaystyle{\{\sqrt{p_{s_n+2}}\Delta_{s_n+1}\} < \{\sqrt{p_{s_n+1}}\Delta_{s_n}\}}$ thereby,
$\displaystyle{\left(\{\sqrt{p_{s_n+1}}\Delta_{s_n}\}\right)_{s_n>1}}$ decreases. Consequently,  $\displaystyle{\left(\{\sqrt{p_{s_n}}\Delta_{s_n}\}\right)_{s_n>1}}$ is increasing and converges to $1$, and
$\displaystyle{\lim_{n\rightarrow\infty}\sqrt{p_{s_n+1}}\Delta_{s_n} = \lim_{n\rightarrow\infty}\sqrt{p_{s_n}}\Delta_{s_n} = s}$.

\section{The Fractional Part of $\displaystyle{\sqrt{p_n}}$}
\setcounter{equation}{0}
\renewcommand{\theequation}{\thesection.\arabic{equation}}
Recall that $\Delta_n := \sqrt{p_{n+1}} - \sqrt{p_n}$,
$D_n := \sqrt{p_{n+1}} + \sqrt{p_n}$, and $d_n := p_{n+1} - p_n$.
Define
$\mu_n$ and $\mu_{n+1}$ to stand for the fractional parts of $\sqrt{p_n}$ and $\sqrt{p_{n+1}}$, respectively.
They are also noted by $\{\sqrt{p_n}\}$ and $\{\sqrt{p_{n+1}}\}$. Thus,
$\displaystyle{\sqrt{p_n} = \lfloor\sqrt{p_n}\rfloor + \mu_n = \lfloor\sqrt{p_n}\rfloor + \{\sqrt{p_n}\}}$,
where $\lfloor\sqrt{p_n}\rfloor$ represents the integral part of $\sqrt{p_n}$.

\vskip 7pt\noindent
Suppose that $d_n = d_m \geq 2$. Then, comparing $D_n$ and $D_m$,
it easily appears that $\Delta_n = \Delta_m$, if and only if, $n = m$.
The fractional part of $\displaystyle{\sqrt{\frac{p_{n+1}}{p_n}}}$ equals $\Delta_n/\sqrt{p_n}$, for
$\displaystyle{\sqrt{\frac{p_{n+1}}{p_n}} - \frac{\Delta_n}{\sqrt{p_n}} = 1}$.
Here are some identities involving $\displaystyle{\sqrt{\frac{p_{n+1}}{p_n}}}$.
One can readily establish them using either Theorem 3.5
or by substitution of $\displaystyle{\sqrt{p_{n+1}}\Delta_n - \{\sqrt{p_{n+1}}\Delta_n\}}$ for $d_n/2$ or
$\displaystyle{\{\sqrt{p_{n+1}}\Delta_n\}}$ for $\displaystyle{1 - \{\sqrt{p_n}\Delta_n\}}$.

\vskip 7pt\noindent
$\displaystyle{1 \ =  \ \sqrt{\frac{p_{n+1}}{p_n}} - \frac{2\{\sqrt{p_{n+1}}\Delta_n\}}{\sqrt{p_n}\Delta_n}
\  =  \  \sqrt{\frac{p_{n+1}}{p_n}} - \frac{\Delta_n}{\sqrt{p_n}} \ = \ \sqrt{\frac{p_{n+1}}{p_n}} + \frac{2\{\sqrt{p_n}\Delta_n\}}{\sqrt{p_n}\Delta_n} - \frac 2{\sqrt{p_n}\Delta_n}.}$

\vskip 7pt\noindent
$\displaystyle{1 \  = \ \frac{d_n}{\sqrt{p_n}\Delta_n} \ - \  \sqrt{\frac{p_{n+1}}{p_n}}
\quad = \quad \frac{D_n}{\sqrt{p_n}} \ - \ \sqrt{\frac{p_{n+1}}{p_n}}
\quad = \quad \frac{d_n}{\sqrt{p_{n+1}}\Delta_n} \ - \  \sqrt{\frac{p_n}{p_{n+1}}}.}$

\vskip 7pt\noindent So
$\displaystyle{\frac 2{\sqrt{p_n}\Delta_n} - \frac{2\{\sqrt{p_n}\Delta_n\}}{\sqrt{p_n}\Delta_n} =
\frac{d_n}{\sqrt{p_n}\Delta_n} - 2 = \left\{\frac{d_n}{\sqrt{p_n}\Delta_n}\right\} = \frac{\Delta_n}{\sqrt{p_n}} \leq \frac{\Delta_4}{\sqrt{p_n}}}$.
In this sequence, either side of each equality signs
represents  $\displaystyle{\left\{\sqrt{\frac{p_{n+1}}{p_n}}\right\}}$.

\vskip 10pt
It is noteworthy that $\displaystyle{\frac{d_n}{\sqrt{p_n}\Delta_n} - 2 = \frac{\Delta_n}{\sqrt{p_n}}}$
implies $\displaystyle{\lim_{n\rightarrow\infty}d_n/\sqrt{p_n} = 2\lim_{n\rightarrow\infty}\Delta_n}$,
an identity met earlier.
Moreover, $\displaystyle{\sqrt{\frac{p_{n+1}}{p_n}} \leq \sqrt{5/3}}$ and
$\displaystyle{\sqrt{\frac{p_{n+1}}{p_n}} = 1 + \frac{\Delta_n}{\sqrt{p_n}}}$ entails
$\frac{\Delta_n}{\sqrt{p_n}} \leq \sqrt{\frac 53} - 1 < \frac 13$.
That means $\{\sqrt{p_{n+1}/p_n}\}$ never exceeds $1/3$.
Remarks on Theorem 2.2 prompt
$\displaystyle{\Delta_n/\sqrt{p_n} = \{\sqrt{p_{n+1}}/\sqrt{{p_n}}\} < 1/4\ }$ $(n\geq 5)$,
 since $\displaystyle{\sqrt{p_{n+1}/p_n} = 1 + \frac{\Delta_n}{\sqrt{p_n}} < \sqrt{\frac 32}}$.
 That is $\displaystyle{\sqrt{\frac{p_{n+1}}{p_n}} < \frac 54} \ (n\geq 5)$. Also,
$\displaystyle{\lim_{n\rightarrow\infty}\frac 1{\sqrt{p_n}\Delta_n} = }$
$\displaystyle{ \lim_{n\rightarrow\infty}\frac{\{\sqrt{p_n}\Delta_n\}}{\sqrt{p_n}\Delta_n} = 0,}$ for
$\displaystyle{\frac{1}{\sqrt{p_n}\Delta_n} - \frac 2{d_n} = \frac{\Delta_n}{d_n\sqrt{p_n}}}$
 and $\displaystyle{\frac 1{\sqrt{p_n}\Delta_n} - \frac{\{\sqrt{p_n}\Delta_n\}}{\sqrt{p_n}\Delta_n}
 = \frac{\Delta_n}{2\sqrt{p_n}}}$.

\subsection{Distance from $p_n$ to the Perfect Square Right Before}
Expressing $\sqrt{p_n}$ as a sum of its integral and fractional parts does not tell that much about direct relations
between these parts. We introduce another quantity that expresses somewhat the relative location of $\mu_n$
within the interval $ 0 < \mu_n < 1$.

\vskip 6pt \noindent {\thm\ }{\sl Given a prime $p_n$ $(n\geq 2)$ wedged between $N^2$ and $(N+1)^2$, there exists
 an integer $h_n$ such that $\displaystyle{\sqrt{p_n} \ = \ \frac{h_n - \mu_n^2}{2\mu_n} + \mu_n \ = \
\frac 12\left(\frac{h_n}{\mu_n} - \mu_n\right) + \mu_n}$. In addition, $N$ and $h_n$ (that stands for
the distance between $N^2$ and $p_n$) have opposite parities. Furthermore,
 $\displaystyle{(h_n - \mu_n^2)/\mu_n = h_n/\mu_n - \mu_n}$ is an even integer.
}

\vskip 10pt \noindent\dem\quad Let $\mu_n = \{\sqrt{p_n}\}$. From now on, there is a prime number between two
consecutive perfect squares. So, $\sqrt{p_n} = N + \mu_n$  implies $p_n - N^2 = 2\mu_n\sqrt{p_n} - \mu_n^2$.
Since $p_n$ is at an integer from $N^2$, we have $h_n := 2\mu_n\sqrt{p_n} - \mu_n^2$
where $h_n$ and $N$ are of different parities.
Thus, $2\sqrt{p_n} = (h_n + \mu_n^2)/\mu_n$ yields the identity of the theorem. And
 this states $\displaystyle{\frac{h_n - \mu_n^2}{\mu_n} = \frac{h_n}{\mu_n} - \mu_n}$ as twice the integral
 part of $\sqrt{p_n}$.
\hspace{\fill}$\Box$

\vskip 10pt\noindent The prime $p_n$ being between $N^2$ and $(N+1)^2$, the quantity $h_n$ is at most equal to $2N$.
Evidently, for parity reasons $h_n$ never takes on the value $N$.

\vskip 12pt As first estimate of $\mu_n$, write
$\displaystyle{h_n = p_n - N^2 = (\sqrt{p_n} - N)(\sqrt{p_n} + N)}$, where $N = \lfloor\sqrt{p_n}\rfloor$. Then,
$\displaystyle{h_n = \mu_n(\sqrt{p_n} + N)}$,
This gives $\mu_n = h_n/(\sqrt{p_n} + N)$, $h_n = 2\mu_nN + \mu_n^2$; and hence
\begin{eqnarray}
\frac{h_n}{2\sqrt{p_n}} \ < \mu_n \ < \ \frac{h_n}{2N} \ \leq \ 1 \mbox{ and }
2\mu_nN  \ < \  h_n  \ <  \ 2\mu_n\sqrt{p_n} = 2\mu_nN + 2\mu_n^2.
\end{eqnarray}

\noindent One may use $\displaystyle{\sqrt{p_n} = N\left(1 + h_n/N^2\right)^{1/2}}$
to obtain $\displaystyle{\frac{h_n}{2N} - \frac{h_n^2}{8N^4} \ < \ \mu_n \ < \ \frac{h_n}{2N}.
}$

\vskip 10pt
The next lemma asserts that two primes $p_n$ and $p_{n+1}$ are on
both sides of a perfect square, if and only if, $\displaystyle{\{\sqrt{p_n}\} > \{\sqrt{p_{n+1}}\}}$.
In particular, $\mu_n - \mu_{n+1} > 1-\Delta_4 = 0,32912\cdots.$

\vskip 10pt \noindent{\lem\ }{\sl $\displaystyle{\Delta_n := \sqrt{p_{n+1}} - \sqrt{p_n} = \mu_{n+1} - \mu_n}$,
if and only if, both primes share the same integral part. Otherwise,
$\displaystyle{\Delta_n = 1 + \mu_{n+1} - \mu_n}$, in which case $\mu_n$ exceeds $\mu_{n+1}$.
}

\vskip 10pt \noindent\dem\quad $\displaystyle{\Delta_n =
\lfloor{\sqrt{p_{n+1}}}\rfloor + \mu_{n+1} - \lfloor\sqrt{p_n}\rfloor - \mu_n }$, that is,
$\displaystyle{\Delta_n - (\mu_{n+1} - \mu_n)}$ equals
$\displaystyle{\lfloor{\sqrt{p_{n+1}}}\rfloor - \lfloor\sqrt{p_n}\rfloor}$.
Then, either $\displaystyle{\lfloor{\sqrt{p_{n+1}}}\rfloor = 1 + \lfloor\sqrt{p_n}\rfloor}$ or
$\displaystyle{\lfloor{\sqrt{p_{n+1}}}\rfloor = \lfloor\sqrt{p_n}\rfloor}$, and the conclusion follows.
\hspace{\fill}$\Box$

\vskip 10pt\noindent Consequently, the fractional part of $\sqrt{p_{n+1}}/\sqrt{p_n}$ turns out to be equal to
either $(\mu_{n+1} - \mu_n)/\sqrt{p_n}$ or $(1 + \mu_{n+1} - \mu_n)/\sqrt{p_n}$, according to whether
$\lfloor{\sqrt{p_n}}\rfloor = \lfloor{\sqrt{p_{n+1}}}\rfloor$ or not.
The largest among all primes between $N^2$ and $(N+1)^2$ has, obviously, the greatest fractional part
which exceeds also the fractional part of the prime right after $(N+1)^2$.
This last fact is the underlying idea of the identity $\Delta_n = 1 + \mu_{n+1} - \mu_n$.
Using the lemma the next result derives from Theorem 4.1.

\vskip 10pt \noindent {\thm\ }{\sl
$\displaystyle{\lfloor{h_n/\mu_n}\rfloor = 2\lfloor{\sqrt{p_n}}\rfloor = 2\sqrt{p_n - h_n}}$ and \\
$\displaystyle{\{h_n/\mu_n\} = \{\sqrt{p_n}\} = \mu_n = \sqrt{p_n} - \sqrt{p_n - h_n}}$. That is,
$\displaystyle{\frac{h_n}{\mu_n} \  = \ 2\sqrt{p_n - h_n} + \mu_n}$. In addition,

$\displaystyle{\mu_n \ = \ \frac{h_n}{\sqrt{p_n} + \sqrt{p_n-h_n}} \ = \
\frac{h_n}{\sqrt{p_n} + \lfloor{\sqrt{p_n}}\rfloor} \ = \ \frac{h_n}{2\sqrt{p_n} - \mu_n}}$,
with the estimates:

\begin{eqnarray}
\frac{h_n}{2\sqrt{p_n}} \ < \ \mu_n \ < \ \left\{ \begin{array}{ll}
\frac{h_n}{D_n - 1} = \frac{h_n}{2\sqrt{p_n} - 1 + \Delta_n} &  \mbox{ if } \quad
\lfloor{\sqrt{p_n}}\rfloor = \lfloor{\sqrt{p_{n+1}}}\rfloor,   \\
 &   \\
\frac{h_n}{2\sqrt{p_n} - 1}  &  \mbox{ if } \quad
 \lfloor{\sqrt{p_n}}\rfloor + 1 = \lfloor{\sqrt{p_{n+1}}}\rfloor.
\end{array} \right.
\end{eqnarray}
}

\vskip 5pt \noindent\dem\quad According to
 Theorem 4.1, $\displaystyle{\frac{h_n}{\mu_n} - \mu_n = 2N }$, i.e.
$\displaystyle{\lfloor{h_n/\mu_n}\rfloor = 2\lfloor{\sqrt{p_n}}\rfloor}$ and $\displaystyle{\{h_n/\mu_n\} = \mu_n}$.
Also, $p_n - N^2 = h_n$ implies $N = \sqrt{p_n - h_n} = \lfloor{\sqrt{p_n}}\rfloor$,
so that $\displaystyle{\mu_n = \sqrt{p_n} - \sqrt{p_n - h_n}}$. Then,

$\displaystyle{\sqrt{p_n} - \sqrt{p_n - h_n} =
\frac{h_n}{\sqrt{p_n} + \sqrt{p_n - h_n}} = \frac{h_n}{\sqrt{p_n} + \lfloor{\sqrt{p_n}}\rfloor} =
\frac{h_n}{2\sqrt{p_n} - \mu_n} >  \frac{h_n}{2\sqrt{p_n}}}$.
Suppose that $\lfloor{\sqrt{p_n}}\rfloor$ is equal to $\lfloor{\sqrt{p_{n+1}}}\rfloor$, then
 Lemma 4.2 induces $\Delta_n = \mu_{n+1} - \mu_n$, so that
$$ \mu_n = \frac{h_n}{2\sqrt{p_n} - \mu_n} = \frac{h_n}{2\sqrt{p_n} + \Delta_n - \mu_{n+1}}
< \frac{h_n}{2\sqrt{p_n} + \Delta_n - 1} = \frac{h_n}{D_n - 1}.$$ Inversely, $\Delta_n = 1 + \mu_{n+1} - \mu_n$ when
$\lfloor{\sqrt{p_n}}\rfloor$ and $\lfloor{\sqrt{p_{n+1}}}\rfloor$ are not equal. So,
$0 < \Delta_n - \mu_{n+1}$ yields
$\displaystyle{\mu_n = \frac{h_n}{2\sqrt{p_n} - \mu_n} = \frac{h_n}{2\sqrt{p_n} - 1 + \Delta_n - \mu_{n+1}}
< \frac{h_n}{2\sqrt{p_n} - 1}}$, say the claimed bound.
\hspace{\fill}$\Box$

\vskip 10pt After (4.1) and (4.2) that give bounds on $\mu_n$, we represent
 $\{\sqrt{p_n}\}$ in terms of a power series.

\vskip 10pt \noindent {\thm\ }{\sl Let $N\geq 2$ and $n\geq 3$ be two integers such that
$\displaystyle{p_n = N^2 + h_n}$. Then,
\\
$\displaystyle{
\mu_n \  = \  \frac{h_n}{2N} - \frac{h_n^2}{8N^3} + \frac{h_n^3}{16N^5} - \frac{5h_n^4}{128N^7} + \cdots
\  = \  \sum_{k=1}^\infty\frac{(-1)^{k-1}(2k)!}{4^k(2k-1)(k!)^2}\cdot\frac{h_n^k}{N^{2k-1}},
}$
where $h_n$ and $N$ are of opposite parities.
}

\vskip 10pt \noindent\dem\quad $\displaystyle{p_n = N^2 + h_n}$ implies $\displaystyle{\sqrt{p_n} = N(1+h_n/N^2)^{1/2}}$.
 And for $N$ greater than $1$, one has $1/N^2\leq h_n/N^2 < 1$ $(n\geq 3)$
so as to represent $\displaystyle{\sqrt{1+h_n/N^2}}$ through a convergent Taylor series.

\noindent
$\displaystyle{\left(1 + \frac{h_n}{N^2}\right)^{1/2} = 1 + \frac{h_n}{2N^2} - \frac{h_n^2}{8N^4} + \frac{h_n^3}{16N^6} - \cdots =
1 + \sum_{k=1}^\infty\frac{(-1)^{k-1}(2k)!}{4^k(2k-1)(k!)^2}\left(\frac{h_n}{N^2}\right)^k.}$
This gives, since
\\
$\sqrt{p_n} = N + \mu_n = N(1+h_n/N^2)^{1/2}$, the desired identity.
\hspace{\fill}$\Box$

\vskip 12pt \noindent There are obviously infinitely many values of $(h_n, N)$ for which the series
converges, but not to a fractional part of the square root of a prime number. For $n =1$, we deal with the unique situation in which
$h_n$ and $N$ are of the same parity, with $h_1 = N = 1$, i.e. $h_n/N^2 = 1$, $\mu_1 = \sqrt{2} - 1$. The case of $n = 2$, $h_2/N^2 = 2$
corresponds to $\mu_2 = \sqrt{3} - 1$. When $h_n = 1$ and $N$ represents an even integer, the series takes the form
$\displaystyle{\sum_{k=1}^\infty\frac{2(-1)^{k-1}(2k)!}{4^{2k}(2k-1)(k!)^2}\times\frac 1{N^{2k-1}}, \ \ N = 1, 2, 3, \cdots.}$

\vskip 10pt Now, we raise some questions about $h_n$, subject to proofs of conjectures 1.1 and 2.5.
\begin{enumerate}
\item[1). ] {\sl Given $m \geq 1$, does there exist a prime $q$ satisfying
$\displaystyle{q - \lfloor\sqrt{q}\rfloor^2  = m}$?} Put differently, does $h_n$ assume any positive integral value?
{\sl Moreover, do an infinity of primes $q$ satisfy the equality?}
\item[2). ] {\sl Given $N \geq 1$, does there exist a pair $(p_n , q)$ of primes such that $3\leq q\leq 4N-1$ and
$(2N)^2 = p_n - q$?} The condition on $q$ surely implies that $\lfloor\sqrt{p_n}\rfloor  = 2N$ and $q$ is the corresponding
$h_n$. To answer in the negative, one should find an integer $N$ such that the set of values of $h_n$, related to the
set of primes between $4N^2$ and $(2N+1)^2$, is composed of non-prime integers only. We narrowly missed such a batch of primes
when $N = 40$, with thirteen primes matching eleven composite $h_n$; the two intruders are $h_i = 73$ and $h_j = 151$.

Yet, circumstances in which the respective batch of $h_n$ is composed only of primes exist. For instance,
if $N = 6$, to all five primes: $149, 151, 157, $  $163$ and $167$ between $12^2$ and $13^2$ correspond five prime
values of $h_n$: $5,$ $7,$ $13,$ $19$ and $23$.
That is the case for $N = 15$ as well; here, eight prime numbers are connected with as many prime values of $h_n$.
The concern whether there exist infinitely many such batches is noticeable. Of course, $4N^2$ will be of the form
$4N^2 = q - p$ where $q$ and $p$ are primes and $\lfloor{\sqrt{q}}\rfloor = 2N$, should there exist at least one prime value of $h_n$
associated with a prime in-between $4N^2$ and $(2N+1)^2$. {\sl Does this representation hold for every even perfect square?}
\item[3). ] Consider the statement: {\sl Given $N\geq 2$, there exist $1\leq h_i , h_j \leq 2N-1$ and $r\geq 0$
such that $2N = (h_i-r) + (h_j+r)$, where $h_i-r$ and $h_j+r$ are prime numbers.} Of course $N^2+h_i$,
and $N^2+h_j$ are also prime numbers. Are there an infinity of integers $N$ validating the statement?
or, on the contrary, do an infinity of integers $N$, like $N=17, 19$ or $46$ deny it?
\end{enumerate}

\vskip 10pt We now approach $\displaystyle{\lim\Delta_n}$ as $n$ goes to infinity.
The idea behind the proof consists of splitting the set of primes into two groups:
$\displaystyle{I := \{p_n: \lfloor{\sqrt{p_n}}\rfloor = \lfloor\sqrt{p_{n+1}}\rfloor\}}$
 and  $\displaystyle{S := \{p_n: \lfloor{\sqrt{p_n}}\rfloor \neq \lfloor\sqrt{p_{n+1}}\rfloor\}}$.
 The quantity $\Delta_n$, in this way, assumes two forms, one on $I$ where $\mu_n < \mu_{n+1}$; and the other on $S$, with
  $\mu_n > \mu_{n+1}$.

\vskip 10 pt \noindent {\thm\ }{\sl
$\displaystyle{\lim_{n\rightarrow\infty}\left(\sqrt{p_{n+1}} - \sqrt{p_{n}}\right) = 0}$.
}

\vskip 10pt \noindent\dem\quad
$\displaystyle{\sqrt{\frac{p_{n+1}}{p_n}} = 1 + \frac{\Delta_n}{\sqrt{p_n}}}$
where $\Delta_n := \sqrt{p_{n+1}} - \sqrt{p_n}$, implies
$\displaystyle{\sqrt{\frac{p_{n+1}}{p_n}} - 1 = o(1)}$, as $n\rightarrow\infty$.
Setting $N = \lfloor{\sqrt{p_n}}\rfloor$, this takes two forms:
$\displaystyle{\frac{N + \mu_{n+1}}{N + \mu_n} - 1 = o(1)}$ or $\displaystyle{\frac{N + 1 +\mu_{n+1}}{N +\mu_n} - 1 = o(1)}$,
according as $\lfloor{\sqrt{p_{n+1}}}\rfloor = N$ or not. The former case gives $N + \mu_{n+1} \sim N + \mu_n$,
i.e. $\mu_{n+1} \sim \mu_n$. So, since $\mu_{n+1} - \mu_n < 1$, $\Delta_n = \mu_{n+1} - \mu_n = o(1)$
as $n$ goes to infinity. The latter case
leads to $N + 1 +\mu_{n+1} \sim N +\mu_n$, i.e $1 +\mu_{n+1} \sim \mu_n$, and $\Delta_n = 1 + \mu_{n+1} - \mu_n = o(1)$
$(n\rightarrow\infty)$. In either case, limit of $\Delta_n$ is zero as $n$ tends to infinity.
\hspace{\fill}$\Box$

\vskip 7pt\noindent The next is a consequence of Theorem 4.5 alongside (3.4).

\vskip 10pt\noindent{\cor\ }{\sl
$\displaystyle{\lim_{n\rightarrow\infty}\frac{d_n}{\sqrt{p_n}} = \lim_{n\rightarrow\infty}\{\sqrt{p_{n+1}}\Delta_n\} = 0}$;
hence, $\displaystyle{\lim_{n\rightarrow\infty}\{\sqrt{p_{n}}\Delta_n\} = 1}$.
}

\vskip 10pt\noindent{\cor\ }{\sl $\displaystyle{\liminf_{n\rightarrow\infty}\mu_n = 0} \ $ and
$ \ \displaystyle{\limsup_{n\rightarrow\infty}\mu_n = 1}$.
}

\vskip 10pt \noindent\dem\quad Consider the sequence of  primes $p_m$ where $p_m$ and $p_{m+1}$
are on both sides of $\displaystyle{(\lfloor{\sqrt{p_n}}\rfloor)^2}$. The first elements being
 $3, 7, 13, 23, \cdots$. For every $p_m$ from the sequence, $\Delta_m = 1 - \mu_m + \mu_{m+1}$. Since
 $\displaystyle{\lim_{m\rightarrow\infty}\Delta_m = 0}$, for any $\epsilon > 0$,
 an integer $m_0$ ($m_0 = m_0(\epsilon)$) exists such that
 $$
 1 - \mu_m + \mu_{m+1} < \epsilon \quad (\mbox{for all } m\geq m_0).
 $$
 That is to say, $1 - \epsilon < 1 - \epsilon + \mu_{m+1} < \mu_m$. Or more succinctly, $0 < 1 - \mu_m < \epsilon$
   holds for infinitely many $m$.
 The implication is $\displaystyle{\lim_{m\rightarrow\infty}\mu_m = 1}$,  that is
 $\displaystyle{\limsup_{n\rightarrow\infty}\mu_n = 1}$. Likewise, on
 using the same sequence, hence based on the displayed inequality, one obtains that
 $0 < \mu_{m+1} < \mu_m - 1 + \epsilon$ \ $(m > m_0)$ for a given $\epsilon > 0$. Now, since
 $\displaystyle{\limsup_{n\rightarrow\infty}\mu_n = 1}$, these inequalities yield
 $\displaystyle{\lim_{m\rightarrow\infty}\mu_{m+1} = 0}$, i.e. $\displaystyle{\liminf_{n\rightarrow\infty}\mu_n = 0}$.
\hspace{\fill}$\Box$

\subsection{Primes of the Shape $\displaystyle{N^2+1}$}
Clearly, a particular value of $h_n$ characterizes $p_n$ as the next integer to a perfect square.
The expression of $\sqrt{p_n}$ in  Theorem 4.1 is  bounded below by
$\displaystyle{\frac 12\left(\frac 1{\mu_n} - \mu_n\right) + \mu_n}$,
with equality holding, if and only if, $h_n = 1$, i.e. $p_n = N^2 + 1$.

\vskip 10pt \noindent {\cor\ }{\sl $\displaystyle{\sqrt{p_n} = \frac{1 - \mu_n^2}{2\mu_n} + \mu_n}$,
if and only if, $p_n$ is one away from a perfect square.
}

\vskip 10 pt \noindent Equivalence between $h_n/\mu_n - \mu_n = 2N$ and $h_n - 1 = 2\mu_nN - (1 - \mu_n^2)$
yields forthwith

\vskip 10pt \noindent {\cor\ }{\sl Consider $p_n$ as above. Then, $\displaystyle{\{2\mu_nN\} = 1 - \mu_n^2}$ and
the integral part of $2\mu_nN$ is equal to $h_n - 1$. Furthermore, $p_n$ equals $N^2+1$,
if and only if, $\mu_n^2 + 2\mu_nN - 1 = 0$.
}

\vskip 10pt Surely, $1 = 2\mu_n\sqrt{p_n} - \mu_n^2$ leads to $\mu_n\sqrt{p_n} = 1/2 + \mu_n^2/2 < 1$,
so one states

\vskip 10pt \noindent {\pro\ }{\sl  A prime $p_n$ is one more than a perfect square, if and only if,
$\displaystyle{\mu_n\sqrt{p_n} = \{\mu_n\sqrt{p_n}\}}$, i.e. if and only if,
$\displaystyle{2\mu_n\sqrt{p_n} - 1 = \{2\mu_n\sqrt{p_n}\}}$.
}

\vskip 10pt The pair $(3,5)$ is the only twin prime pair $(p_n , p_n+2)$
such that $p_{n+1} = p_n + 2$ is one unit away from a perfect square. Actually, if $N =  \lfloor{\sqrt{p_n}}\rfloor$,
since $\displaystyle{p_{n+1} - (N+1)^2 = d_n + (p_n - N^2) - (2N+1)}$, one has
$\displaystyle{p_{n+1} = (N+1)^2 + 1}$, if and only if, $\displaystyle{p_n - N^2 = 2N = h_n}$. Also,
$\displaystyle{p_n = (N+1)^2 - 1 =}$ $\displaystyle{ N(N+2)}$ gives $N = 1$ and $p_n = 3$.
Therefore, if $p_{n+1} = p_n + 2$ and $n\geq 3$,
then $p_{n+1}$ and $p_n$ are on the same side of any perfect square, i.e.
 $\lfloor{\sqrt{p_{n+1}}}\rfloor = \lfloor{\sqrt{p_n}}\rfloor$.

\vskip 10pt \noindent {\thm\ }{\sl Let $m > n \geq 2$ be two integers such that $p_n$ is
$1$ more than a perfect square, say, $N^2$; and $p_m$ is also $1$ more than another perfect square, say, $M^2$.
 Then, $\mu_m < \mu_n$. In particular, $\mu_r < 1/4$
for every integer $r$ for which $p_r$ is the next integer to a perfect square. So,
$\{4\sqrt{p_r}\} = 4\{\sqrt{p_r}\} = 4\mu_r$. Moreover, $\sqrt{p_r}$ and $\sqrt{p_{r+1}}$ have the same
integral part.
}

\vskip 10 pt \noindent\dem\quad We have $\frac 1{\mu_n} -\mu_n = 2N$ and
$\frac 1{\mu_m} -\mu_m = 2M \geq 2(N+2)$.
This implies $\frac 1{\mu_m} - \mu_m \geq \frac 1{\mu_n} - \mu_n + 4 > \frac 1{\mu_n} - \mu_n + 1$
that yields
$\displaystyle{0 < \mu_m - \mu_n + 1 < \frac 1{\mu_m} - \frac 1{\mu_n} = \frac{\mu_n - \mu_m}{\mu_m\mu_n}}$, i.e.
$\mu_m < \mu_n$. For $n = 3$, $\mu_3 = 0.236\cdots < 1/4,$ thus $\mu_n\leq\mu_3$ and
$\{4\sqrt{p_n}\} = 4\mu_n = 4\{\sqrt{p_n}\} < 1.$ Suppose, lastly, that $p_n = 1+N^2$ and $p_{n+1}$ is after $(N+1)^2$,
(i.e. $\sqrt{p_n}$ and $\sqrt{p_{n+1}}$ have different integral parts).
Then, $\sqrt{p_{n+1}} - \sqrt{p_n} = 1 + \mu_{n+1} - \mu_n > 3/4$, for $\mu_n < 1/4$. And $3/4$ exceeds the upper bound
in Andrica's problem. As a result, $\sqrt{p_{n+1}}$ and $\sqrt{p_n}$ must share their integral part,
and this ends the proof.
\hspace{\fill}$\Box$

\vskip 10pt \noindent
According to Theorem 4.11, should there be infinitely many primes of the form $1 + N^2$, it would set a proof of
$\displaystyle{\liminf_{n\rightarrow\infty}\mu_n = 0}$. The first statement of the theorem, of course, holds also for primes
$p_n$ and $p_m$ that are $h_0 \geq 2$ away from respectively $N^2$ and $M^2$. That means, $\mu_m < \mu_n$ when
$h_n = h_m = h_0 \geq 2$ with $\sqrt{p_n} = N$ and $\sqrt{p_m} = M$.

\section{Relations Between $\mu_n$, $\{\mu_n\sqrt{p_n}\}$ and the Parity of $\displaystyle{\lfloor\sqrt{p_n}\rfloor}$}
\setcounter{equation}{0}
\renewcommand{\theequation}{\thesection.\arabic{equation}}
Obviously, $\displaystyle{\{2\mu_n\sqrt{p_n}\} = \mu_n^2 < \mu_n}$, regardless of the parity of $h_n$.
The identity $2\mu_n\sqrt{p_n} - \mu_n^2 = h_n$ implies that
$h_n$ is an even integer, if and only if, $\displaystyle{\{\mu_n\sqrt{p_n}\} = \mu_n^2/2}$.
It also implies $\displaystyle{\mu_n\sqrt{p_n} - \frac{1 + \mu_n^2}2}$
$\displaystyle{= \frac{h_n - 1}2}$, so that $h_n$
is an odd integer, if and only if, $\displaystyle{\{\mu_n\sqrt{p_n}\} = \frac{1 + \mu_n^2}2}$. That is,
$\mu_n^2 = 2\{\mu_n\sqrt{p_n}\} - 1$. Clearly $p_n = 2$ satisfies this last equality.
Therefore, we state the following.

\vskip 10pt \noindent {\thm\ }{\sl Given a prime $p_n$,
 $\displaystyle{2\{\mu_n\sqrt{p_n}\} = \{2\mu_n\sqrt{p_n}\}}$, if and only if, $h_n$ is an even integer.
Viewed otherwise, $h_n$ is an even integer, if and only if, $\displaystyle{\{\mu_n\sqrt{p_n}\} < 1/2}$.
In particular, $\lfloor{\mu_n\sqrt{p_n}}\rfloor = h_n/2$. In other words,
$\lfloor\sqrt{p_n}\rfloor$ is an odd integer greater than $1$, if and only if,
$\displaystyle{\mu_n = \sqrt{2\{\mu_n\sqrt{p_n}\}}}$.
}

\vskip 10pt \noindent\dem\quad Firstly, $\displaystyle{\mu_n = \sqrt{\{2\mu_n\sqrt{p_n}\}}}$ stems from
$2\mu_n\sqrt{p_n} - \mu_n^2 = h_n$, for $h_n$ is an integer and $\mu_n^2 < 1$.
Then, out of $\mu_n\sqrt{p_n} - \mu_n^2/2 = h_n/2$, we see that $\displaystyle{2\{\mu_n\sqrt{p_n}\} = \mu_n^2}$, if and only if,
$h_n$ is an even integer. Comparing this to the expression of $\mu_n$  above that involves the square root sign,
appears the first identity of the theorem. Secondly, from
$\displaystyle{2\{\mu_n\sqrt{p_n}\} = \mu_n^2 < 1}$, comes $\{\mu_n\sqrt{p_n}\} < 1/2$. Finally,
since $\lfloor\sqrt{p_n}\rfloor$ and $h_n$ are of different parities, once more,
$\displaystyle{2\{\mu_n\sqrt{p_n}\} = \mu_n^2}$, say, the last identity.
\hspace{\fill}$\Box$

\vskip 10pt \noindent {\cor\ }{\sl $h_n$ represents an even integer, (i.e. $\lfloor\sqrt{p_n}\rfloor$ stands for
an odd integer), if and only if, $\displaystyle{\mu_n > 2\{\mu_n\sqrt{p_n}\}}$.
}

\vskip 10pt \noindent\dem\quad If $h_n$ is even, then $\displaystyle{\{\mu_n\sqrt{p_n}\} = \mu_n^2/2}$,
so that $\displaystyle{2\{\mu_n\sqrt{p_n}\} = \mu_n^2 < \mu_n}$.
Conversely, if $\mu_n$ exceeds $2\{\mu_n\sqrt{p_n}\}$ and $h_n$ were odd, then
calling on the comment at the beginning of this section, we would have
$\displaystyle{0 > 2\{\mu_n\sqrt{p_n}\} - \mu_n = (1+\mu_n^2) - \mu_n}$. We would, therefore, come by the nonsense
$1 + \mu_n^2 < \mu_n$. Consequently, $h_n$ must be even.
\hspace{\fill}$\Box$

\vskip 10pt The proof of the next result is similar to that of the last theorem. In fact,
$\displaystyle{\mu_n^2 = 2\{\mu_n\sqrt{p_n}\} - 1}$, if and only if,
$\displaystyle{2\{\mu_n\sqrt{p_n}\} - 1 = 2\mu_n\sqrt{p_n} - h_n}$. And this is rewritten as
$\displaystyle{(h_n - 1)/2 = \mu_n\sqrt{p_n} - \{\mu_n\sqrt{p_n}\}}$.

\vskip 10pt \noindent {\pro\ }{\sl Given a prime $p_n$ as above,
$\displaystyle{2\{\mu_n\sqrt{p_n}\} - 1 = \{2\mu_n\sqrt{p_n}\}}$, if and only if, $h_n$ is an odd integer.
In another formulation, $\displaystyle{\{\mu_n\sqrt{p_n}\} > 1/2}$ is a necessary and sufficient condition for
$h_n$ to be an odd integer. That is, the integral part of $\sqrt{p_n}$ is an even integer, if and only if,
$\displaystyle{\mu_n = \sqrt{2\{\mu_n\sqrt{p_n}\} - 1}}$.  In this case, $\lfloor{\mu_n\sqrt{p_n}}\rfloor = (h_n - 1)/2$.
}

\vskip 10pt \noindent The next result is an easy consequence of this proposition. Its proof
 proceeds from the simple inequality $(\mu_n-1)^2 > 0$. Indeed, $h_n$ is an odd integer, if and only if,
$\displaystyle{\{\mu_n\sqrt{p_n}\} = (1+\mu_n^2)/2}$, so that $\{\mu_n\sqrt{p_n}\}$ is greater than $\mu_n$, if and only if
$\mu_n^2 - 2\mu_n + 1 = (\mu_n - 1)^2$ is positive.

\vskip 10pt \noindent {\cor\ }{\sl $\displaystyle{\mu_n < \{\mu_n\sqrt{p_n}\}}$,
if and only if, $h_n$ stands for an odd integer. That is, this inequality is a necessary and sufficient condition for
$\lfloor{\sqrt{p_n}}\rfloor$ to be an even integer.
}

\vskip 10pt Like $\sqrt{p_n}$, $\mu_n\sqrt{p_n}$ increases throughout $N < \sqrt{p_n} < (N+1)$. Since
$\{\mu_n\sqrt{p_n}\} = \mu_n^2/2$ or $(1+\mu_n^2)/2$ whether $N$ is odd or not, the fractional part
$\{\mu_n\sqrt{p_n}\}$, similarly to $\mu_n$, increases also throughout the same range. So, one can
take advantage of anyone of Theorem 5.1 and Proposition 5.3. Really, by the fact that
$\{\mu_n\sqrt{p_n}\} < 1/2$ when $N$ is odd and exceeds $1/2$ otherwise,
we set a more precise statement.

\vskip 7pt \noindent {\thm\ }{\sl For primes between two consecutive odd perfect squares,
the fractional part $\{\mu_n\sqrt{p_n}\}$ is increasing. That is for odd $N$, in the range
$N^2 < p_n < p_{n+1} < (N+2)^2$, one has $\{\mu_n\sqrt{p_n}\} < \{\mu_{n+1}\sqrt{p_{n+1}}\} \ (n\geq 2)$.
To explain further, as $p_n$ increases to $(N+1)^2$ and beyond, $\{\mu_n\sqrt{p_n}\}$ increases through $1/2$ to near $1$.
}

\vskip 10pt The theorem has a forthright consequence.

\vskip 10pt \noindent {\cor\ }{\sl $\displaystyle{\{\mu_{n+1}\sqrt{p_{n+1}}\} - \{\mu_n\sqrt{p_n}\}
< 1/2 }$ $\displaystyle{ (n\geq 2)}$ holds, whenever
 $\lfloor{\mu_n\sqrt{p_n}}\rfloor = \lfloor{\mu_n\sqrt{p_{n+1}}}\rfloor$.
}

\vskip 10pt \noindent\dem\quad Let $N = \lfloor{\sqrt{p_n}}\rfloor = \lfloor{\sqrt{p_{n+1}}}\rfloor$.
We have $0 < \mu_n^2 < 1$, and $\mu_n^2/2 = \{\mu_n\sqrt{p_n}\}$ if $N$ is an odd integer, i.e. $h_n$ and $h_{n+1}$ are even integers.
This leads to $\{\mu_n\sqrt{p_n}\} < 1/2 \ (n\geq 2)$, and hence the inequality of the corollary.
If conversely, $N$ is an even integer, then $\{\mu_n\sqrt{p_n}\} = (1+\mu_n^2)/2$, so that
$1/2 < \{\mu_n\sqrt{p_n}\} < 1 \ (n>1)$. This case also yields the conclusion of the corollary.
\hspace{\fill}$\Box$

\subsection{$\displaystyle{\sqrt{p_n}}$ and $\displaystyle{\sqrt{p_{n+1}}}$ with a Common Integral Part}
Suppose $\displaystyle{N^2 < p_n < p_{n+1} < (N+1)^2}$ for a given integer $N\geq 1$. Then
$\sqrt{p_n} = N + \mu_n$, $\sqrt{p_{n+1}} = N + \mu_{n+1}$ so as to induce $D_n = 2N + \mu_{n+1} + \mu_n$.
Obviously, $p_{n+1}$ cannot be of the shape $N^2 + 1$. Then, the relative location
of $\mu_{n+1} + \mu_n$ to unity affects the parity of $\lfloor{D_n}\rfloor$. With these constraints on $p_n$
and $p_{n+1}$, we state the next results.

\vskip 10pt \noindent {\pro\ }{\sl $\displaystyle{D_n - \left(\mu_{n+1} + \mu_n\right)}$ equals twice
the common integral part of $\sqrt{p_n}$ and $\sqrt{p_{n+1}}$. Therefore,
$\displaystyle{\sqrt{p_n} = \frac{D_n}2 - \frac 12\left(\mu_{n+1} + \mu_n\right) + \mu_n}$ \ and \
 $\displaystyle{\sqrt{p_{n+1}} = \frac{D_n}2 - }$   $\displaystyle{\frac 12\left(\mu_{n+1} + \mu_n\right) + \mu_{n+1}}$.
}

\vskip 10pt \noindent {\cor\ }{\sl $\lfloor{D_n}\rfloor$ is even, if and only if, $\mu_{n+1} + \mu_n < 1$.
In this case,
$\displaystyle{\{D_n\} \ = \ \{\sqrt{p_{n+1}}\} +
\{\sqrt{p_n}\}  =}$  $\displaystyle{\mu_{n+1} + \mu_n.}$ That is to say \
$\displaystyle{D_n \ = \ 2(\sqrt{p_n} - \mu_n) \ + \ \mu_{n+1} + \mu_n \ = \
2(\sqrt{p_{n+1}} - \mu_{n+1}) + \mu_{n+1} + \mu_n.}$
 In contrast, $D_n$ has an odd integral part, if and only if
$\displaystyle{\{D_n\}  = \mu_{n+1} + \mu_n - 1 =  \Delta_n + }$  $\displaystyle{ 2\mu_n - 1.}$ That is,
\\
$\displaystyle{D_n  =  \frac{d_n}{\Delta_n} \ = \ 1 + 2(\sqrt{p_n} - \mu_n) + (\mu_{n+1} + \mu_n - 1)
\ = \ 1 + 2(\sqrt{p_{n+1}} - \mu_{n+1}) + (\mu_{n+1} + \mu_n - 1)}$.
}

\vskip 10pt Since the square roots of the primes $p_n$ and $p_{n+1}$ share their integral parts,
we have, in the first place, $\Delta_n = \mu_{n+1} - \mu_n$
so that $\{D_n\} = \mu_{n+1} + \mu_n$, if and only if, $0 < \Delta_n = \mu_{n+1} - \mu_n < 1 - 2\mu_n$.
In the second place, oppositely, $\{D_n\} = \mu_{n+1} + \mu_n - 1$, if and only if,
$1 - \Delta_n < 2\mu_n < 2 - \Delta_n$. In either situation, we used $0 < \{D_n\} < 1$.
Then the following is equivalent to the corollary.

\vskip 10pt \noindent {\cor\ }{\sl Let $p_n$ and $p_{n+1}$ be as above.

1. $\lfloor{D_n}\rfloor$ is an even integer, if and only if, $2\mu_n < 1 - \Delta_n$. In particular, $\mu_n < 1/2$.

\vskip 7pt
2. $\lfloor{D_n}\rfloor$ represents an odd integer, if and only if,
$2\mu_n > 1 - \Delta_n. $
}

\vskip 10pt \noindent {\pro\ }{\sl $d_n = h_{n+1} - h_n$, $\displaystyle{\quad \{h_n/\mu_n\} = \{\sqrt{p_n}\} }$,
whenever $\lfloor{\sqrt{p_n}}\rfloor = \lfloor{\sqrt{p_{n+1}}}\rfloor$; \ and,

1. $\displaystyle{\Delta_n = \mu_{n+1} - \mu_n = \frac{h_{n+1}}{\mu_{n+1}} - \frac{h_n}{\mu_n} =
\left\{\frac{h_{n+1}}{\mu_{n+1}}\right\} - \left\{\frac{h_n}{\mu_n}\right\}}$

\vskip 7pt
2. $\displaystyle{\frac{d_n}2 = \lfloor{\mu_{n+1}\sqrt{p_{n+1}} - \mu_n\sqrt{p_n}}\rfloor}$,
and $\displaystyle{(\mu_{n+1}^2 - \mu_n^2)/2 = \{\mu_{n+1}\sqrt{p_{n+1}} - \mu_n\sqrt{p_n}\}}$. That is,
\\
$\displaystyle{\qquad\qquad \frac{d_n}2 \ \ = \ \ \mu_{n+1}\sqrt{p_{n+1}} - \mu_n\sqrt{p_n} \ - \
(\mu_{n+1}^2 - \mu_n^2)/2.}$
}

\vskip 10pt \noindent\dem\quad
The statement involving $h_n/\mu_n$ stems from $\{h_n/\mu_n\} = \mu_n$,
for $h_n/\mu_n = 2(\sqrt{p_n} - \mu_n) + \mu_n$.
As to the other assertion, we firstly
have $d_n = h_{n+1} - h_n$. Secondly,
using $h_n = 2\mu_n\sqrt{p_n} - \mu_n^2$ leads to
$h_{n+1} - h_n = 2(\mu_{n+1}\sqrt{p_{n+1}} - \mu_n\sqrt{p_n}) + \mu_n^2 - \mu_{n+1}^2$, and this establishes the identity.
\hspace{\fill}$\Box$

\vskip 10pt\noindent  The next consequence of the proposition holds, no matter what the parity of $h_n$.

\vskip 10pt \noindent {\cor\ }{\sl $\displaystyle{\{\mu_{n+1}\sqrt{p_{n+1}} - \mu_n\sqrt{p_n}\} = \{\mu_{n+1}\sqrt{p_{n+1}}\} - \{\mu_n\sqrt{p_n}\}}$, if and only if,
$\displaystyle{\lfloor{\sqrt{p_n}}\rfloor = \lfloor{\sqrt{p_{n+1}}}\rfloor}$. Similarly,
$\displaystyle{\lfloor{\mu_{n+1}\sqrt{p_{n+1}} - \mu_n\sqrt{p_n}}\rfloor}$ equals both $d_n/2$ and the difference
$\displaystyle{\lfloor{\mu_{n+1}\sqrt{p_{n+1}}}\rfloor - \lfloor{\mu_n\sqrt{p_n}}\rfloor}$.
}

\vskip 10pt \noindent\dem\quad Since
$\displaystyle{\mu_{n+1}\sqrt{p_{n+1}} - \mu_n\sqrt{p_n} - \frac{(\mu_{n+1}^2 - \mu_n^2)}2 = \frac{(h_{n+1} - h_n)}2}$,
its right-hand side equals an integer, if and only if, $h_{n+1}$ and $h_n$ are of the same parity. Thus,
the second statement of the proposition yields,
$\displaystyle{\{\mu_{n+1}\sqrt{p_{n+1}} - \mu_n\sqrt{p_n}\} = \frac{\mu_{n+1}^2}2 - \frac{\mu_n^2}2}$
which equals $\displaystyle{= \frac{\mu_{n+1}^2+1}2 - \frac{\mu_n^2+1}2 = }$
$\displaystyle{ \{\mu_{n+1}\sqrt{p_{n+1}}\} - \{\mu_n\sqrt{p_n}\}.}$
So it follows that,
$\displaystyle{\lfloor{\mu_{n+1}\sqrt{p_{n+1}} - \mu_n\sqrt{p_n}}\rfloor = \lfloor{\mu_{n+1}\sqrt{p_{n+1}}}\rfloor -
\lfloor{\mu_n\sqrt{p_n}}\rfloor = d_n/2.}$
\hspace{\fill}$\Box$

\subsection{Two Successive Primes on Either Side of a Perfect Square}
In this subsection, the opposite condition is set on $p_n$ and $p_{n+1}$:
requiring the integral parts of $\sqrt{p_n}$ and $\sqrt{p_{n+1}}$
to be successive integers. This comes down to $\Delta_n = 1 + \mu_{n+1} - \mu_n$.

\vskip 7pt\noindent {\pro\ }{\sl For $p_n$ and $p_{n+1}$ on either side of $(N+1)^2$,

1. $\displaystyle{\Delta_n = \mu_{n+1} + 1 - \mu_n = \frac{h_{n+1}}{\mu_{n+1}} - 1 - \frac{h_n}{\mu_n} =
\left\{\frac{h_{n+1}}{\mu_{n+1}}\right\} + 1 - \left\{\frac{h_n}{\mu_n}\right\}}$.

\vskip 7pt
2. $\displaystyle{d_n/2 - (\sqrt{p_n} - \mu_n)  \ =
\ \mu_{n+1}\sqrt{p_{n+1}} - \mu_n\sqrt{p_n} - (\mu_{n+1}^2 - \mu_n^2 - 1)/2 \quad (n\geq 2).}$
}

\vskip 10pt \noindent\dem\quad
$\displaystyle{h_{n+1} = 2\mu_{n+1}\sqrt{p_{n+1}} - \mu_{n+1}^2}$ implies that
$\displaystyle{\frac{h_{n+1}}{\mu_{n+1}} = 2\sqrt{p_{n+1}} - \mu_{n+1}}$, i.e. \\
$\displaystyle{\frac{h_{n+1}}{\mu_{n+1}} - \mu_{n+1} \ = \ 2(\sqrt{p_{n+1}} - \mu_{n+1}) \ = \
2(N + 1) =  \frac{h_n}{\mu_n} - \mu_n + 2}$. So,
$\displaystyle{\mu_n = \left\{\frac{h_n}{\mu_n}\right\}}$,
$\displaystyle{\mu_{n+1} = \left\{\frac{h_{n+1}}{\mu_{n+1}}\right\}}$ and
$\displaystyle{\frac{h_{n+1}}{\mu_{n+1}} - \frac{h_n}{\mu_n} - 1 = \mu_{n+1} + 1 - \mu_n = \Delta_n}$,
namely, the first statement. Next, squaring the identities
\\
$\displaystyle{\sqrt{p_{n+1} - h_{n+1}} =  N + 1 = 1 + \sqrt{p_n - h_n}}$ gives
$$
d_n \ = \ h_{n+1} - h_n + 1 + 2\sqrt{p_n - h_n} \ = \ h_{n+1} -
h_n + 1 + 2(\sqrt{p_n} - \mu_n).
$$
Again, $\displaystyle{h_{n+1} = 2\mu_{n+1}\sqrt{p_{n+1}} - \mu_{n+1}^2}$ leads to
$\displaystyle{h_{n+1} - h_n = 2(\mu_{n+1}\sqrt{p_{n+1}} - \mu_n\sqrt{p_n}) + \mu_n^2 - \mu_{n+1}^2}$.
Writing this down into the display formula yields the second statement and completes the proof.
\hspace{\fill}$\Box$

\vskip 10pt \noindent {\cor\ }{\sl Let $p_n$ and $p_{n+1}$ be like in the proposition. Then,

1. $\displaystyle{\left\{\mu_n\sqrt{p_n} - \mu_{n+1}\sqrt{p_{n+1}}\right\} = 1 +
\left\{\mu_n\sqrt{p_n}\right\} - \left\{\mu_{n+1}\sqrt{p_{n+1}}\right\}}$, if and only if, $h_n$
is an even integer.

$\quad$ Furthermore,
$\displaystyle{ \lfloor{\mu_n\sqrt{p_n} - \mu_{n+1}\sqrt{p_{n+1}}}\rfloor =
\lfloor{\mu_n\sqrt{p_n}}\rfloor - \lfloor{\mu_{n+1}\sqrt{p_{n+1}}}\rfloor - 1}$
then holds.

\vskip 8pt
2. $\displaystyle{\left\{\mu_n\sqrt{p_n} - \mu_{n+1}\sqrt{p_{n+1}}\right\} =
\left\{\mu_n\sqrt{p_n}\right\} - \left\{\mu_{n+1}\sqrt{p_{n+1}}\right\}}$ otherwise; hence

\vskip 5pt $\quad $
$\displaystyle{\lfloor{\mu_n\sqrt{p_n} - \mu_{n+1}\sqrt{p_{n+1}}}\rfloor =
\lfloor{\mu_n\sqrt{p_n}}\rfloor - \lfloor{\mu_{n+1}\sqrt{p_{n+1}}}\rfloor}$.
}

\vskip 10pt \noindent\dem\quad Since only one of the two integers $h_{n+1}$ and $h_n$ is odd, we have
$$
\mu_n\sqrt{p_n} \ - \ \mu_{n+1}\sqrt{p_{n+1}} \ + \ (1 + \mu_{n+1}^2 - \mu_n^2)/2 \ = \ (h_n - h_{n+1} + 1)/2.
$$
So, $\displaystyle{\left\{\mu_n\sqrt{p_n} - \mu_{n+1}\sqrt{p_{n+1}}\right\} = 1 - (1 + \mu_{n+1}^2 - \mu_n^2)/2}$, for
$(h_n - h_{n+1} + 1)/2$ represents an integer and $1/2 > (1 + \mu_{n+1}^2 - \mu_n^2)/2 > 0$.
Therefore, if $h_n$ is even, then \\
$\displaystyle{\mu_n\sqrt{p_n} - \mu_{n+1}\sqrt{p_{n+1}} - \left(1 + \frac{\mu_n^2}2 - \frac{1 +
\mu_{n+1}^2}2\right) = \frac{h_n- h_{n+1}+1}2 - 1 = \frac{h_n}2 - \frac{h_{n+1}-1}2 - 1}$,
which is the first statement. The case in which $h_n$ is odd gives
$$
\mu_n\sqrt{p_n} - \mu_{n+1}\sqrt{p_{n+1}} - \left(\frac{1 + \mu_n^2}2 - \frac{\mu_{n+1}^2}2
\right) = \frac{h_n- h_{n+1}+1}2 - 1 = \frac{h_n - 1}2 - \frac{h_{n+1}}2
$$
say, the second statement of the corollary.
\hspace{\fill}$\Box$

\vskip 10pt \noindent{\pro\ }{\sl Under the hypothesis of the proposition,
$\displaystyle{\sqrt{p_n} - \mu_n  =  \frac{d_n}{2\Delta_n} - \frac 12\left(\mu_{n+1} + \mu_n + 1\right) =}$ \\
 $\displaystyle{
\frac{D_n}2 - \frac 12\left(\mu_{n+1} + \mu_n + 1\right) = \frac{D_n}2 - \frac 12\left(\Delta_n + 2\mu_n\right)}$.
 Correspondingly,
$\displaystyle{\sqrt{p_{n+1}} =  \frac{d_n}{2\Delta_n} - \frac 12\left(\mu_{n+1} + \mu_n - 1\right) + }$
\\
$\displaystyle{\mu_{n+1} =
\frac{D_n}2 - \frac 12\left(\mu_{n+1} + \mu_n - 1\right) +
\mu_{n+1} =  \frac{D_n}2 - \frac 12\left(\Delta_n + 2\mu_n - 2\right) + \mu_{n+1}}$.
}

\vskip 10pt \noindent\dem\quad As the primes are on both sides of a perfect square, say, $(N+1)^2$, one has \\
$\displaystyle{\sqrt{p_{n+1}} = N + 1 + \mu_{n+1} = \sqrt{p_n} + 1 + \mu_{n+1} - \mu_n}$, i.e.
 $\displaystyle{1 - (\sqrt{p_{n+1}} - \sqrt{p_n}) = 1 - \Delta_n = \mu_n - \mu_{n+1}}$.  Out of this, along with
$\displaystyle{\Delta_n^2/2 = 1 - \{\sqrt{p_n}\Delta_n\} = d_n/2 - \sqrt{p_n}\Delta_n}$,  follows
$\displaystyle{\mu_n\Delta_n - \mu_{n+1}\Delta_n = \Delta_n - \Delta_n^2 = \Delta_n + }$
$\displaystyle{2\{\sqrt{p_n}\Delta_n\} - 2 = \Delta_n + 2\sqrt{p_n}\Delta_n - d_n.}$
Thus, $\displaystyle{\sqrt{p_n} = d_n/(2\Delta_n) + \mu_n/2 - \mu_{n+1}/2 - 1/2}$, that we rewrite,
recalling $d_n = D_n\Delta_n$, to obtain the identities concerning $\sqrt{p_n}$.
To prove the identities related to $\sqrt{p_{n+1}}$, we use those involving $\sqrt{p_n}$ with
$\displaystyle{\sqrt{p_{n+1}} - \sqrt{p_n} = 1 + \mu_{n+1} - \mu_n}$.
\hspace{\fill}$\Box$

\vskip 10pt Obviously, $\displaystyle{\sqrt{p_n} = N + \mu_n}$ and
$\displaystyle{\sqrt{p_{n+1}} = N + 1 + \mu_{n+1}}$, so that $D_n = 2N + 1 + \mu_{n+1} + \mu_n$.
Thus, $\displaystyle{\{D_n\} = \mu_{n+1} + \mu_n}$, if and only if, $\mu_{n+1} + \mu_n < 1$, i.e.
if and only if, $\lfloor{D_n}\rfloor = 2N + 1$.

\vskip 10pt\noindent{\cor\ }{\sl Let $p_n$ and $p_{n+1}$ be on both sides of a perfect square.
The integral part of $D_n$ is an even integer, if and only if, $\mu_{n+1} + \mu_n > 1$, i.e. if and only if,
$\displaystyle{\{D_n\} =  \mu_{n+1} + \mu_n - 1 = \Delta_n + 2\mu_n - 2}$. Then, one has,
$\displaystyle{D_n = 2 + 2(\sqrt{p_n} - \mu_n) + (\mu_{n+1} + \mu_n - 1) =
2(\sqrt{p_{n+1}} - \mu_{n+1}) + (\mu_{n+1} + \mu_n - 1)}$.

\vskip 10pt\noindent In contrast, $\lfloor{D_n}\rfloor$ represents an odd integer, if and only if,
$\mu_{n+1} + \mu_n < 1$, i.e. if and only if, \\
$\displaystyle{\{D_n\} =  \mu_{n+1} + \mu_n = \Delta_n + 2\mu_n - 1}$. This situation induces
$\displaystyle{D_n = 1 + 2(\sqrt{p_n} - \mu_n) + (\mu_{n+1} + \mu_n) = }$ \\
 $\displaystyle{2(\sqrt{p_{n+1}} - \mu_{n+1}) - 1 + (\mu_{n+1} + \mu_n)}$.
}

\vskip 10pt In the following, we set bounds on $\mu_n$ and $\mu_{n+1}$, when $p_n$ 
is the largest prime before $(N+1)^2$ and $p_{n+1}$ the first prime after. 
The bounds depend on the parity of the integral part of $D_n$.

\vskip 10pt\noindent{\cor\ }{\sl Suppose that $p_n$ and $p_{n+1}$ are as above.

1. $D_n$ has an even integral part, if and only if, $2\mu_{n+1} > \Delta_n$.
That is, $\mu_n > 1 - \Delta_n/2 \geq 1 - \Delta_4/2 = $

$\quad 0.6645\cdots .$

\vskip 7pt
2. The integral part of $D_n$ represents an odd integer, if and only if, $2\mu_{n+1} < \Delta_n$.
That is, $\mu_n < 1 - $

$\quad \Delta_n/2.$ In this case, $\mu_{n+1} < \Delta_4/2 = 0.3354\cdots .$
}

\vskip 10pt \noindent\dem\quad Assume that $D_n$ has an even integral part. Then, according to the previous corollary,
$\mu_{n+1} + \mu_n > 1$, so that Lemma 4.2 gives $2\mu_{n+1} > 1 - \mu_n + \mu_{n+1} = \Delta_n$, i.e. $\mu_{n+1} > \Delta_n/2$.
Again, $\mu_{n+1} = \Delta_n - 1 + \mu_n > \Delta_n/2$ implies $\mu_n > 1 - \Delta_n/2 \geq 1 - \Delta_4/2$.
Conversely, since $D_n = 2\lfloor\sqrt{p_n}\rfloor + 1 + \mu_{n+1} + \mu_n$, the estimate $2\mu_{n+1} > \Delta_n$ yields
$\{D_n\} = \mu_{n+1} + \mu_n - 1$, thus $\lfloor D_n\rfloor$ is even. This proves the first proposition.
To prove the second statement, we proceed on the same lines with $\mu_{n+1} + \mu_n < 1$.
\hspace{\fill}$\Box$

\vskip 10pt {\sl Does $\mu_n > 2\mu_{n+1}$ hold for two consecutive primes on both sides of a perfect square?}
This would be deduced from $1 - \Delta_4 > \mu_{n+1}$. Surely, this and
$1 + \mu_{n+1} - \mu_n = \Delta_n < \Delta_4$ imply that $\mu_n > (1 - \Delta_4) + \mu_{n+1} > 2\mu_{n+1}$.
Unfortunately, the estimate $1 - \Delta_4 > \mu_{n+1}$ does not hold for $n = 9$ i.e. for $p_n = 23$,
 ($\mu_{n+1} = 0.3851\cdots > 1 - \Delta_4 = 0.3293\cdots)$, even though $\mu_9 > 2\mu_{10}$. Note that
$\displaystyle{\limsup_{n\rightarrow\infty}\mu_n = 1}$ certainly induces the existence of an infinity of $n$ such that
$\mu_n > 2\mu_{n+1}$.

\section{$\Delta_n$ Versus $\mu_n$,  $d_n$ Versus $h_n$}
\setcounter{equation}{0}
\renewcommand{\theequation}{\thesection.\arabic{equation}}
$\displaystyle{d_n = \Delta_nD_n = \Delta_n (2\sqrt{p_n} + \Delta_n)}$ implies
$\displaystyle{\Delta_n  = \frac{d_n - \Delta_n^2}{2\sqrt{p_n}}}$.
Suppose that $\lfloor{\sqrt{p_n}}\rfloor = \lfloor{\sqrt{p_{n+1}}}\rfloor = N$.
Then, $\displaystyle{d_n = h_{n+1} - h_n}$ and $\displaystyle{\Delta_n = \mu_{n+1} - \mu_n}$, so
\begin{eqnarray}
\Delta_n \ = \ \mu_{n+1} - \mu_n \ = \ \frac{h_{n+1} - h_n}{2\sqrt{p_n}} \ - \ \frac{\Delta_n^2}{2\sqrt{p_n}} \ = \
\frac{d_n - \Delta_n^2}{2\sqrt{p_n}}.
\end{eqnarray}

The case of $h_{n+1} = p_{n+1} - (N+1)^2$ and $h_n = p_n - N^2$ yields
\begin{eqnarray}
d_n = h_{n+1} - h_n + 2N + 1 = h_{n+1} - h_n + 1 + 2(\sqrt{p_n} - \mu_n) = \Delta_n (2\sqrt{p_n} + \Delta_n).
\end{eqnarray}
 That is, $d_n - \Delta_n^2 = h_{n+1} - h_n + 1 + 2\sqrt{p_n} - (2\mu_n + \Delta_n^2) = 2\sqrt{p_n}\Delta_n$.
 This comes to
\begin{eqnarray}
\Delta_n \ = \ 1 + \mu_{n+1} - \mu_n \ = \ 1 - \ \frac{h_n - h_{n+1} - 1}{2\sqrt{p_n}} \ - \ \frac{2\mu_n + \Delta_n^2}{2\sqrt{p_n}}
\ = \ \frac{d_n - \Delta_n^2}{2\sqrt{p_n}}.
\end{eqnarray}

\vskip 10pt We sum up the two expressions of $d_n$ in terms of $h_n$.
\begin{eqnarray*}
d_n \quad = \quad \left\{ \begin{array}{ll} h_{n+1} - h_n, &
\mbox{ if \ $\lfloor{\sqrt{p_n}}\rfloor = \lfloor{\sqrt{p_{n+1}}}\rfloor$}   \\
2\lfloor{\sqrt{p_n}}\rfloor + 1 + h_{n+1} - h_n, & \mbox{ otherwise.} \end{array} \right.
\end{eqnarray*}
Yet, $d_n$ never equals $h_{n+1}$. In one case, the $n$th prime gap is greater than $h_{n+1}$, and less
than $h_{n+1}$ in the other. With regard to $h_n$, the equality $d_n = h_n$ occurs. For instance, when $n = 6$, one has
$p_6 = 13$ with $d_6 = 4 = h_6$. That is the case for $n = 24$ as well, $p_{24} = 89$ and $d_{24} = h_{24} = 8.$
In fact, if $\lfloor{\sqrt{p_n}}\rfloor$ $ (n\geq 1)$ stands for an odd integer,
then, $d_n = h_n$ if and only if either proposition holds: $2h_n = h_{n+1} $ or
$2h_n = 2\lfloor{\sqrt{p_n}}\rfloor + 1 + h_{n+1}.$ This latter equality should hold for only a
finite number of integers $n$. The former implies $2\mu_n > \mu_{n+1}$. Indeed,
$h_n = 2\mu_n\lfloor{\sqrt{p_n}}\rfloor + \mu_n^2$ is used to express $2h_n = h_{n+1}$
in terms of $\mu_n$ and $\mu_{n+1}$. Thus,
$2\mu_{n+1}\lfloor{\sqrt{p_n}}\rfloor + \mu_{n+1}^2 = 2(2\mu_n\lfloor{\sqrt{p_n}}\rfloor + \mu_n^2)$.
That is $(2\lfloor{\sqrt{p_n}}\rfloor + \mu_{n+1} + 2\mu_n)(\mu_{n+1} - 2\mu_n) + 2\mu_n^2 = 0$;
then the estimate $\mu_{n+1} - 2\mu_n < 0$ emerges.

\subsection{Bounds on $d_n$ When $\lfloor{\sqrt{p_n}}\rfloor + 1 = \lfloor{\sqrt{p_{n+1}}}\rfloor$}
Clearly $h_{n+1} < d_n$ if and only if $p_n$ and $p_{n+1}$
are on either side of a perfect square. Let $N = \lfloor{\sqrt{p_n}}\rfloor$. If $N\geq 2$ is even,
then $\displaystyle{h_n - h_{n+1} \leq 2N - 1 - h_{n+1}}$.
This, along with (6.2) induces $d_n \geq 2 + h_{n+1}$. When $N\geq 3$ is odd, surely
$\displaystyle{h_n - h_{n+1} \leq 2N - 2 - h_{n+1}}$ so that $d_n \geq 3 + h_{n+1}$, thereby
\begin{eqnarray}
d_n \quad \geq \quad \left\{ \begin{array}{ll} 2 + h_{n+1}, &  \mbox{ if \ $N\geq 2$ \ is even}   \\
3 + h_{n+1}, & \mbox{ if \ $N\geq 3$ \ is odd}   \\
2,  & \mbox{ if \ $N = 1$}
\end{array} \right.
\end{eqnarray}
holds for $\sqrt{p_{n+1}}$ and $\sqrt{p_n}$ located on both sides of $N+1$.
In that context, $d_n$ is less than or equal to $2\lfloor{\sqrt{p_n}}\rfloor$ because of Legendre's conjecture.
Also, as $h_{n+1} \geq 1$, one has $d_n \geq 2\lfloor{\sqrt{p_n}}\rfloor + 2 - h_n$. Then,
$2\lfloor{\sqrt{p_n}}\rfloor - (h_n - 2) \leq d_n \leq 2\lfloor{\sqrt{p_n}}\rfloor$.

\vskip 6 pt \noindent {\thm\ }
{\sl $2\lfloor\sqrt{p_n}\rfloor$ represents the maximum value of $d_n$. Also,
$\lfloor\sqrt{p_n}\rfloor + 1 = \lfloor\sqrt{p_{n+1}}\rfloor$ is a necessary condition for $d_n$
to hit $2\lfloor\sqrt{p_n}\rfloor$.
 }

\vskip 10pt\noindent\dem\quad If $d_n \geq 2N = 2\lfloor\sqrt{p_n}\rfloor$, then $(N+1)^2$
is necessarily located between $p_n$ and $p_{n+1}$.
Suppose  $d_n \geq 2N+2$.  This implies $2N+2 < \sqrt{2p_{n+1}} = \sqrt{2}(N+1+\mu_{n+1}).$
That is, $2N+2 < \sqrt{2}(N+2)$ or $N < (\sqrt{2} - 1)(2 + \sqrt{2}) < 3/2.$
We verify that for $N = 1$, $d_2 = 2 = 2\lfloor\sqrt{p_2}\rfloor$ is less than $2N+2$.
Consequently, for every integer $n$, $d_n \leq 2\lfloor\sqrt{p_n}\rfloor$.
\hspace{\fill}$\Box$

\vskip 10pt \noindent It appears from the proof that when $N =1$, we have $h_2 = 2 = 2N$ and
$h_3 = 1 = 2N-1$. In particular for $n > 2$,  $d_n \leq 2\lfloor{\sqrt{p_n}}\rfloor - 2$ and $h_n \geq 4$.
 One quickly derives, on using (6.2), the upcoming inequality.

\vskip 10pt \noindent {\cor\ }
{\sl Consider $p_n$ and $p_{n+1}$  as above. Then, $h_n \geq h_{n+1}+1$ for every integer $n\geq 2.$
 }

\vskip 10pt\noindent The following result characterizes the integers $N$ for
which $d_n = 2\lfloor\sqrt{p_n}\rfloor$; especially, since $d_n = 2\lfloor\sqrt{p_n}\rfloor$
if and only if $h_n = h_{n+1} + 1$.

\vskip 10pt \noindent {\thm\ }
{\sl Let $N$ be an integer such that $\displaystyle{1 \leq N^2 < p_n < (N+1)^2 < p_{n+1}}$.
Then, $h_n = h_{n+1} + 1$, if and only if, $N = 1$ or $N = 2$; i.e. $n = 2$ or $n = 4$. }

\vskip 10 pt\noindent\dem\quad If $N=1$ then, $h_2 = 2$, $h_3 = 1$ with $d_2 = 2 = 2N = h_{n+1} + h_n - 1$.
When $N = 2$, we have $h_4 = 3$, $h_5 = 2$ with $d_4 = 4 = 2N = h_{n+1} + h_n - 1$. In each case, $h_n = h_{n+1} + 1$,
so $N = 1$ and $N = 2$ satisfy the theorem. First, $d_n \leq \sqrt{2p_{n+1}} =
\sqrt{2}(N+1+\mu_{n+1}) < \sqrt{2}(N+2)$, where $N+1$ represents the integral part of $\sqrt{p_{n+1}}$.
Applying this to the leftmost equality of (6.2), along with the condition $h_n - h_{n+1} - 1 = 0$,
yields $N < 2(1 + \sqrt{2})$. That is, $N$ is less than or equal to $4$; and only $N = 1$ and $N = 2$
are such that the identity $h_n = h_{n+1} + 1$ is satisfied. This completes the proof.
\hspace{\fill}$\Box$

\vskip 10pt Neither $2N = h_{n+1} + h_n - 1$ nor $d_n = h_{n+1} + h_n - 1$
is equivalent to the theorem, although these identities are true for
$N = 1$ and $N = 2$. Actually, if $N = 7$ one has $h_{18} = 12$, $h_{19} = 3$
with $2N = 14 = h_{n+1} + h_n - 1$. Similarly, for $N = 5$, $h_{11} = 6$, $h_{12} = 1$
so that $d_{11} = 6 = h_{n+1} + h_n - 1$. In both examples
$h_n = h_{n+1} + 1$ is untrue, hence as is $d_n = 2N$. But, having at once $h_n = N + 1$ and
$h_{n+1} = N$ is obviously equivalent to the theorem.

\vskip 6 pt \noindent{\cor\ }{\sl $d_n = 2\lfloor\sqrt{p_n}\rfloor$, if and only if, $h_n =
\lfloor\sqrt{p_n}\rfloor + 1$ and $h_{n+1} = \lfloor\sqrt{p_n}\rfloor$.
 }

\vskip 10pt \noindent\dem\quad Suppose $N = \lfloor\sqrt{p_n}\rfloor$. If $N = 1$ we have from the proof of
the theorem $h_2 = 2 = N+1$ and $h_3 = 1 = N$. When $N = 2$, $h_4 = 3 = N+1$ and $h_5 = 2 = N$. Conversely,
if $h_n = N+1$ and $h_{n+1} = N$, then $h_n = h_{n+1} + 1$, which implies $d_n = 2N$ according to (6.2).
\hspace{\fill}$\Box$

\vskip 6pt The $n$th prime gap $d_n = (2N + 1 - h_n) + h_{n+1}$ is maximum when both $2N + 1 - h_n$ and
$h_{n+1}$ are maximum. That is, when $h_n$ meets its minimum value and $h_{n+1}$ does so with its maximum,
in the range between $1$ and $2N+1$.

\vskip 10pt \noindent{\thm\ }
{\sl Consider two prime numbers $p_{n+1}$ and $p_n$ on both sides of $(N+1)^2$,
for a given integer $N\geq 2$.  Then, $\displaystyle{h_{n+1} \leq N < \sqrt{p_n} < N+1 \leq h_n}$. }

\vskip 10pt \noindent\dem\quad Again set $N = \lfloor\sqrt{p_n}\rfloor$.
Since $h_n - (h_{n+1} + 1)$ is non negative, the identity $d_n = 2N - (h_n - h_{n+1} - 1)$ shows that
$d_n$ meets its maximum value just as $h_n - h_{n+1} - 1$ does so with its minimum value, in the range
$1\leq h_{n+1} \leq 2N+1$.
This maximum value of $d_n$ occurs when $p_n$ and $p_{n+1}$ are further away from each other,
 and on both sides of $(N+1)^2$. In other words, in accordance with Corollary 6.4,
 $d_n$ reaches $2N$ concurrently with $h_n$ and $h_{n+1}$ hitting their minimum and maximum
 values $N+1$ and $N$ respectively. Thus, $h_n\geq N+1$ and $h_{n+1}\leq N$.
\hspace{\fill}$\Box$

\vskip 10pt Instant effects of $h_{n+1} \leq N$ and $N+1 \leq h_n$ are $p_{n+1} < (N+3/2)^2$ and
$p_n > (N+1/2)^2$.

\vskip 10pt\noindent{\thm\ }{\sl Let the primes $p_{n+1}$ and $p_n$ be as above. Then,
$\mu_n$ is bounded below by $1/2$ at the same time that $\mu_{n+1}$ is less than $1/2$.}

\vskip 10pt\noindent Theorem 6.5 implies that $h_n = 1 + 2N - d_n + h_{n+1}\geq N + 1$. i.e. $d_n \leq N + h_{n+1}$.
This, together with (6.4) yield the coming result.

\vskip 10pt\noindent{\cor\ }{\sl $\displaystyle{2 + h_{n+1} \leq d_n \leq N + h_{n+1} < \sqrt{p_n} + h_{n+1} \ (n\geq 4)}$,
under the hypothesis of the theorem.}

\vskip 10 pt\noindent As seen in the proof of Theorem 6.3, $n=2$ and $n = 4$ meet the upper bound $N+h_{n+1}$.
  We did have for $N=1$, $d_2 = 2 = N + h_{n+1}$, and for $N=2$,
$d_4 = 4 = N + h_{n+1}$. Incidentally, on other values of $n$ as well, for instance, $n=6 \ (N=3)$ and $n=11 \ (N=5)$,
the bound $N + h_{n+1}$  is reached. Nevertheless, for these values, $d_n$ is less than $2N$. On the other hand,
there is a relationship between a lower bound on $h_n$ (i.e. on $\mu_n$) and the upper bound $N+h_{n+1}$.
Obviously, if for example, $\mu_n$ is greater than $3/4$ (i.e. $h_n \geq 3N/2 + 1$),
then $d_n$ is bounded above by $N/2 + h_{n+1}$.

\vskip 10pt Finally, $d_n^2 > p_{n+1}$ induces
$\displaystyle{\Delta_n > \frac{\sqrt{p_{n+1}}}{\sqrt{p_{n+1}} + \sqrt{p_n}} = 1 - \frac{\sqrt{p_n}}{\sqrt{p_{n+1}} + \sqrt{p_n}}
> 1 - \frac{\sqrt{p_n}}{2\sqrt{p_n}} = \frac 12}$. However, there are other integers like $n = 2, 6 \mbox{ and } 11$,
satisfying $d_n^2 < p_{n+1}$ and $\Delta_n > 1/2$ at a time. Unlike the first ones, the corresponding $p_n$ have no even
multiple between $(p_n - d_n)p_{n+1}$ and $p_n^2$. In all cases, these six numbers are such that $p_n$ and
$p_{n+1}$ are on both sides of $(\lfloor{\sqrt{p_n}}\rfloor + 1)^2$, in addition to satisfying $d_n > \sqrt{p_n}$.
And above all we do not know
if $n =2, 4, 6, 9, 9, 11 \mbox{ and } 30$ are the only integers
with $\Delta_n > 1/2$. If so, one can deduce $d_n < \sqrt{p_n}$ for $n\notin\{2, 4, 6, 9, 11, 30\}$.

\subsection{Oppermann's and Brocard's Conjectures}
A first consequence of Theorem 6.6 sets to $2$ the lower bound on the number of
primes between the squares of two consecutive integers.

\noindent{\thm\ }{\sl There are at least two prime numbers between successive perfect squares.}

\vskip 10 pt \noindent\dem\quad The theorem is certainly true for $N^2 = 1, 4$ or $9$.
Suppose that for a fixed integer $N\geq 3$ there exits a unique prime $p_m$ $(m\geq 5)$
between $N^2$ and $(N+1)^2$. Then, calling in Theorem 6.6, the fractional part of $\sqrt{p_m}$
must exceed $1/2$ as $p_m$ is the largest prime less than $(N+1)^2$.
It must also be less than $1/2$ as the very first prime after $N^2$. These contradictory facts show that
another prime different from $p_m$ must exist between $N^2$ and $(N+1)^2$.
\hspace{\fill}$\Box$

\vskip 7pt A straightforward implication of this theorem is the following.

\vskip 10pt \noindent{\cor\ }{\sl $\pi(N^2)\geq 2(N-1)$ \ $(N\geq 1)$ \  and \
$\pi(p_n^2)\geq 2(p_n-1)$ \ $(n\geq 1)$.}

\vskip 10pt\noindent\dem\quad Suppose $N\geq 2$.
By summation of $\displaystyle{\pi\left((j+1)^2\right) - \pi\left(j^2\right) \geq 2}$,
for $j = 1, 2, \cdots , N$ we obtain $\pi((N+1)^2) \geq 2N$, the expected estimate for $N+1$.
Finally, $N=1$ satisfies the corollary.
\hspace{\fill}$\Box$

\vskip 7pt Another consequence is the existence of a prime between $(N-1/2)^2$ and $N^2$,
like between $N^2$ and $(N+1/2)^2$, provided $N\geq 2$; that is to say, a conjecture by Oppermann.

\vskip 10pt \noindent {\thm\ }{\sl $\pi(N^2 - N) < \pi(N^2) < \pi(N^2 + N)$ \ $(N\geq 2)$. }

\vskip 10pt\noindent The number $N$ in Oppermann's conjecture must be an integer.
Really, by setting $N = \sqrt 7 > 2$, the rightmost inequality fails short;
it does so when $N = \sqrt 8$. Therefore, the conjecture does not induce the existence of a prime number
between $x$ and $x + \sqrt x$ for every $x\geq 4$; take $x = 8$ or $x = 24$, for example.

\vskip 10pt\noindent An effect of the presence of a prime between $n^2$ and $(n+1/2)^2$ is that, for every given
integer $k\geq 1$, there exists a sequence of $k$ consecutive integers not all composite. For instance,
 $k^2+1$, $k^2+2$, $\cdots$, $k^2+k$ are successive and all between $k^2$ and $(k+1/2)^2$ and so one of them
 must be prime.

\vskip 7pt We state a subsequent result which refines Brocard's conjecture.

\vskip 10pt \noindent{\pro\ }
{\sl $\displaystyle{\pi\left(p_{n+1}^2\right) - \pi\left(p_n^2\right) \geq 2d_n \ \ (n\geq 1)}$.
 }

\vskip 10pt\noindent\dem\quad Since
$\displaystyle{p_n^2 < \left(p_n + 1\right)^2 < \left(p_n + 2\right)^2 < \cdots < \left(p_n + d_n-1\right)^2 <
\left(p_n + d_n\right)^2}$, we may insert at least two prime
numbers between two consecutive squares.
\hspace{\fill}$\Box$

\section{Bounds on $d_n$ When $\lfloor{\sqrt{p_n}}\rfloor = \lfloor{\sqrt{p_{n+1}}}\rfloor$}
\setcounter{equation}{0}
\renewcommand{\theequation}{\thesection.\arabic{equation}}
In this section, the prime numbers are on the same side of any perfect square. Therefore, $d_n$
takes on the form $h_{n+1} - h_n$. Unlike the preceding situation, we consider two cases related to the
parity of the common integral part of the square roots of the primes. Considering the upper bound $N + h_{n+1}$
from Corollary 6.7, we are entitled to expect that $d_n$
cannot exceed $N=\lfloor\sqrt{p_n}\rfloor$ when $\lfloor\sqrt{p_n}\rfloor = \lfloor\sqrt{p_{n+1}}\rfloor$.
In this case, the quantity $h_{n+1}$ of the corollary is irrelevant.

\vskip 7pt \noindent {\thm\ }{\sl Consider $n$ and $N$ such that $N^2 < p_n < p_{n+1} < (N+1)^2$. Then,
$d_n \leq N$, if and only if, $\Delta_n = \mu_{n+1} - \mu_n < 1/2$. Therefore, $d_n < \sqrt{p_n}$.
}

\vskip 10pt \noindent\dem\quad
$d_n = h_{n+1} - h_n = (\sqrt{p_{n+1}} + \sqrt{p_n})\Delta_n = (2N + \mu_{n+1} + \mu_n)\Delta_n$, for
$N = \lfloor{\sqrt{p_n}}\rfloor = \lfloor{\sqrt{p_{n+1}}}\rfloor$. So, $h_{n+1} - h_n \leq N$ implies
$\displaystyle{\Delta_n < \frac{N}{2N + \mu_{n+1} + \mu_n} = \frac 12
\left(1 - \frac{\mu_{n+1} + \mu_n}{2N + \mu_{n+1} + \mu_n}\right) < 1/2}$.
Conversely, if $\Delta_n < 1/2$, then $h_{n+1} - h_n < N + (\mu_{n+1} + \mu_n)/2$, i.e. $h_{n+1} - h_n \leq N$.
\hspace{\fill}$\Box$

\subsection{Case of $\lfloor{\sqrt{p_n}}\rfloor$ representing an Even Integer}
Let $N\geq 2$ represent an even integer such that
$N = \lfloor{\sqrt{p_n}}\rfloor = \lfloor{\sqrt{p_{n+1}}}\rfloor$.
Since $N$ is even, the largest value of $h_{n+1}$ is $2N-1$, while $d_n$ is less than or equal to $2N-2$.

\vskip 10pt \noindent {\lem\ }{\sl Let $\lfloor{\sqrt{p_n}}\rfloor = N \geq 2$ be an even integer.
In the configuration of the primes $p_n$
and $p_{n+1}$, $d_n = 2N - 2$, if and only if, $N = 2$, i.e. $n = 3$.
}

\vskip 10pt \noindent\dem\quad
$d_n = 2N-2 < \sqrt{2}\sqrt{p_{n+1}}$ $ = $ $\sqrt{2}(N + \mu_{n+1}) < \sqrt{2}(N + 1)$.
 It implies $N < (2 + \sqrt{2})^2/2 < 6$. Then, only for $N = 2$ we have $d_3 = 2N-2$.
\hspace{\fill}$\Box$

\vskip 10pt \noindent {\thm\ }{\sl Consider $p_n$ and $p_{n+1} \ (n\geq 2)$ such that
$\lfloor{\sqrt{p_n}}\rfloor = \lfloor{\sqrt{p_{n+1}}}\rfloor = N$ is an even integer.
Then, $d_n = N$, if and only if, $h_{n+1} = 2N-1$ and $h_n = N-1$.
}

\vskip 10pt \noindent\dem\quad From the proof of Lemma 7.2,
any value of $n$ on which $d_n$ reaches its maximum  is less than or equal to $10$ i.e. $N <6$.
When $N = 2$, $d_3 = 2 = N$
with $h_3 = 1 = N-1$ and $h_4 = 3 = 2N-1$. So, $n = 3$ satisfies the theorem. As to $N=4$,
 second even value of $N<6$, one gets $d_8 = 4 = N$. Thus, $h_9 = 7 = 2N-1$ and $h_8 = 3 = N-1$.
This value of $N$ does not satisfy the lemma. However, in both cases ($N=2$ and $N=4$),
$d_n = N$, $h_{n+1} = 2N-1$ and $h_n = N-1$. Conversely, subtracting $h_n = N-1$ from $h_{n+1} = 2N-1$
 gives $d_n = N$.
\hspace{\fill}$\Box$

\vskip 6pt \noindent {\cor\ }{\sl Under the hypothesis of the lemma, $d_n = \lfloor{\sqrt{p_n}}\rfloor$,
if and only if, $n=3$ or $n=8$. In particular, $d_n = h_{n+1} - h_n \leq N < \sqrt{p_n}$.
}

\vskip 6pt  In light of this, if $d_n >\sqrt{p_n}$ then
$\lfloor{\sqrt{p_n}}\rfloor + 1 = \lfloor{\sqrt{p_{n+1}}}\rfloor$. In addition,
the above corollary leads to an estimate of the difference $\mu_{n+1} - \mu_n$, where
the involved primes are subject to the conditions of the preceding theorem.

\vskip 10pt \noindent {\cor\ }{\sl $\displaystyle{0 < \mu_{n+1} - \mu_n \ < \ \frac 12 -  \frac{\Delta_n^2}{2\sqrt{p_n}}}$
}

\vskip 6pt \noindent\dem\quad We rewrite (6.1) as $\displaystyle{\mu_{n+1} - \mu_n = \left(\frac12 -
\frac{\Delta_n^2}{2\sqrt{p_n}}\right) - \left(\frac12 - \frac{d_n}{2\sqrt{p_n}}\right)}$. Since the last corollary induces
$\frac 12 - \frac{d_n}{2\sqrt{p_n}} > 0$, this establishes the estimate.
\hspace{\fill}$\Box$

\vskip 10pt Further consequence of the theorem, based on the fact that $d_n$ is even whereas $\lfloor{\sqrt{p_n}}\rfloor$
and $h_n$ are of opposite parities, is next.

\vskip 6pt \noindent {\cor\ }{\sl Consider $N$, $p_n$ and $p_{n+1}$ as in the theorem such that
$h_{n+1} < N$, i.e. $p_{n+1}$ is less than $(N+1/2)^{1/2}$. Then, $d_n < \sqrt{p_n} - 1 - h_n$.
On the contrary, if $h_n > \sqrt{p_n}$ then $d_n < h_{n+1} - \sqrt{p_n}$.
}

\vskip 10pt \noindent {\thm\ }{\sl Let $N = \lfloor{\sqrt{p_n}}\rfloor$ be an even integer.
Suppose in addition that $n$ does not assume the value $8$ (i.e. $p_n\neq 19)$. Then
\begin{eqnarray}
d_n \quad \leq \quad N - 2 \ < \ \left\{ \begin{array}{ll} \sqrt{p_{n-1}} - 2, &
\mbox{ if \ $\lfloor{\sqrt{p_{n-1}}}\rfloor = \lfloor{\sqrt{p_n}}\rfloor$}   \\
\sqrt{p_{n-1}} - 1, & \mbox{ otherwise.} \end{array} \right.
\end{eqnarray}
}

\vskip 10pt \noindent\dem\quad Since $\lfloor{\sqrt{p_n}}\rfloor = N \geq 4$ is even, and $n\neq 8$, Corollary 7.4 implies
$d_n \leq N-2 < \sqrt{p_n}-2$. If $\lfloor{\sqrt{p_{n-1}}}\rfloor = N$, then $N-2 < \sqrt{p_{n-1}}-2$. While, $N-2 < \sqrt{p_{n-1}}-1$
when $p_n$ is the smallest prime down to $N^2$.
\hspace{\fill}$\Box$

\vskip 10pt \noindent {\thm\ }{\sl Let $\lfloor{\sqrt{p_n}}\rfloor = \lfloor{\sqrt{p_{n+1}}}\rfloor = N \geq 4$
  be an even integer. If $h_n = 1$, then, $\pi(N^2 + N) - \pi(N^2) \geq 2$.
Whereas, $\pi(N^2 + 2N) - \pi(N^2 + N) \geq 2 \quad (N>4)$, if $h_{n+1} = 2N-1$.
}

\vskip 10pt \noindent\dem\quad When $N=4$, we have $h_7 = 1$ so that $\pi(N^2 + N) - \pi(N^2)= 2$. Suppose $N>4$
and $h_n = 1$. Then, $n > 8$, so that $d_n = h_{n+1} - 1$ is less than or equal to $N-2$.
That is, $h_{n+1}\leq N-1<N$ which means, like $p_n$,
$p_{n+1}$ is before $N^2+N-1$. Therefore, $\pi(N^2 + N) \geq 2 + \pi(N^2)$. If $h_{n+1} = 2N-1$ and $N=4$, we know that
$d_n = 4 = N$. Assume $N>4$ and $h_{n+1} = 2N-1$. Here, $d_n = 2N-1 - h_n\leq N-2$. That is, $h_n \geq N+1$,
 locating $p_{n+1}$ between $N^2+N+1$ and $N^2+2N-3$, hence $\pi(N^2 + 2N) \geq 2 + \pi(N^2 + N).$ This completes the proof.
\hspace{\fill}$\Box$

\subsection{Case of $\lfloor{\sqrt{p_n}}\rfloor$ representing an Odd Integer}
By contrast, if $N$ is an odd integer, then $h_n$ and $h_{n+1}$ assume even values and the
maximum for $h_{n+1}$ is $2N-2$ so that $d_n \leq 2N-4$.

\vskip 10pt \noindent {\lem\ }{\sl Consider an odd integer $N\geq 3$ satisfying
$N = \lfloor{\sqrt{p_n}}\rfloor = \lfloor{\sqrt{p_{n+1}}}\rfloor$. Then,
 $d_n = 2N - 4$, if and only if, $N = 3$, i.e. $n = 5$.
}

\vskip 10pt \noindent The proof presents no additional difficulties compared to the preceding.
Here, $h_6 = 4 = 2N-2$, $h_5 = 2 = N-1$ and $d_5 = 2 = 2N-4 = N-1$. However, $d_n = 2N-4 < \sqrt{2p_{n+1}} < \sqrt{2}(N+1)$,
like in the proof of Lemma 7.2, we come up with $N < (4 + \sqrt{2})(2 + \sqrt{2})/2 < 10$.
For $N = 3, 7$ and $9$, we have respectively $n = 5, 16$ and $24$  with $p_{16} = 53$ and $p_{24} = 89$.
Thus, $d_5 = 2 = N-1$, $h_5 = 2$ and
$h_6 = 4$; $d_{16} = 6 = N-1$, $h_{16} = 4$ and $h_{17} = 10$;
finally, $d_{24} = 8 = N-1$, $h_{24} = 8$ and $h_{25} = 16$. In each case, $d_n = N - 1.$

\vskip 10pt \noindent {\thm\ }{\sl Let $\lfloor{\sqrt{p_n}}\rfloor = \lfloor{\sqrt{p_{n+1}}}\rfloor = N$ represent
an odd integer. Then, $d_n = N - 1$, if and only if, $N = 3, 7$ or $9$, i.e. $n = 5, 16$ or $24$. In particular,
$d_n = h_{n+1} - h_n \leq N - 1 < \sqrt{p_n} - 1 < \sqrt{p_{n-1}}$.
}

\vskip 7pt \noindent We derive $\displaystyle{\frac 12 - \frac{d_n}{2\sqrt{p_n}} > \frac 1{2\sqrt{p_n}}}$
from the theorem, and along with (6.1) to state

\vskip 10pt \noindent {\cor\ }{\sl $\displaystyle{0 < \mu_{n+1} - \mu_n < \frac 12 - \frac{\Delta_n^2}{2\sqrt{p_n}}
- \frac 1{2\sqrt{p_n}}}$.
}

\vskip 7pt \noindent {\cor\ }{\sl Let $\lfloor{\sqrt{p_n}}\rfloor = \lfloor{\sqrt{p_{n+1}}}\rfloor$ represent an odd integer.
Then
\begin{eqnarray}
d_n \quad \leq \quad \lfloor{\sqrt{p_n}}\rfloor - 1 \ < \ \left\{ \begin{array}{ll} \sqrt{p_{n-1}} - 1 &
\mbox{ if \ $\lfloor{\sqrt{p_{n-1}}}\rfloor = \lfloor{\sqrt{p_n}}\rfloor$}   \\
\sqrt{p_{n-1}} & \mbox{ otherwise.} \end{array} \right.
\end{eqnarray}
}

\vskip 10pt The following statement is the like of Corollary 7.6 for odd integers $N$.

\vskip 10pt \noindent {\cor\ }{\sl Consider $N$, $p_n$ and $p_{n+1}$ as in the theorem. If
$h_{n+1} < \lfloor{\sqrt{p_n}}\rfloor = N$, then $d_n < \sqrt{p_n} - 1 - h_n$.
In terms of $\lfloor{\sqrt{p_{n-1}}}\rfloor$, one has
$d_n < \sqrt{p_{n-1}} - 1 - h_n$, if $\lfloor{\sqrt{p_{n-1}}}\rfloor = \lfloor{\sqrt{p_n}}\rfloor$;
and $d_n < \sqrt{p_{n-1}} - h_n$ otherwise. When $h_n$ exceeds $\sqrt{p_n}$, one has $d_n < h_{n+1} - \sqrt{p_n}$.
}

\vskip 10pt \noindent\dem\quad
If $h_{n+1} < N$, then $d_n = h_{n+1} - h_n \leq N - 1 - h_n < \sqrt{p_n} - 1 - h_n$. If in addition,
$\lfloor{\sqrt{p_{n-1}}}\rfloor = \lfloor{\sqrt{p_n}}\rfloor$,
then $\lfloor{\sqrt{p_n}}\rfloor - 1 - h_n = \lfloor{\sqrt{p_{n-1}}}\rfloor - 1 - h_n < \sqrt{p_{n-1}} - 1 - h_n$.
On the contrary, $\lfloor{\sqrt{p_{n-1}}}\rfloor = \lfloor{\sqrt{p_n}}\rfloor - 1$ implies $p_n$ and $p_{n-1}$ are
on both sides of $N^2$, and according to Theorem 6.6, $\mu_n < \mu_{n-1}$. That is to say,
$(N + \mu_n) - 1 - h_n < (N - 1 - \mu_{n-1}) - h_n$, i.e. $\sqrt{p_n} - 1 - h_n < \sqrt{p_{n-1}} - h_n$.
The last inequality of the corollary comes from $h_n > \sqrt{p_n}$.
\hspace{\fill}$\Box$

\vskip 10pt \noindent By calling on both Corollaries 7.5 and 7.11, since they induce $d_n < \sqrt{p_n}$,
one strengthens Theorem 2.14. There are at least two other primes $p$ and $q$ in
$p_n < p < q < p_n + 2\sqrt{p_n} + d_n/\sqrt{p_n} = (\sqrt{p_n} + 1)^2 - (1 - d_n/\sqrt{p_n})$, provided $p_n$ and
$p_{n+2}$ are both between two consecutive perfect squares.

\vskip 10pt We conclude this section raising a concern.
It hints at a change in the parity of $\lfloor{D_n}\rfloor$,
when $p_n$ varies  between two consecutive perfect squares, starting as an even integer.

\vskip 10pt \noindent {\thm\ }{\sl Let $N \geq 2$ be an integer, and $p_n$
the smallest prime between $N^2$ and $(N+1)^2$.
Then, $\lfloor{D_n}\rfloor := \lfloor{\sqrt{p_{n+1}} + \sqrt{p_n}}\rfloor$ represents an even integer.
}

\vskip 10pt \noindent\dem\quad Since $p_n$ is the prime right after $N^2$, one has $\mu_n < 1/2$, according to
Theorem 6.6 in which $h_n$ takes on the role of $h_{n+1}$. We next deduce via Theorem 6.8 that
$\lfloor{\sqrt{p_n}}\rfloor = \lfloor{\sqrt{p_{n+1}}}\rfloor$.
Finally, the conclusion comes from Corollary 5.9.
\hspace{\fill}$\Box$

\vskip 10pt \noindent Clearly, as soon as $p_n$ passes through $(N+1/2)^2$ to $(N+1)^2$, $\lfloor{D_n}\rfloor$
represents an odd integer provided $p_{n+1}$ is less than $(N+1)^2$.
However, as regard to the largest prime in-between $N^2$ and $(N+1)^2$, a question remains.
{\sl Are there infinitely many integers $N$ for which the largest prime $p_n$ before $(N+1)^2$ is such
that $\lfloor{D_n}\rfloor$ stands for an even integer?} The answer is probably in the affirmative.

\section{Primes Between Consecutive Powers}
\setcounter{equation}{0}
\renewcommand{\theequation}{\thesection.\arabic{equation}}
To extend Legendre's conjecture beyond perfect squares, Andrica's estimate is once again suited to
prove the presence of primes between successive powers. Of course, the results here are proven conditionally on Conjecture 2.5.
Let start giving on a few values of $k\geq 3$, a lower bound for
\begin{eqnarray}
\pi((n+1)^k)  \ - \ \pi(n^k).
\end{eqnarray}

\vskip 5pt\noindent {\thm\ }{\sl Let $n\geq 1$ and define $f(n) := n(n+1)/2$. Then,
$$
\pi\left(n^3\right) < \pi\left(n^3 + f(n)\right) < \pi\left(n^3 + 2f(n)\right)
< \pi\left(n^3 + 4f(n)\right) < \pi\left(n^3 + 6f(n)\right).
$$ In simpler terms, there are at least four primes between $n^3$ and $(n+1)^3$. }

\vskip 10 pt\noindent\dem\quad The inequalities hold for $n = 1, 2, 3, 4, 5$ and $n=6$. On using
$\sqrt{p_{n+1}} - \sqrt{p_n} < \frac{\sqrt 2}2$, we prove that each of the four intervals
$]n^3 , n^3 + f(n)[$; $]n^3 + f(n) , n^3 + 2f(n)[$; $]n^3 + 2f(n) , n^3 + 4f(n)[$ and $]n^3 + 4f(n) , n^3 + 6f(n)[$
contains a prime for $n\geq 7$. In that respect, define four increasing functions of one real variable $x\geq 1$
\begin{eqnarray*}
u(x) & := &  \sqrt{x^3 + f(x)} - \sqrt{x^3} \ = \ \frac{f(x)}{\sqrt{x^3 + f(x)} + \sqrt{x^3}}, \\
v(x)  & := & \sqrt{x^3 + 2f(x)} - \sqrt{x^3 + f(x)}  \ = \  \frac{f(x)}{\sqrt{x^3 + 2f(x)} + \sqrt{x^3 + f(x)}}, \\
w(x)  & := & \sqrt{x^3 + 4f(x)} - \sqrt{x^3 + 2f(x)}  \ = \  \frac{2f(x)}{\sqrt{x^3 + 4f(x)} + \sqrt{x^3 + 2f(x)}} \quad \mbox{and} \\
t(x)  & := & \sqrt{x^3 + 6f(x)} - \sqrt{x^3 + 4f(x)} \ = \ \frac{2f(x)}{\sqrt{x^3 + 6f(x)} + \sqrt{x^3 + 4f(x)}}.
\end{eqnarray*}
Moreover, for $x\geq 2$, $\displaystyle{v(x) < u(x) < t(x) < w(x)}$ and $v(7) = 0.713\cdots > \frac{\sqrt 2}2$.
Consequently, if for instance, the third interval $n^3 + 2f(n) < m < n^3 + 4f(n)$
contains no primes for a certain integer $n\geq 7$,
then, since $w(x)$ increases, there exits an integer $k$ such that
$w(7) \leq w(n) = \sqrt{n^3 + 4f(n)} - \sqrt{n^3 + 2f(n)} < \sqrt{p_{k+1}} - \sqrt{p_k}$.
That is, $\frac{\sqrt 2}2 < 0.713\cdots = v(7) < w(7) < \sqrt{p_{k+1}} - \sqrt{p_k}$,
yielding the needed contradiction. Therefore, the third interval contains a prime and so do the three others,
for $\frac 1{\sqrt 2} < v(7) < u(7) < t(7)$.
\hspace{\fill}$\Box$

\vskip 6pt \noindent{\cor\ }{\sl $\displaystyle{\pi(p_{n+1}^3) - \pi(p_n^3) \geq 4d_n\ }$
 and \  $\pi(n^3)\geq 4(n-1)$ \ $(n\geq 1)$.
 }

\vskip 7pt We move on to the fourth power.

\vskip 10pt\noindent {\thm\ }{\sl
$\displaystyle{\pi((n + 1)^4) - \pi(n^4) \geq 6 \ \ (n\geq 1)}.$ This comes from $\pi(2^4) = 6$
and for $n\geq 2$, $\displaystyle{\pi(n^4) < }$
$\displaystyle{ \pi(n^4+n^3/2) < \pi(n^4+n^3) < \pi(n^4+2n^3) < \pi(n^4+3n^3) < \pi(n^4+4n^3) < \pi((n + 1)^4)}.$
}

\vskip 10pt\noindent\dem\quad The inequalities are valid for $n = 2$ and $ n = 3$.
When $n\geq 4$, we use the same idea like in last theorem.
Define on $[4 , \infty[$ the next increasing functions: $\displaystyle{u_1(x) = \sqrt{x^4+x^3/2} - x^2 ; \ \ }$
 $\displaystyle{u_2(x) = \sqrt{x^4+x^3} - }$ $\displaystyle{\sqrt{x^4+x^3/2}; \ }$  $\displaystyle{u_3(x) =
 \sqrt{x^4+2x^3} - \sqrt{x^4+x^3} ; \ }$ $\displaystyle{u_4(x) = \sqrt{x^4+3x^3} - \sqrt{x^4+2x^3}; \ \ }$
$\displaystyle{u_5(x) = \sqrt{x^4+4x^3} - }$ \\
 $\displaystyle{\sqrt{x^4+3x^3}; \ \ }$ $\displaystyle{u_6(x) = (x+1)^2 - \sqrt{x^4+4x^3}.}$
On $[4 , \infty[$, they assume values that are greater than or equal to $u_2(4) = 0.917\cdots > \sqrt 2/2$.
The sequel of the proof is the same as  the previous theorem.
\hspace{\fill}$\Box$

\vskip 10pt \noindent{\cor\ }{\sl $\displaystyle{\pi(p_{n+1}^4) - \pi(p_n^4) \geq 6d_n\ }$
 and \ $\pi(n^4)\geq 6(n-1)$  \ $(n\geq 1)$.
}

\vskip 10pt The instance of $k=5$ in (8.1) is next.

\vskip 2pt\noindent {\thm\ }{\sl Let $f(x) := (5x^4 + 10x^3 + 10x^2 + 5x + 1)/11$
and $x\geq 2$. Then,
$\displaystyle{\pi\left(n^5\right) < \pi\left(n^5 + f(n)\right) < }$
$\displaystyle{ \pi\left(n^5 + 2f(n)\right) < \pi\left(n^5 + 3f(n)\right) <
\pi\left(n^5 + 4f(n)\right) < }$
$\displaystyle{\pi\left(n^5 + 5f(n)\right) < \pi\left(n^5 + 6f(n)\right) < }$
$\displaystyle{ \pi\left(n^5 + 7f(n)\right) }$ \\
$\displaystyle{ <  \pi\left(n^5 + 8f(n)\right) < \pi\left(n^5 + 9f(n)\right) <
\pi\left(n^5 + 10f(n)\right) < \pi\left((n+1)^5\right) \ (n\geq 2).
}$
To be more precise, there are at least eleven primes between $n^5$ and $(n+1)^5$. }

\vskip 10 pt\noindent\dem\quad For $n=1$, $\pi(2^5) = 11$ and $ \pi(3^5) - \pi(2^5) = 42$. On
 $[3 , \infty[$, we define eleven functions $u_1, u_2, u_3, \cdots, u_{11}$ of a real variable:
\begin{eqnarray*}
u_1(x) & = & \sqrt{x^5 + f(x)} - \sqrt{x^5} \ = \ \frac{f(x)}{\sqrt{x^5 + f(x)} + \sqrt{x^5}}, \\
u_2(x) & = & \sqrt{x^5 + 2f(x)} - \sqrt{x^5 + f(x)} \ = \ \frac{f(x)}{\sqrt{x^5 + 2f(x)} + \sqrt{x^5 + f(x)}}, \\
u_3(x) & = & \sqrt{x^5 + 3f(x)} - \sqrt{x^5 + 2f(x)} \ = \ \frac{f(x)}{\sqrt{x^5 + 3f(x)} + \sqrt{x^5 + 2f(x)}}, \\
u_4(x) & = & \sqrt{x^5 + 4f(x)} - \sqrt{x^5 + 3f(x)} \ = \ \frac{f(x)}{\sqrt{x^5 + 4f(x)} + \sqrt{x^5 + 3f(x)}}, \\
u_5(x) & = & \sqrt{x^5 + 5f(x)} - \sqrt{x^5 + 4f(x)} \ = \ \frac{f(x)}{\sqrt{x^5 + 5f(x)} + \sqrt{x^5 + 4f(x)}}, \\
u_6(x) & = & \sqrt{x^5 + 6f(x)} - \sqrt{x^5 + 5f(x)} \ = \ \frac{f(x)}{\sqrt{x^5 + 6f(x)} + \sqrt{x^5 + 5f(x)}}, \\
u_7(x) & = & \sqrt{x^5 + 7f(x)} - \sqrt{x^5 + 6f(x)} \ = \ \frac{f(x)}{\sqrt{x^5 + 7f(x)} + \sqrt{x^5 + 6f(x)}}, \\
u_8(x) & = & \sqrt{x^5 + 8f(x)} - \sqrt{x^5 + 7f(x)} \ = \ \frac{f(x)}{\sqrt{x^5 + 8f(x)} + \sqrt{x^5 + 7f(x)}}, \\
u_9(x) & = & \sqrt{x^5 + 9f(x)} - \sqrt{x^5 + 8f(x)} \ = \ \frac{f(x)}{\sqrt{x^5 + 9f(x)} + \sqrt{x^5 + 8f(x)}}, \\
u_{10}(x) & = & \sqrt{x^5 + 10f(x)} - \sqrt{x^5 + 9f(x)} \ = \ \frac{f(x)}{\sqrt{x^5 + 10f(x)} + \sqrt{x^5 + 9f(x)}}
\quad \mbox{and} \\
u_{11}(x) & = & \sqrt{(x+1)^5} - \sqrt{x^5 + 10f(x)} \ = \ \frac{f(x)}{\sqrt{(x+1)^5} + \sqrt{x^5 + 10f(x)}}.
\end{eqnarray*}
Each of them is increasing. In fact, in their rational form, the numerator is a polynomial of degree four with
positive leading coefficient $5/\pi(2^5)$, dominating the denominator that is positive on $[2 , \infty[$.
On the other hand, the functions $u_j$ share the same numerator, as a result, $u_j(x) > u_{j+1}(x)$
with $u_j(x) > u_{11}(3) = 1.12\cdots > \sqrt{2}/2$, for every $x\geq 3$ $(1\leq j \leq 10)$.
Each function $u_j$ is used to prove the presence of a prime within the $j$th interval.
All this leads to the conclusion.
\hspace{\fill}$\Box$

\vskip 7pt \noindent{\cor\ }{\sl $\displaystyle{\pi(p_{n+1}^5) - \pi(p_n^5) \geq 11d_n\ }$
and \ $\pi(n^5)\geq 11(n-1)$  \ $(n\geq 1)$.
 }

\vskip 10pt A generalization to all values of $k > 5$: $\displaystyle{\pi((n + 1)^k) - \pi(n^k) \geq \pi(2^k)}$,
requires $x_0\geq 2$, so as to define the functions $u_j$ on $[x_0 , \infty[$.
Then, the increasing nature of $u_j$ is combined with any bound $b$ ($b = \sqrt{2}/2; 7/10 \mbox{ or }
\Delta_4:= \sqrt{11} - \sqrt{7} = 0.6708\cdots$) to prove the existence of
a prime into the $j$th interval. Therefore, the value $\pi(2^k)$ of $j$ can be used to determine $x_0$.
Indeed, for $j = \pi(2^k)$ one has $\displaystyle{u_{\pi(2^k)}(x) = \sqrt{(x+1)^k} - \sqrt{x^5 + (\pi(2^k)-1)f(x)}}$,
and we can choose $x_0$ to equal the smallest integer exceeding $1$ such
that $\displaystyle{u_{\pi(2^k)}(x) > b}$.

\vskip 7pt \noindent{\pro\ }{\sl $\displaystyle{\pi(2^{k+1}) \geq 2 + \pi(2^k) \quad (k\geq 2)}.$ Then,
$\displaystyle{\pi(2^k) \geq 2(k-1)}$.
 }

\vskip 10pt\noindent\dem\quad Consider $k\geq 4$, then Corollary 2.7
yields $\displaystyle{\pi(2^k) < \pi(2^k + \sqrt{2^{k+1}})}$, and this holds for $k = 2$ or $3$ as well.
When $k=2$ the quantity $\displaystyle{2^k + \sqrt{2^{k+1}}}$ is less than $7$,
while for $k\geq 3$ it is after $8$. Thus, again one calls on the corollary:
$\displaystyle{\pi(2^k) < \pi(2^k + \sqrt{2^{k+1}}) < \pi(2^k + \sqrt{2^{k+1}} + \sqrt{2(2^k + \sqrt{2^{k+1}})})}$.
To complete the argument, we must compare $2^{k+1}$ to
$\displaystyle{2^k + \sqrt{2^{k+1}} + \sqrt{2(2^k + \sqrt{2^{k+1}})}}$.
That comes to compare $\displaystyle{\sqrt{2 + 2^{(3-k)/2}}}$ and $\displaystyle{2^{k/2}-2}$.
Yet, this last quantity exceeds $\displaystyle{\sqrt{2 + 2^{(3-k)/2}}}$ for $k\geq 4$.
Consequently, for $k\geq 4$, one has $\displaystyle{\pi(2^k) < \pi(2^k + \sqrt{2^{k+1}}) <
\pi(2^k + \sqrt{2^{k+1}} + \sqrt{2(2^k + \sqrt{2^{k+1}})}) \leq \pi(2^{k+1})}$.
\hspace{\fill}$\Box$

\vskip 10pt The estimate $\displaystyle{\pi((n + 1)^k) - \pi(n^k) \geq 2(k - 1) \ (k\geq 2 , n\geq 1)}$
is certainly true, but we do better with the next result.

\vskip 10pt \noindent{\thm\ }{\sl Given $k\geq 2$, there are $\pi(2^k)$ primes between two consecutive $k$th powers.
That is, $\displaystyle{\pi((n + 1)^k) - \pi(n^k) \geq \pi(2^k)}$; hence
$\displaystyle{\pi(p_{n+1}^k) - \pi(p_n^k) \geq \pi(2^k)d_n\ }$ and \ $\pi(n^k)\geq \pi(2^k)(n-1)$ \ $(n\geq 1)$.
Recall that  $\displaystyle{\pi(2^k) = 2^{k-1} + 1 - \phi_c(2^k)}$, where $\displaystyle{\phi_c(n)}$
counts the number of composite integers $m < n$ and coprime to $n$.
 }

\vskip 10pt\noindent\dem\quad The theorem holds for $k\leq 5$, as seen through preceding results. Suppose then, $k\geq 6$
and consider for $x\geq 2$ the polynomial $f_k$ of degree $k-1$ defined by
$\displaystyle{f_k(x) := \frac{(x+1)^k - x^k}{\pi(2^k)} = \frac 1{\pi(2^k)}\sum_{r=0}^{k-1}{k\choose r}x^r}$.
On the one hand, since $x$ is greater than or equal to $2$ and $k$ is after $5$, $x^k$ is definitely greater than $8$.
So we refer to Corollary 2.7 to write: $\displaystyle{\pi\left(x^k + jf_k(x)\right) <
\pi\left(x^k + jf_k(x) + \sqrt{2(x^k + jf_k(x))}\right)}$,\quad $j = 0, 1, 2, \cdots, \pi(2^k)-1$.
On the other hand, the identity $\displaystyle{2(x^k + jf_k(x))} = 2x^k(1 + jf_k(x)/x^k)$ implies
that $f_k(x)$ as polynomial of degree $k-1$ dominates
$\displaystyle{\sqrt 2x^{k/2}\sqrt{1 + jf_k(x)/x^k}}$. Consequently, since the range
$\displaystyle{x^k + jf_k(x) < m < x^k + jf_k(x) + \sqrt{2\left(x^k + jf_k(x)\right)}}$
shelters a prime, so does the wider interval $\displaystyle{x^k + jf_k(x) < m < x^k + (j+1)f_k(x)}$. Thus,
$\displaystyle{\pi\left(x^k + jf_k(x)\right) < \pi\left(x^k + (j+1)f_k(x)\right)}$, \ for $k\geq 2$, $n\geq 1$ and
$j = 0, 1, \cdots, \pi(2^k)-1$,  where for $j = \pi(2^k)-1$, one has $x^k + (j+1)f_k(x) = (x+1)^k$.
\hspace{\fill}$\Box$

\section{Twin Primes and $\sqrt{p_{n+1}}\Delta_n$}
\setcounter{equation}{0}
\renewcommand{\theequation}{\thesection.\arabic{equation}}
Recall that $\displaystyle{\Delta_n := \sqrt{p_{n+1}} - \sqrt{p_n} < 1}$ and
$\displaystyle{D_n := \sqrt{p_{n+1}} + \sqrt{p_n}}$.

\noindent
{\thm\ }{\sl Given an integer $n\geq 2$, the following statements are equivalent. \\
1. \ $d_n = 2\qquad $
2. \ $\displaystyle{\sqrt{p_n}\Delta_n = \{\sqrt{p_n}\Delta_n\}}$ \qquad
3. \ $\displaystyle{\sqrt{p_{n+1}}\Delta_n = 1 + \{\sqrt{p_{n+1}}\Delta_n\}}$ \quad
4. \ $\displaystyle{\Delta_n^2 + 2\sqrt{p_n}\Delta_n = 2}$.
}

\vskip 10pt \noindent\dem\quad
The identity $\displaystyle{d_n/2 - 1 = \lfloor{\sqrt{p_n}\Delta_n}\rfloor = \sqrt{p_n}\Delta_n - \{\sqrt{p_n}\Delta_n\}}$
establishes the equivalence between the first two statements of the theorem. The second assertion implies that \\
$0 = \lfloor{\sqrt{p_n}\Delta_n}\rfloor =  \lfloor{\sqrt{p_{n+1}}\Delta_n}\rfloor -1 = - 1 +
\sqrt{p_{n+1}}\Delta_n - \{\sqrt{p_{n+1}}\Delta_n\}$, namely the third one. We also have
$\sqrt{p_{n+1}}\Delta_n - \{\sqrt{p_{n+1}}\Delta_n\} = 1 = \lfloor{\sqrt{p_n}\Delta_n}\rfloor + 1 = \sqrt{p_n}\Delta_n + 1 -
\{\sqrt{p_n}\Delta_n\} = \sqrt{p_n}\Delta_n + $ $ \{\sqrt{p_{n+1}}\Delta_n\}$,
giving $1 = \sqrt{p_n}\Delta_n + \{\sqrt{p_{n+1}}\Delta_n\} = \sqrt{p_n}\Delta_n + \Delta_n^2/2$, i.e. the fourth proposition.
From this proposition comes $2(1 - \sqrt{p_n}\Delta_n) = \Delta_n^2$ which means $\sqrt{p_n}\Delta_n < 1$.
This goes back to $\lfloor{\sqrt{p_n}\Delta_n}\rfloor = d_n/2 - 1 = 0$, say the very first statement. This completes
the equivalence loop, hence the proof.
\hspace{\fill}$\Box$

\vskip 7pt Similarly, one may derive merely out of the fourth assertion that
$\displaystyle{1/\Delta_n \ - \ 1/D_n = \sqrt{p_n}}$, that is
$\displaystyle{D_n = \frac{\Delta_n}{1 - \sqrt{p_n}\Delta_n}}$, if and only if $d_n = 2$.
And from the third proposition, that $d_n = 2$ is equivalent to
$\displaystyle{\sqrt{p_{n+1}}\Delta_n + \{\sqrt{p_n}\Delta_n\} = 2}$.

\vskip 10pt \noindent {\thm\ }
{\sl Let $2 \leq m_1 < m_2$ be two integers such that $\displaystyle{d_{m_1} = d_{m_2} = 2}$.
Then, we have
$$ 1 \leq \sqrt{\frac{p_{m_2}}{p_{m_1+1}}} < \frac 1{\sqrt{p_{m_1+1}}\Delta_{m_2}} < \frac 2{D_{m_1}\Delta_{m_2}} =
\frac{\Delta_{m_1}}{\Delta_{m_2}} = \frac{D_{m_2}}{D_{m_1}} < \frac 1{\sqrt{p_{m_1}}\Delta_{m_2}} = \frac{D_{m_2}}{2\sqrt{p_{m_1}}}.
$$
  }

\vskip 7pt \noindent\dem\quad Divide the equalities
 $2 = \Delta_{m_1}D_{m_1} = \Delta_{m_2}D_{m_2}$ by $D_{m_1}\Delta_{m_2}$ gives:
 $\displaystyle{\frac 2{D_{m_1}\Delta_{m_2}} = \frac{\Delta_{m_1}}{\Delta_{m_2}} = \frac{D_{m_2}}{D_{m_1}}}$,
 where either fraction certainly equals $D_{m_2}\Delta_{m_1}/2$.
The estimates $\displaystyle{2\sqrt{p_{m_1}} < D_{m_1} < 2\sqrt{p_{m_1+1}}}$, in turn, induce
$\displaystyle{\frac 1{\sqrt{p_{m_1+1}}\Delta_{m_2}} < \frac 2{D_{m_1}\Delta_{m_2}} < \frac 1{\sqrt{p_{m_1}}\Delta_{m_2}}}$;
hence the rightmost inequality. As to the others, recalling
$\displaystyle{\{\sqrt{p_{m_2}}\Delta_{m_2}\} = \sqrt{p_{m_2}}\Delta_{m_2}}$, we have
$\displaystyle{\frac 1{\sqrt{p_{m_1+1}}\Delta_{m_2}} > \frac{\{\sqrt{p_{m_2}}\Delta_{m_2}\}}{\sqrt{p_{m_1+1}}\Delta_{m_2}} =
\frac{\sqrt{p_{m_2}}\Delta_{m_2}}{\sqrt{p_{m_1+1}}\Delta_{m_2}} = \sqrt{\frac{p_{m_2}}{p_{m_1+1}}} \geq 1}$,
namely, the last two leftmost estimates to prove.
\hspace{\fill}$\Box$

\vskip 7pt  Note that $m_1$ and $m_2$ are subject to no other constraint,
apart from $d_{m_1} = d_{m_2} = 2$. So, there may exist one or more integers $n$, in the range between $m_1$ and $m_2$,
such that $d_n = 2$ as well.

\vskip 10pt \noindent {\cor\ }{\sl
$\displaystyle{\sqrt{p_{m_1+1}}\Delta_{m_2} < 1 < \sqrt{p_{m_2+1}}\Delta_{m_2} < \sqrt{p_{m_1+1}}\Delta_{m_1} <
\sqrt{p_{m_2+1}}\Delta_{m_1}}$, where $ m_1$ and $m_2$ are like in the theorem. Furthermore,
$\displaystyle{\sqrt{p_{m_1+1}}\Delta_{m_1} < 5/4}$.
}

\vskip 6pt \noindent\dem\quad The leftmost estimate proceeds from the theorem, while the next ones
are straightforward. Since $\displaystyle{\{\sqrt{p_{m_1+1}}\Delta_{m_1}\}}$ is less than $1/4$,
follows the last statement of the corollary.
\hspace{\fill}$\Box$

\vskip 10pt \noindent {\cor\ }{\sl If $m_1$ and $m_2$ are as above, then,
$\displaystyle{2\Delta_{m_2} \ < \ \Delta_{m_1}\left(1 \ + \ \sqrt{\frac{p_{m_1}}{p_{m_1+1}}} \right)}$.
}

\vskip 10pt \noindent\dem\quad $2 = \Delta_{m_1}D_{m_1}$ implies
$\displaystyle{\frac 2{\sqrt{p_{m_1+1}}\Delta_{m_2}} = \frac{\Delta_{m_1}}{\Delta_{m_2}}\times\frac{D_{m_1}}{\sqrt{p_{m_1+1}}}
= \frac{\Delta_{m_1}}{\Delta_{m_2}}\left(1 + \sqrt{\frac{p_{m_1}}{p_{m_1+1}}}\right)}$ the leftmost side of which, according to
the theorem, is greater than $2$. The conclusion yet appears.
\hspace{\fill}$\Box$

\subsection{Successive Twin Prime Pairs}
We prove that $\Delta_m < \Delta_n$ for every $n < m$ provided $d_m = 2$.

\vskip 7pt \noindent {\thm\ }{\sl Let $m\geq 5$ satisfy $d_m = 2$. Then,
$\displaystyle{\Delta_n - \Delta_m > 0 \quad (1 \leq n < m)}$.
}

\vskip 10pt \noindent\dem\quad
$\displaystyle{\Delta_n = \frac{d_n}{D_n} \leq \Delta_4}$. Also,
$\displaystyle{\Delta_n = \Delta_m + \left(\frac{d_n}{D_n} - \frac{d_m}{D_m}\right).}$
If $n > 1$ is less than $m$ $(m\geq 5)$, then
$\displaystyle{\Delta_4 = \frac{4}{D_4} \geq \Delta_n = \frac{d_n}{D_n} \geq \frac{d_m}{D_n} > \frac{d_m}{D_m} =
\Delta_m}$, for $d_n \geq d_m = 2$. Moreover,
$\displaystyle{\Delta_1 = \frac{1}{D_1} > \frac{2}{D_m} = \Delta_m}$.
Therefore, $\displaystyle{\Delta_n - \Delta_m = \frac{d_n}{D_n} - \frac{d_m}{D_m}}$ is positive in the range
 $(1 \leq n < m)$. The theorem is thus proved.
\hspace{\fill}$\Box$

\vskip 10pt\noindent {\cor\ }
{\sl $\displaystyle{\{\sqrt{p_{n+1}}\Delta_n\} > \{\sqrt{p_{m+1}}\Delta_m\} \ (1\leq n < m)}$,
under the hypothesis of the theorem. And correspondingly,
$\{\sqrt{p_m}\Delta_m\}$ is greater than $\{\sqrt{p_n}\Delta_n\}$. }

\vskip 10pt\noindent\dem\quad The identities $\Delta_n^2/2 = \{\sqrt{p_{n+1}}\Delta_n\}$ and
$\Delta_m^2/2 = \{\sqrt{p_{m+1}}\Delta_m\}$ combined with the preceding theorem meet the estimate of the corollary.
The other statement of the corollary is induced by $\Delta_n^2/2 = 1 - \{\sqrt{p_n}\Delta_n\}$ and
$\Delta_m^2/2 = 1 - \{\sqrt{p_m}\Delta_m\}$ along with the theorem.
\hspace{\fill}$\Box$

\vskip 10pt Consider $m\geq 3$ such that $d_m = 2$. Then, clearly $d_{m+1}\geq 4$,
and since $\Delta_{m+1}< 1$ and $\Delta_m < 1$,
\begin{eqnarray*}
\frac{D_m}{D_{m+1}} & = & 1 - \frac{\sqrt{p_{m+2}} - \sqrt{p_m}}{D_{m+1}}
   =  1 \ - \ \frac{\Delta_{m+1} + \Delta_m}{D_{m+1}}   >  1 - \frac 2{D_{m+1}} \ > \ 1 - \frac 1{\sqrt{p_{m+1}}}.
\end{eqnarray*}
Thus, $\displaystyle{\frac{D_m}{D_{m+1}} > 1 - \frac 1{\sqrt{p_4}} \
= \ 0.622035\cdots > \frac 12 \geq \frac 2{d_{m+1}}}$. Therefore,
$\displaystyle{\frac 2{D_m} < \frac{d_{m+1}}{D_{m+1}}}$ i.e. $\Delta_m < \Delta_{m+1}$.

\vskip 10pt \noindent {\lem\ }
{\sl Given $n_2 > n_1 \geq 4$ satisfying $d_{n_1} \geq 4$ and $d_{n_2} \geq 4$, suppose there exists
a unique integer $m$ $(n_1< m < n_2)$ such that $d_m = 2$. Then, $\Delta_n \geq \Delta_m$ for every
integer $n$ belonging to $[n_1 , n_2]$, with equality holding if and only if $n=m$.
}

\vskip 10pt \noindent\dem\quad If $n$ is in $[n_1 , m]$, then the lemma is induced by Theorem 9.5.
On the other hand, from (3.3), $\displaystyle{\Delta_n = \frac{d_n}{D_n} = \frac{d_n/2 - 1}{\sqrt{p_n}} + \frac{\{\sqrt{p_n}\Delta_n\}}{\sqrt{p_n}}}$, so that $\Delta_n$ is bounded below
by $\displaystyle{\frac{\{\sqrt{p_n}\Delta_n\}}{\sqrt{p_n}}}$. This bound is met if and only if $d_n = 2$.
In consequence, for $n\in[m , n_2]$,  $\displaystyle{\Delta_n = \frac{d_n}{D_n} \geq \frac 4{D_n}
\geq \frac{\{\sqrt{p_{m}}\Delta_{m}\}}{{\sqrt{p_{m}}}} = \Delta_m}$.
\hspace{\fill}$\Box$

\vskip 7pt We use the lemma to deduce a property of $\Delta_n$ on successive pairs of twin prime number.

\vskip 10pt \noindent {\thm\ }
{\sl Let $2\leq m_1 < m_2$ be two integers corresponding to two consecutive twin prime pairs.
That is to say, $d_{m_1} = d_{m_2} = 2$ and $d_n \geq 4$ for every $n\in ]m_1 , m_2[$. Under these
circumstances, $\displaystyle{d_nD_{m_1} > 2D_n}$; that is, $\Delta_{m_1}$
 is less than $\Delta_n$. }

\vskip 10pt \noindent\dem\quad $d_n\geq 4$ if $n\in ]m_1 , m_2[$ means $\Delta_n$ does not hit
the bound $\displaystyle{\frac{\{\sqrt{p_n}\Delta_n\}}{\sqrt{p_n}}}$. In this range with no
twin primes, one has, $\displaystyle{\frac{d_n}{D_n} \geq \frac 4{D_n} = \frac 1{\sqrt{p_n}} +
\frac{\{\sqrt{p_n}\Delta_n\}}{\sqrt{p_n}} > \frac{\{\sqrt{p_{m_1}}\Delta_{m_1}\}}{{\sqrt{p_{m_1}}}} =
2/D_{m_1} > 2/D_n > 2/D_{m_2}}$.
This gives $\displaystyle{\Delta_n > \Delta_{m_1} > \frac{2\Delta_n}{d_n} > \Delta_{m_2},}$
yielding $\Delta_n - \Delta_{m_1} > 0$.
\hspace{\fill}$\Box$

\vskip 10pt
$\displaystyle{\Delta_n > \Delta_{m_1} > 2/D_n = 2\Delta_n/d_n > \Delta_{m_2}}$, throughout $]m_1,m_2[$,
i.e. $\Delta_n - \Delta_{m_2}$ exceeds $\Delta_n - \Delta_{m_1} >0$. So,
$\displaystyle{\Delta_n - \Delta_{m_2} = \frac{d_n}{D_n} - \frac{2}{D_{m_2}} >
\frac{d_n}{D_n} - \frac{2}{D_{m_1}} = \Delta_n - \Delta_{m_1} > 0}$,
then the first consequence

\vskip 2pt \noindent {\cor\ }{\sl $\displaystyle{\Delta_n > \Delta_{m_1} > \Delta_{m_2}}$
 $ \ (m_1 < n < m_2)$, if $n$, $m_1$ and $m_2$ are under
the same restrictions as in the theorem. }

\vskip 10pt \noindent The second consequence of the theorem represents the converse of Theorem 9.5.

\vskip 10pt \noindent {\cor\ }
{\sl Suppose $m\geq 5$ such that $\displaystyle{\Delta_n > \Delta_m}$ for every $n\in [1 , m[$. Then,
 $d_m = 2$.}

\vskip 10pt Note that if $2\leq n < m$ are such that $d_n = d_m$, then $\Delta_n > \Delta_m$, for $D_m > D_n$.
 At present, we prove that if a given even integer greater than $2$ is assumed by $d_n$ for
 an infinity of integers $n$, then the twin prime conjecture is true. It is based on the previous theorem and
$\displaystyle{\lim_{n\rightarrow\infty}\Delta_n = 0}$.

\vskip 2pt \noindent {\thm\ }{\sl Suppose there are only a finite number of twin primes. Then, no other even integer stands
for prime gap for infinitely many primes.
 }

\vskip 10pt \noindent\dem\quad In the context of the hypothesis, let $m_0$ satisfy $d_{m_{0}} = 2$ and such that $d_n \geq 4$ for
every integer $n > m_0$. Suppose there exists an integer $k\geq 2$ such that $d_n = 2k$ for infinitely many
integers $n_i := n_i(k)$ $(i = 1, 2, \cdots)$. In particular,
$\displaystyle{\Delta_{n_1} > \Delta_{n_2} > \cdots > \Delta_{n_i} > \Delta_{n_{i+1}} > \cdots}$;
i.e. $\left(\Delta_{n_i}\right)_{i\geq 1}$ is a decreasing sequence bounded by $0$.
But, as stated in Theorem 9.8, $\Delta_{m_0} < \Delta_{n_i}$ for all $n_i > m_0$;
hence $\displaystyle{\Delta_{m_0}\leq 0 = \lim_{i\rightarrow\infty}\Delta_{n_i}}$.
This is a non sense since $\Delta_{m_0} > 0$. Consequently, only on a finite number of $n$,
$d_n$ takes $2k$ as value.
\hspace{\fill}$\Box$

\vskip 7pt \noindent {\thm\ }{\sl There are infinitely many twin prime numbers.  }

\vskip 10pt \noindent\dem\quad Suppose the twin prime conjecture is false, and let $m_0$ be the largest integer such that
$d_{m_0} = 2$. That is to say $d_n \geq 4$ for each $n> m_0$. Theorem 9.8 says that $\Delta_n > \Delta_{m_0}$
for every integer $n > m_0$. In particular,
$\displaystyle{0 < \Delta_{m_0} \leq\lim_{n\rightarrow\infty}\Delta_n = 0}$ and
$\Delta_{m_0} = 0$ a nonsense. In other words, for any $\epsilon$ in the interval $0 < \epsilon < \Delta_{m_0}$,
there is no integer $n$ satisfying $\Delta_n <\epsilon$. This contradictory fact to
$\displaystyle{\lim_{n\rightarrow\infty}\Delta_n = 0}$, shows that the number of twin primes is infinite.
\hspace{\fill}$\Box$

\subsection{A Postulate Taste Amongst Twin Primes}
{\sl  Is the estimate $2\Delta_{m_2} > \Delta_{m_1}$, where $2 < m_1 < m_2$ and $d_{m_1} = d_{m_2} = 2$,
a necessary and sufficient condition for the
conclusion of Theorem 9.8?} Indeed, one has the equality
\begin{eqnarray}
\frac{d_n}{D_n} - \frac 2{D_{m_2}} = \frac{d_n}{D_n} - \frac 2{D_{m_1}} + \left(\frac 2{D_{m_1}} - \frac 2{D_{m_2}}\right)
\quad (m_1 < n < m_2)
\end{eqnarray}
the left-hand side of which is greater than or equal to $\displaystyle{4/D_n - 2/D_{m_2} > 2/D_{m_2}}$.
Suppose that $2/D_{m_2}$ exceeds $2/D_{m_1} - 2/D_{m_2}$, i.e. $\Delta_{m_2} > \Delta_{m_1} - \Delta_{m_2}$.
Then, the right-hand side of (9.1) satisfies \\
$\displaystyle{\frac{d_n}{D_n} - \frac 2{D_{m_1}} + \left(\frac 2{D_{m_1}} - \frac 2{D_{m_2}}\right) > \frac 2{D_{m_2}} >
\frac 2{D_{m_1}} - \frac 2{D_{m_2}}}$, i.e.
$\displaystyle{\frac{d_n}{D_n} - \frac 2{D_{m_1}} > \frac 4{D_{m_2}} - \frac 2{D_{m_1}} > 0}$. This yields
 $\Delta_n - \Delta_{m_1} > 2\Delta_{m_2} - \Delta_{m_1} > 0$, hence Theorem 9.8.
Conversely, as regards the estimate $d_n/D_n - 2/{D_{m_1}} > 0$ fulfilling the condition:
$\Delta_{m_2}>\Delta_{m_1}-\Delta_{m_2}$, this seems to be a more difficult task.

\vskip 10pt
When $m_1$ and $m_2$ are related to consecutive twin prime pairs, numerical evidence backs
the inequality $\Delta_{m_1}/\Delta_{m_2} < 3/2$, that implies $2\Delta_{m_2} > \Delta_{m_1}$.
Furthermore, the condition $2/D_{m_2} > 1/D_{m_1}$
points to a possible statement of Bertrand's postulate type for pairs of twin primes.
The leftmost and rightmost non trivial inequalities within
$\displaystyle{\frac 23 < \frac{\sqrt 2}2 \stackrel{?}{<} \frac{\Delta_{m_2}}{\Delta_{m_1}} = \frac{D_{m_1}}{D_{m_2}} < 1 <
\frac{\Delta_{m_1}}{\Delta_{m_2}} = \frac{D_{m_2}}{D_{m_1}} \stackrel{?}{<} \sqrt 2 < \frac 32}$ $ \ (2\leq m_1 < m_2)$
need to be proven or disproved.
Obviously, proving one side of the middle inequalities, implies the other one.
Finally, and above all, {\sl does
$\displaystyle{\Delta_{m_2} > \Delta_{m_1} - \Delta_{m_2}}$ hold under the hypothesis of Theorem 9.8?}

\vskip 10pt \noindent {\thm\ } {\sl Let $3\leq m_1 < m_2$ be two integers corresponding to two
pairs of twin primes. That is, $d_{m_1} = d_{m_2} = 2$. Suppose in addition that
$\lfloor{\sqrt{p_{m_1}}}\rfloor = \lfloor{\sqrt{p_{m_2}}}\rfloor$.
Then, $31p_{m_1} > 25p_{m_2}$ i.e. $p_{m_2} - p_{m_1} < p_{m_1}/4 $.
}

\vskip 10pt \noindent The first instance of this situation is met at $m_1 = 26$, $p_{26} = 101$, $m_2 = 28$, $p_{28} = 107$
with $10^2 < p_{26} < p_{28} < 11^2.$ This happens also for $m_1 = 41$, $p_{41} = 179$, $m_2 = 43$, $p_{43} = 191$
and $\lfloor{\sqrt{179}}\rfloor = \lfloor{\sqrt{p_{191}}}\rfloor = 13.$ Similarly, for $m_1 = 57$, $p_{57} = 269$,
$m_2 = 60$, $p_{60} = 281$, with $\lfloor{\sqrt{p_{57}}}\rfloor = \lfloor{\sqrt{p_{60}}}\rfloor = 16.$ The instance of
$p_{26} = 101$ will be used in the proof of the theorem.

\vskip 10pt \noindent\dem\quad Since $\lfloor{\sqrt{p_{m_1}}}\rfloor = \lfloor{\sqrt{p_{m_2}}}\rfloor$ we have
$\mu_{m_1} := \{\sqrt{p_{m_1}}\} < \mu_{m_2} := \{\sqrt{p_{m_2}}\} < 1$ and
$$
\sqrt{\frac{p_{m_1}}{p_{m_2}}} = \frac{\lfloor{\sqrt{p_{m_1}}}\rfloor + \mu_{m_1}}
{\lfloor{\sqrt{p_{m_1}}}\rfloor + \mu_{m_2}} =  1 - \frac{\mu_{m_2} - \mu_{m_1}}{\lfloor{\sqrt{p_{m_1}}}\rfloor + \mu_{m_2}}
= 1 - \frac{\mu_{m_2} - \mu_{m_1}}{\sqrt{p_{m_2}}}.
$$
But,
$\displaystyle{\frac{\mu_{m_2} - \mu_{m_1}}{\sqrt{p_{m_2}}} < \frac 1{\sqrt{p_{m_2}}} <
\frac 1{\lfloor{\sqrt{107}}\rfloor} = \frac 1{10}}$. Thus,
$\displaystyle{\sqrt{\frac{p_{m_1}}{p_{m_2}}} = 1 - \frac{\mu_{m_2} - \mu_{m_1}}{\sqrt{p_{m_2}}} > \frac 9{10}}$;
i.e. $9\sqrt{p_{m_2}} < 10\sqrt{p_{m_1}}$, and the claimed estimate, since $100/81 < 31/25$.
\hspace{\fill}$\Box$

\vskip 10pt The new restriction on $p_{m_1}$ and $p_{m_2}$ is that
$\lfloor{\sqrt{p_{m_2}}}\rfloor = k + \lfloor{\sqrt{p_{m_1}}}\rfloor$. The integer $k$ is undoubtedly
bounded in order to emerge with the like of preceding theorem. Particularly, if any two consecutive pairs of twin
primes should satisfy a postulate-type estimate, an elusive case.

\subsection{Bounds on the Number of Twin Primes up to $p_n$}
We establish an identity involving $\sqrt{p_n}$ and come up with an equality
between $\sqrt{2p_{n+1}}$ and the number of twin prime pairs up to $p_{n+1}$.
For $i\geq 1$, define the real number $\alpha_i$ by $\alpha_i := \sqrt 2/2 - \Delta_i$. Therefore,
$\alpha_n - \alpha_m = \Delta_m - \Delta_n$ ($1\leq n, m $); and if $\Delta_n$ is greater than $\Delta_m$
then $\alpha_m$ exceeds $\alpha_n$, and vice versa. Then, by summation over $i\geq n$ leads to
$\displaystyle{\sqrt{p_{n+1}} = \frac{\sqrt{2}}2n + \sqrt{2} - \sum_{i=1}^n\alpha_i}$ which
furnishes an upper bound for $\sqrt{p_n}$, albeit less accurate:
$\displaystyle{\sqrt{2n - 1} \ \leq \ \sqrt{p_n} \ \leq \ (n + 1)\sqrt{2}/2 \quad (n\geq 1),}$
where equality holds in the right-hand side for $n = 1$.

On the other hand, $\displaystyle{\sqrt{p_{n+1}} = \frac{\sqrt{2}}2(n + 2) - \sum_{i=1}^n\alpha_i}$
 implies that
\begin{eqnarray*}
\frac{\sqrt{2p_{n+1}}}2 =  1 + \frac n2 - \frac{\sqrt 2}2\sum_{i=1}^n\alpha_i = 1 +
\frac 12\sum_{i=1}^n(1-\sqrt 2\alpha_i)  = (n+1) - \frac 12\sum_{i=1}^n(1+\sqrt 2\alpha_i).
\end{eqnarray*}
Set $j_1:=0$ and denote the number of twin prime pairs up to $p_n$ by
$\displaystyle{j_n \ := \ n - 1 - \sum_{\scriptstyle i = 1 \atop d_i \neq 2}^{n-1}1 \quad (n\geq 2)}.$
 So,
\begin{eqnarray*}
\frac{\sqrt{2p_{n+1}}}2 & = & 1 + n - \frac 12\sum_{i=1}^n(1+\sqrt{2}\alpha_i) \quad = \quad 1 + j_{n+1}
\ + \ \sum_{\scriptstyle i = 1 \atop d_i \neq 2}^n1 \quad - \quad \frac 12\sum_{i=1}^n(1+\sqrt{2}\alpha_i)  \\
& = &  1 \ + \ j_{n+1} \ + \ \sum_{\scriptstyle i = 1 \atop d_i \neq 2}^n1 \quad - \quad
\frac 12\sum_{\scriptstyle i = 1 \atop d_i \neq 2}^n(1+\sqrt{2}\alpha_i) \quad - \quad
\frac 12\sum_{\scriptstyle i = 1 \atop d_i = 2}^n(1+\sqrt{2}\alpha_i) \\
& = &  1 \ + \ j_{n+1} \ + \ \frac 12\sum_{\scriptstyle i = 1 \atop d_i \neq 2}^n(1-\sqrt{2}\alpha_i) \quad - \quad
\frac 12\sum_{\scriptstyle i = 1 \atop d_i = 2}^n(1+\sqrt{2}\alpha_i)
\end{eqnarray*}
and then
$\displaystyle{\sqrt{2p_{n+1}}  =  2  +  2j_{n+1}  -  (B_n - A_n) \ }$
where $\displaystyle{A_0 = B_0 := 0}$, and for $n\geq 1$,
$\displaystyle{A_n := \sum_{\scriptstyle i = 1 \atop d_i \neq 2}^n(1-\sqrt{2}\alpha_i)}$ and
$\displaystyle{B_n := \sum_{\scriptstyle i = 1 \atop d_i = 2}^n(1+\sqrt{2}\alpha_i)}$.
Establishing $\displaystyle{B_n > A_n}$ points to the possible lower bounds on $j_{n+1}$:
\begin{eqnarray}
\frac{d_n}2 \ < \  \sqrt{p_{n+1}}\Delta_n \ < \ \frac{\sqrt{2p_{n+1}}}2 \
\stackrel{?}{<} \ j_{n+1} \quad (n\geq 5).
\end{eqnarray}

It is quite possible that $B_n$ exceeds $A_n$ when we consider the wording of some results
in terms of $\alpha_n$ instead of $\Delta_n$. In fact, Theorem 9.8 and Corollary 9.9
have respectively the following equivalent forms.

\vskip 10pt\noindent{\pro\ }
{\sl Let $2\leq m_1 < m_2$ be two integers corresponding to two consecutive twin prime pairs.
That is to say, $d_{m_1} = d_{m_2} = 2$ and $d_n \geq 4$ for every $n\in ]m_1 , m_2[$. Under these
circumstances, $\displaystyle{\alpha_{m_1} \ = \ \alpha_n + (\Delta_n - \Delta_{m_1}) \ =
\ \alpha_n + \left(d_n/D_n - 2/D_{m_1}\right)}$, where $\displaystyle{d_nD_{m_1} > 2D_n}$;
that is, $\alpha_{m_1} > \alpha_n$. }

\vskip 2pt\noindent{\pro\ }{\sl $\displaystyle{\alpha_{m_2} > \alpha_n + \Delta_n\left(1 - \frac 2{d_n}\right) =
\alpha_n + \frac{d_n - 2}{D_n} > \alpha_{m_1} > \alpha_n }$ $ \ (m_1 < n < m_2)$, if $n$, $m_1$ and $m_2$ are under
the same restrictions as in the theorem. }

\vskip 10pt\noindent Implication mainly from the corollary is that $\alpha_{m_1}$ is greater than
every $\alpha_n$ with $d_n\geq 4$. This holds, regardless of the relative location of $n$ with respect to
$m_1$.

\vskip 10pt\noindent
The twin prime conjecture is expressed in (9.2) by the overlapped signs (question mark and
the inequality sign) between either $d_n/2$ and $j_{n+1}$, $\sqrt{p_{n+1}}\Delta_n$ and $j_{n+1}$
or $\sqrt{2p_{n+1}}/2$ and $j_{n+1}$. Is the $n$th prime gap less than the number
of twin primes up to the $(n+1)$th prime, for $n\geq 5$?
Any of the next three inequalities is an issue put forward for $n$ equals $5$ and beyond:
$d_n < 2j_{n+1}$; $\sqrt{p_{n+1}}\Delta_n < j_{n+1}$ and $\sqrt{p_{n+1}} < \sqrt 2j_{n+1}$.

\vskip 10pt We rephrase the relation $\sqrt{2p_n}/2 < j_n$ by means of elementary functions.
Dusart [4] proved that $p_n \geq n(\log n + \log\log n - 1)$ for $n\geq 2$.
This imposes a lower bound on $\sqrt{p_n}$, thus the following concern
$\displaystyle{n\left(\log n + \log\log n - 1\right) \ < \  2j_n^2 \quad (n\geq 3)? }$
An argument in the affirmative for either questioning will be a decisive statement
for lower bounds on $j_n$.

\section{Rational Numbers as Accumulation Points of $\displaystyle{\left(\{\sqrt{p_n}\}\right)_{n\geq 1}}$}
\setcounter{equation}{0}
\renewcommand{\theequation}{\thesection.\arabic{equation}}
\noindent Anyone of the estimates (4.1) and (4.2) implies the next two theorems.
Recall $h_n = p_n - (\lfloor{\sqrt{p_n}}\rfloor)^2$.

\vskip 5pt\noindent {\pro\ }{\sl Given an integer $h\geq 1$, suppose that $h_n = h$
for infinitely many values of $n$. Then, $\displaystyle{\liminf_{n\rightarrow\infty}\mu_n = 0}$
where $\displaystyle{\mu_n = \{\sqrt{p_n}\}}$. Rephrasing it, $\liminf\mu_n$ as $n$ goes to infinity
is zero provided the set of primes $p_{n_i}$ of the form $N_i^2+h$ \ $(i\geq 1)$ is unbounded.
In this case, the sequence $\displaystyle{(\mu_{n_i})_{i\geq 1}}$
is decreasing and $\displaystyle{\lim_{i\rightarrow\infty}\mu_{n_i} = 0}$
}

\vskip 10pt\noindent\dem\quad Suppose there are an infinity of integers
$N_i$ $(i=1, 2, 3, \cdots )$ such that $\displaystyle{p_{n_i} = N_i^2 + h}$. Which is to say,
$\displaystyle{\sqrt{p_{n_i}} = N_i + \mu_{n_i}}$ where
$\displaystyle{\frac{h}{2\sqrt{p_{n_i}}} < \mu_{n_i} < \frac h{2N_i}}$ \ $(i \geq 1)$.
The sequence $\displaystyle{\left(\frac h{N_i}\right)_{i\geq 1}}$ decreases and vanishes, consequently
 inducing the conclusion of the proposition.
\hspace{\fill}$\Box$

\vskip 10pt\noindent The instance of $h = 1$ is about prime numbers of the shape $m^2+1$. For example,
$\{\sqrt{5}\} = 0.23606\cdots$, $\{\sqrt{257}\} = 0.03121\cdots$,
$\{\sqrt{16901}\} = 0.00384\cdots$ and  $\{\sqrt{50177}\} = 0.00223\cdots$.

\vskip 10pt Likewise, suppose $h_n = 2N-1$ for infinitely many even integers $N$. Then
there exists for each such integer $N_i$, an integer $n_i$ such that $h_{n_i} = 2N_i-1$.
By (4.1) this leads to
$\displaystyle{\frac{2N_i-1}{2\sqrt{p_{n_i}}} = 1 - \frac{2\mu_{n_i}+1}{2\sqrt{p_{n_i}}} < }$ \\
$\displaystyle{\quad\mu_{n_i} < \frac{2N_i-1}{2N_i} = 1 - \frac 1{2N_i}}$,
hence $\displaystyle{\lim_{i\rightarrow\infty}\mu_{n_i} = 1}$, i.e.
$\displaystyle{\limsup_{n\rightarrow\infty}\mu_n = 1.}$ The prime numbers, in this case, are from
$N^2+2N-1$. They are $2$ units up to a perfect square. Thus we state the following result.

\vskip 10pt\noindent {\pro\ }{\sl Let the polynomial $N^2+2N-1$ generate an infinite number of primes $p_n$.
Then, $\displaystyle{\limsup_{n\rightarrow\infty}\{\sqrt{p_n}\} = 1}$. Similarly, if there exits
a fixed integer $h\geq 1$
such that $N^2-h$ represents primes for infinitely many values $N_i \ (i\geq 1)$ of $N$, then
$\displaystyle{\limsup_{n\rightarrow\infty}\{\sqrt{p_n}\} = 1}$. The corresponding $p_{n_i} = N_i^2 - h$
are such that $\displaystyle{(\mu_{n_i})_{i\geq 1}}$ increases.
}

\vskip 10pt\noindent The assumed values by $h_n$ and $\mu_n$ are not only interrelated but depend on the location of
$p_n$ in $]N^2 , (N+1)^2[$. Actually, suppose for example that $N$ is even. Surely, the prime $p_n$ lies between
any two successive squares among the following: $N^2$, $(N+1/4)^2$, $(N+1/2)^2$, $(N+3/4)^2$ and $(N+1)^2$. Then, since
$\mu_n$ is irrational, it is in either $0 < \mu_n < 1/4$; $1/4 < \mu_n < 1/2$; $1/2 < \mu_n < 3/4$ or $3/4 < \mu_n < 1$.
There are, correspondingly, as many $h_n-$related sub-intervals:
$[1 , N/2]$, $]N/2 , N[$, $]N , 3N/2]$ and $]3N/2 , 2N-1]$, respectively. It is to be noted that $h_n$
is an odd integer, since $N$ is even. Each related pair of
$h_n-$interval and $\mu_n-$interval is connected with a specific sub-interval of $]N^2 , (N+1)^2[$, holding the prime $p_n$:
$]N^2 , N^2+N/2[$, $]N^2+N/2 , N^2+N[$, $]N^2+N , N^2+3N/2[$ and $]N^2+3N/2 , N^2+2N-1]$. Similarly,
the situation in which $N > 3$ is odd, i.e. $h_n$ assumes even values, is barely different.
The $\mu_n-$related sub-intervals stay the same, while the $h_n-$related sub-intervals become:
$[2 , (N+1)/2[$, $[(N+1)/2 , N-1]$, $[N+1 , (3N-1)/2]$ and $[(3N+1)/2 , 2(N-1)]$.

\vskip 10pt A real number $x$ is an accumulation point (or cluster point) of a set $S$, if every neighborhood of $x$
contains infinitely many elements of $S$. The preceding two propositions present $0$ and $1$ as plausible accumulation
point of $\mu_n$.

\subsection{The case of $1/2$ as Credible Accumulation Point}
Assume an infinity of $\displaystyle{N_i\geq 2}$ $(i=1, 2, \cdots)$ such
that $\displaystyle{N_i - 1 = h_{n_i}}$. Namely, $\displaystyle{p_{n_i} = N_i^2 + N_i - 1}$. So
there exists  $k_i\geq 1$ satisfying $\displaystyle{h_{n_{i+1}} = h_{n_i} + k_i}$, i.e.
$\displaystyle{N_{i+1} = N_i + k_i}$. Yet,
$\displaystyle{h_{n_{i+1}} = 2\mu_{n_{i+1}}(N_i + k_i) + \mu_{n_{i+1}}^2}$ and
$\displaystyle{h_{n_i} = 2\mu_{n_i}N_i + \mu_{n_i}^2}$ give
$\displaystyle{h_{n_{i+1}} - h_{n_i} = }$  $\displaystyle{2N_i(\mu_{n_{i+1}} - \mu_{n_i}) + 2\mu_{n_{i+1}}k_i +
\mu_{n_{i+1}}^2 - \mu_{n_i}^2 = k_i}$. That is to say,
\begin{eqnarray*}
(\mu_{n_{i+1}} - \mu_{n_i})(2N_i + \mu_{n_{i+1}} + \mu_{n_i}) = (1 - 2\mu_{n_{i+1}})k_i \ \ (i \geq 1).
\end{eqnarray*}
Since $\displaystyle{h_{n_{i+1}} = N_{i+1} - 1}$ yields
$\mu_{n_{i+1}} < 1/2$, one has $1 - 2\mu_{n_{i+1}} > 0$, so each side of the displayed identity is positive.
Therefore, after a certain rank $i_0$,
$\displaystyle{\left(\mu_{n_i}\right)_{i\geq i_0}}$ is an increasing sequence bounded above by $1/2$. Then
estimates (4.1) result in $\displaystyle{\frac{N_i -1}{2\sqrt{p_{n_i}}} =
\frac 12 - \frac{1 + \mu_{n_i}}{2\sqrt{p_{n_i}}} < \mu_{n_i} < \frac 12 - \frac 1{2N_i}}$, hence
$\displaystyle{\lim_{i\rightarrow\infty}\mu_{n_i} = 1/2}$.

\vskip 10pt As illustration, $N^2 + N - 1$, for
$N = 3, 6, 9, 16, 20, 21, 114$ and $131$ gives respectively, $p_n = 11$, $41$, $89$, $271$, $419$,
$461$, $13109$ and $17291$; with corresponding $\mu_n = 0.31662\cdots$, $0.403124\cdots$,
$0.43398\cdots$, $0.46207\cdots$, $0.46948\cdots$, $0.47091\cdots$, $0.49445\cdots$ and $0.49524\cdots$.
 They approach $1/2$ from the left.

\vskip 10pt\noindent
In a similar fashion, speculation about $1/2$ can be carried out
through $\displaystyle{N^2 + N + 1}$ producing infinitely many primes. In this case,
$\displaystyle{p_{n_i} = N_i^2 + N_i + 1} > (N_i+1/2)^2$, i.e. $\displaystyle{\mu_{n_i} > 1/2}$. Thus
$1 - 2\mu_{n_i+1} < 0$ and $\displaystyle{\left(\mu_{n_i}\right)_{i\geq 1}}$ is a decreasing sequence bounded below
by $1/2$. Here $N^2+N + 1$, for $N = 3$, $6$, $20$, $21$, $69$, $110$ and $131$, provides
respectively the following primes $p_n = 13, 43, 421, 463, 4831, 12211$ and $17293$.
Meanwhile, the $\mu_n's$ coincide with: $0.60555\cdots$, $0.55743\cdots$,
$0.51828\cdots$, $0.51743\cdots$, $0.50539\cdots$, $0.50339\cdots$ and $0.50285\cdots$, respectively.
These values, in contrast, near $1/2$ from the right.

\subsection{The General Case}
We explore two other examples. Some rational numbers have decimal representation that becomes periodic,
i.e. with repeat digits in decimal form. With
 $N^2+2/3N+1$ and $N^2+9/11N+1$, we present the cases of $1/3 = 0.333\cdots$ and $9/22 = 0.4090909\cdots$.

\vskip 5pt\hskip 20pt
\begin{tabular}{|l|l|l|l|l|l|l|}  \hline
    \multicolumn{7}{|c|}{$N^2+2/3N+1 \qquad (1/3 = 0.333\cdots)$}  \\  \hline
        $N$ & $6$ & $36$ & $90$ & $402$ & $612$  & $\cdots$  \\  \hline
        $p_n$ & $41$ & $1321$ & $8161$ & $161873$ & $374953$  & $\cdots$  \\  \hline
        $\mu_n$ & $0.403\cdots$ & $0.345\cdots$ & $0.338\cdots$ & $0.3344\cdots$
         & $0.33405\cdots$ & $\cdots$  \\   \hline
\end{tabular}

\vskip 15pt\hskip 5pt
\begin{tabular}{|l|l|l|l|l|l|l|}  \hline
    \multicolumn{7}{|c|}{$N^2+9/11N+1 \qquad (9/22 = 0.4090909\cdots)$}  \\  \hline
        $N$ & $11$ & $33$ & $121$ & $451$ & $715$  & $\cdots$  \\  \hline
        $p_n$ & $131$ & $1117$ & $14741$ & $203771$ & $511811$  & $\cdots$  \\  \hline
        $\mu_n$ & $0.44552\cdots$ & $0.42154\cdots$ & $0.41251\cdots$ & $0.41003\cdots$
         & $0.40967\cdots$ & $\cdots$  \\  \hline
\end{tabular}

\vskip 10pt
The next example deals with $1/3$ and $3/10$ to illustrate the case of two rational numbers close to each other.
The chart above already shows data for $1/3$. Here is a table for $3/10$ with $N^2+3/5N + 1$.
Of course, the rational numbers being different so are the polynomials.

\vskip 5pt\hskip 5pt
\begin{tabular}{|l|l|l|l|l|l|l|l|}  \hline
    \multicolumn{8}{|c|}{$N^2+3/5N+1 \qquad (3/10 = 0.3)$}  \\  \hline
        $N$ & $10$ & $30$ & $60$ & $100$  & $500$ & $\cdots$  & $\cdots$  \\  \hline
        $p_n$ & $107$ & $919$ & $3637$ & $10061$  & $250301$ & $\cdots$  & $\cdots$  \\  \hline
        $\mu_n$ & $0.3440\cdots$ & $0.3150\cdots$ & $0.3075\cdots$ & $0.3045\cdots$ &
        $0.3009\cdots$  & $\cdots$ & $\cdots$  \\  \hline
\end{tabular}

\vskip 15pt
Suppose the existence of an infinity of primes from a specific polynomial.
Given integers $a$ and $b$ $(0 < a < b)$ that are relatively prime, let
$\displaystyle{N^2+\frac{2a}bN \pm 1}$ generate an infinite number of primes $p_{n_i}$ $(i=1, 2, \cdots)$. That is,
$\displaystyle{p_{n_i} = N_i^2 + \frac{2a}bN_i\pm 1}$ where
$\displaystyle{N_i^2 + \frac{2a}bN_i- 1 < (N_i+a/b)^2}$ when the primes $p_{n_i}$
are generated by $\displaystyle{p_{n_i} = N_i^2 + \frac{2a}bN_i - 1}$, and  $\mu_{n_i} < a/b$.
Whereas $\displaystyle{N_i^2 + \frac{2a}bN_i+ 1 > (N_i+a/b)^2}$ for
 $\displaystyle{p_{n_i} = N_i^2 + \frac{2a}bN_i + 1}$; hence $\mu_{n_i} > a/b$. Clearly,
$\displaystyle{N_i = \sqrt{p_{n_i}} - \mu_{n_i} > \sqrt{p_{n_i}} - 1}$
implies $\displaystyle{\frac{2a}bN_i\pm 1 > \frac{2a}b(\sqrt{p_{n_i}} - 1) \pm 1 }$.
Thus, from $\displaystyle{\frac ab - \frac a{b\sqrt{p_{n_i}}}
\pm \frac 1{2\sqrt{p_{n_i}}} < \frac{2aN_i \pm b}{2b\sqrt{p_{n_i}}} <
\mu_{n_i} < \frac{2aN_i \pm b}{2bN_i} = \frac ab \pm \frac 1{2N_i}}$, stems
$\displaystyle{\lim_{i\rightarrow\infty}\mu_{n_i} = a/b}$.

\vskip 7pt \noindent{\cnj\ }{\sl Any rational number $r$ $(0\leq r \leq 1)$ is an
accumulation point for the set of fractional parts of the square roots of primes. In other words,
given a rational number $r$ $(0\leq r < 1)$, the quadratic $N^2 +2rN + 1$ produces an unbounded set
$\displaystyle{S_{r_n}}$ of primes $\displaystyle{p_{r_n}}$
satisfying $\displaystyle{\lim_{n\rightarrow\infty}\{\sqrt{p_{r_n}}\} = r}$. So does
$N^2 +2rN - 1$ $(0 < r \leq 1)$.
}

\vskip 10pt\noindent The conjecture embodies Propositions 10.1 and 10.2 as instantiations.
Moreover, suppose by somehow one knows that none of the polynomials $N^2 + 2r N \pm 1$
produces infinitely many primes. Then, replace $1$ by a prime $q$,
and evaluate the polynomials $N^2 + 2r N \pm q$ on values of $N$ not divisible by $q$.
Even if the latter generate finitely many prime numbers for every fixed $q$,
since the primes $q$ are infinite, the conjecture about rational $r$ as accumulation point still holds.

\vskip 10pt Given two rational numbers $r$ and $s$ $(0\leq r < s\leq 1)$, the polynomials $N^2+2rN+1$ and
$N^2+2sN+1$ do not share a same integer. In particular, there is no prime number generated, at a time, by
both polynomials. In fact, suppose they reproduce a same integer, say, $m\geq 2$. Clearly, this cannot happen
on the same value of $N$, unless $r=s$. Thus, consider two different values $N$ and $M$ $(N < M)$ for which
$m = M^2+2sM+1 = N^2+2rN+1$. There exists, evidently, an integer $k\geq 1$ so that $M=N+k$, which
we transfer to the identities between $m$ and the polynomials. This yields
$2kN+k^2 + 2s(N+k) = 2rN$. Since $2rN$ is less than $2N$, we obtain the contradiction
$2kN+k^2 + 2s(N+k) < 2N$, for $k\geq 1$.

\vskip 20pt
\noindent {\sc Institut de la statistique du Qu\'ebec, Direction des technologies de l'information,}
\noindent {\sc 200 chemin Sainte-Foy, 9e \'etage  Qu\'ebec (Qu\'ebec), G1R 5T4 Canada}

\vskip 10pt
\noindent {\it Email address:} jacques.grah@stat.gouv.qc.ca  \\


\end{document}